\documentclass[11pt,letterpaper]{amsart}

\usepackage[toc,page]{appendix}
\usepackage[T1]{fontenc}
\usepackage[latin1]{inputenc}
\usepackage[frenchb, english]{babel}
\usepackage{amsthm}
\usepackage{amscd}
\usepackage{tensor}
\usepackage{amssymb,amsmath,bbm}
\usepackage{tikz-cd}
\usepackage[all]{xy}
\usepackage{mathrsfs}
\usepackage{hyperref}
\usepackage{verbatim}

\numberwithin{equation}{section}

\theoremstyle{plain}
\newtheorem{theorem}{Theorem}[subsection]

\newtheorem{lemma}[theorem]{Lemma}
\newtheorem{proposition}[theorem]{Proposition}
\newtheorem{prop}[theorem]{Proposition}
\newtheorem{corollary}[theorem]{Corollary}

\newtheorem*{theorem*}{Theorem}

\newcounter{alphalabels}

\newtheorem{theoremx}[alphalabels]{Theorem}

\theoremstyle{definition}
\newtheorem{definition}[theorem]{Definition}
\newtheorem{example}[theorem]{Example}

\theoremstyle{remark}
\newtheorem{remark}[theorem]{Remark}

\usepackage{tikz}
\usetikzlibrary{calc}

\def\Q{{\bb Q}}

\def\Z{{\bb Z}}

\def\zp{{\Z_p}}

\def\zpe{{\mathbb{Z}_p^\times}}
\def\qpe{{\mathbb{Q}_p^\times}}

\def\qp{{\Q_p}}

\def\Hom{\mathrm{Hom}}

\DeclareMathOperator{\Func}{Func}

\def\D{{\n{D}}}

\def\epsilon{\varepsilon}
\def\det{\mathrm{det}}

\def\GL{\mathrm{GL}}

\def\Cond{\mathrm{Cond}}

\def\iHom{\underline{\mathrm{Hom}}}

\def\sol{{{\scalebox{0.6}{$\square$}}}}

\def\Kos{\mathrm{Kos}}

\DeclareMathOperator{\Mod}{Mod}

\DeclareMathOperator{\Spa}{Spa}

\DeclareMathOperator{\Rep}{Rep}

\DeclareMathOperator{\Cont}{{Cont}}

\DeclareMathOperator{\Extdis}{Extdis}

\DeclareMathOperator{\CondAb}{CondAb}
\DeclareMathOperator{\CondSet}{CondSet}
\DeclareMathOperator{\CondRing}{CondRing}

\DeclareMathOperator{\id}{id}

\DeclareMathOperator{\LMod}{LMod}
\DeclareMathOperator{\RMod}{RMod}

\DeclareMathOperator{\Gal}{Gal}
\DeclareMathOperator{\Ind}{Ind}

\DeclareMathOperator{\End}{End}

\DeclareMathOperator{\alg}{alg}
\DeclareMathOperator{\lalg}{lalg}
\DeclareMathOperator{\SpecAn}{AnSpec}
\DeclareMathOperator{\AnSpec}{AnSpec}
\DeclareMathOperator{\iCont}{\underline{Cont}}

\newcommand{\n}[1]{\mathcal{#1}}
\newcommand{\bb}[1]{\mathbb{#1}}
\newcommand{\bbf}[1]{\mathbf{#1}}
\newcommand{\f}[1]{\mathfrak{#1}}
\newcommand{\s}[1]{\mathscr{#1}}

\def\Dist{\mathcal{D}}

\title{Solid locally analytic representations}

\author{Joaqu\'in Rodrigues Jacinto}
\author{Juan Esteban Rodr\'iguez Camargo}

\makeatletter
\def\@tocline#1#2#3#4#5#6#7{\relax
  \ifnum #1>\c@tocdepth 
  \else
    \par \addpenalty\@secpenalty\addvspace{#2}%
    \begingroup \hyphenpenalty\@M
    \@ifempty{#4}{%
      \@tempdima\csname r@tocindent\number#1\endcsname\relax
    }{%
      \@tempdima#4\relax
    }%
    \parindent\z@ \leftskip#3\relax \advance\leftskip\@tempdima\relax
    \rightskip\@pnumwidth plus4em \parfillskip-\@pnumwidth
    #5\leavevmode\hskip-\@tempdima
      \ifcase #1
       \or\or \hskip 1em \or \hskip 2em \else \hskip 3em \fi%
      #6\nobreak\relax
    \dotfill\hbox to\@pnumwidth{\@tocpagenum{#7}}\par
    \nobreak
    \endgroup
  \fi}
\makeatother
\usepackage{hyperref}

\begin{document}

\begin{abstract}We develop the $p$-adic representation theory of $p$-adic Lie groups on solid vector spaces over a complete non-archimedean extension of $\qp$. More precisely, we define and study categories of solid, solid locally analytic and solid smooth representations. We show that the category of solid locally analytic representations of a compact $p$-adic Lie group is equivalent to that of quasi-coherent modules over its algebra of locally analytic distributions, generalizing a classical result of Schneider and Teitelbaum. For arbitrary $G$, we prove an equivalence between solid locally analytic representations and quasi-coherent sheaves over certain locally analytic classifying stack over $G$. We also extend our previous cohomological comparison results from the case of a compact group defined over $\qp$ to the case of an arbitrary group, generalizing results of Lazard and Casselman-Wigner. Finally, we study an application to the locally analytic $p$-adic Langlands correspondence for $\GL_1$.
\end{abstract}

\maketitle
\tableofcontents

\section{Introduction} \label{s:Intro}

Let $p$ be a prime number, $G$ be a $p$-adic Lie group defined over a finite extension $L$ of $\qp$ and let  $\n{K} = (K, K^+)$ be a complete non-archimedean extension of $L$. In this article we give new foundations of the theory of locally analytic representations of $G$ on $\n{K}$-solid vector spaces through the use  of condensed mathematics, generalizing our previous work \cite{SolidLocAnRR}, where the case $G$ compact and $L = \qp$ was studied.  We also obtain new results in the theory of locally analytic representations, such as new comparison theorems of group cohomologies, and a generalization of the classical (anti-)equivalence of Schneider and Teitelbaum between admissible locally analytic representations of $G$ and coadmissible modules over the locally analytic distribution algebra $\n{D}^{la}(G, K)$. We also make the first steps towards geometrizing the theory of locally analytic representations, by establishing an equivalence of cateogories between solid locally analytic representations of $G$ and solid quasi-coherent sheaves over certain locally analytic classifying space of $G$. Finally, as our main application, we state and prove  the locally analytic $p$-adic categorical Langlands correspondence for $\GL_1$. \\ 

\subsection{Motivation} \label{ss:IntroMotiv}

The classical theory of locally analytic representations was developed by Schneider and Teitelbaum (\cite{SchTeitDist}, \cite{SchTeitGl2}) and has had crucial applications, e.g., in the $p$-adic Langlands program \cite{ColmezPhiGamma} and in the study of families of $p$-adic modular forms \cite{EmertonInterpolation}. This extremely rich theory has certain (notably categorical and homological) limitations, coming mostly from the theory of $p$-adic functional analysis. For instance, one needs to put a finiteness condition (admissibility) in order to construct an abelian category of locally analytic representations. In  previous work, we gave new foundations of the theory of locally analytic representations of a compact $p$-adic Lie group $G$ defined over $\Q_p$ through the use of condensed mathematics. We showed that this approach has many advantages and that it puts the theory in a good context to make homological algebra. \\

Recently, in the works \cite{LuePan,LuePanII}, the theory of locally analytic representations has been applied to relate $p$-adic Hodge theory, $p$-adic modular forms, and the theory of $p$-adic differential equations over rigid spaces. The first of these works has been generalized in \cite{camargo2022locally,camargo2023geometric} to arbitrary Shimura varieties, where our  theory of solid locally analytic representations plays a key role. On the other hand, the current conjectural statements of the locally analytic categorical $p$-adic Langlands correspondence \cite{EGH} require the construction of certain derived categories of locally analytic representations. In particular, these works show the need of better categorical foundations of the subject. The main purpose of this work is to provide such foundations for any $p$-adic Lie group.   \\

Our first goal is to define and study enhancements of classical representation categories attached to $p$-adic Lie groups. There are at least three of them, namely continuous, smooth and locally analytic representations. Using the formalism of condensed mathematics, we construct and study the ($\infty$-)categories of solid, solid smooth and solid locally analytic representations of $G$. We denote them, respectively, by $\Rep_{\n{K}_\sol}(G)$, $\Rep^{sm}_{\n{K}_\sol}(G)$, $\Rep^{la}_{\n{K}_\sol}(G)$. These categories arise as the derived category of a corresponding abelian category of representations. Furthermore, these abelian categories contain fully faithfully all the classical categories of continuous, smooth and locally analytic representations on complete  compactly generated locally convex $K$-vector spaces. One of the main advantages of our approach is that many of the difficulties appearing in fundamental constructions in classical representation theory, such as Hochschild-Serre, Shapiros's lemma, duality, etc., are easily overcome with the use of homological algebra when one works in a solid framework. \\

We now explain the main features of the theory, before moving into a detailed account of our results. The first result is an equivalence, for $G$ a compact group, between the (derived) category of solid locally analytic representations of $G$ and the category of solid quasi-coherent sheaves over certain non-commutative adic Stein space associated with $G$. This can be seen as a generalization of a classical anti-equivalence of Schneider and Teitelbaum \cite{SchTeitDist}, which can be recovered from our equivalence when restricting to the (abelian) subcategory of admissible representations after applying a duality functor. This result can also be seen as a step towards geometrizing the category of solid locally analytic representations. Our second  result has also a geometric flavor; we prove that the category of solid locally analytic representations of an arbitrary $p$-adic Lie group $G$ is equivalent to the category of solid quasi-coherent sheaves of certain classifying stack $[* / G^{la}]$. Our third main result is an extension of the cohomological comparison theorems for  solid representations from the case where $G$ is compact and defined over $\qp$ obtained in \cite{SolidLocAnRR} to the general case, extending also the non compact version \cite{CasselmanWigner} of Lazard's isomorphisms \cite{Lazard} from the case of finite dimensional representations to arbitrary solid representations. The main novelty of our approach to the comparison results is that we deduce them in a completely formal way from adjunctions between certain functors. Finally, as an application, we state and prove the locally analytic categorical $p$-adic Langlands correspondence for $\GL_1$ confirming the expectations of \cite[\S 7.1]{EGH}.  \\

\subsection{Main results} Let us now explain our results in more detail. Let $G$ be a $p$-adic Lie group over $L$. Let $\n{K}_\sol[G]$ be the Iwasawa algebra of $G$ over $\n{K}_\sol$, i.e. the free $\n{K}_\sol$-vector space generated by $G$. If $\n{K} = (K, \mathcal{O}_K)$ is a finite extension of $\qp$, then $\n{K}_\sol[G]$ is the classical Iwasawa algebra of $G$, i.e. the dual of the space $C(G, K)$ of continuous functions on $G$. Let $\n{D}^{la}(G, K)$ denote the locally analytic distribution algebra of $G$, i.e. the dual of the space $C^{la}(G, K)$ of locally analytic functions on $G$. We denote by $\Mod_{\n{K}_\sol}(\n{D}^{la}(G, K))$ and $\Mod_{\n{K}_\sol}(\n{K}_\sol[G])$ the ($\infty$-)categories of $\n{D}^{la}(G, K)$ and $\n{K}_\sol[G]$-modules on $\n{K}_\sol$-vector spaces, respectively.

\subsubsection{The cateogory of solid locally analytic representations}

Classically, if $V$ is a (say, $\Q_p$-Banach) continuous representation of $G$, locally analytic vectors of $V$ are defined as those $v \in V$ whose orbit map $o_v : G \to C(G, V), g \mapsto g \cdot v$ is locally analytic as a function of the $p$-adic manifold $G$ with values in $V$ (here, $ C(G, V)$ denotes the continuous functions on $G$ with values in $V$).
One observes also that the orbit map is fixed by the diagonal action of $G$ on $C(G, V)$ given by the left regular action on itself and the action on $V$ (i.e. $(g \cdot \varphi)(-) = g \varphi(g^{-1} \cdot -)$); moreover one can show that $V^{la} \cong C^{la}(G, V)^G$ as $G$-representations, where the action of $G$ on the right hand term is given by the right regular action.

Inspired by this definition, we define the functor of locally analytic vectors $(-)^{Rla}:\Mod_{\n{K}_{\sol}}(\n{D}^{la}( G,K)) \to \Mod_{\n{K}_{\sol}}(\n{D}^{la}( G,K))$ as 
\[
V^{Rla}=R\iHom_{\n{D}^{la}(G,K)}(K, C^{la}(G,V)_{\star_{1,3}})
\]
where we see $V^{Rla}$ endowed with the $\star_{2}$-action of $\n{D}^{la}(G, K)$. We say that an object $V\in \Mod_{\n{K}_{\sol}}(\n{D}^{la}(G,K))$ is (derived) locally analytic if the natural arrow $V^{Rla}\to V$ is an equivalence. We define the ($\infty$-)category of locally analytic representations, denoted as $\Rep^{la}_{\n{K}_{\sol}}(G)$,   to be the full subcategory of $\Mod_{\n{K}_{\sol}}(\n{D}^{la}(G,K))$ whose objects are locally analytic representations of $G$. The following result summarizes our construction of the category of solid locally analytic representations and its main properties (c.f. Propositions \ref{PropositionStableColimitsAndTensor}, \ref{PropPropertiesCategoryLocAn}, \ref{PropReplaGrothendieck} and Corollary \ref{CorollaryFiniteCohoDim} (3)).

\begin{theoremx} \label{TheoIntroAA}
The full subcategory $\Rep^{la}_{\n{K}_\sol}(G) \subset \Mod_{\n{K}_\sol}(\n{D}^{la}(G, K))$ of solid locally analytic representations of $G$ on $\n{K}_\sol$-vector spaces is stable under tensor products and colimits, where the inclusion has a right adjoint given by (derived) locally analytic vectors $V \mapsto V^{Rla}$. Moreover, the following properties are satisfied.
\begin{enumerate}
    \item An object $V \in \Mod_{\n{K}_\sol}(\n{D}^{la}(G, K))$ is locally analytic if and only if $H^i(V)$ is (non-derived) locally analytic for every $i \in \Z$. In particular, $\Rep^{la}_{\n{K}_\sol}(G)$ has a natural $t$-structure.
    \item $\Rep^{la}_{\n{K}_\sol}(G)$ is the derived category of its heart.
    \item The functor of locally analytic vectors satisfies the projection formula, namely, for any $V, W \in \Mod_{\n{K}_\sol}(\n{D}^{la}(G, K))$, one has $(V^{Rla} \otimes_{\n{K}_\sol}^L W)^{Rla} = V^{Rla} \otimes_{\n{K}_\sol}^L W^{Rla}$.
\end{enumerate}
\end{theoremx}

\begin{remark} \leavevmode
    \begin{enumerate}
    \item  Let $V$ be a  locally $L$-analytic representation of $G$ on an $LB$ space in the classical sense. Then point (1) implies that $V$ is an object in $\Mod_{\n{K}_\sol}(\n{D}^{la}(G, K))$  that is derived locally analytic. In particular, classical locally analytic representation theory lives naturally in $\Rep^{la}_{\n{K}_\sol}(G)$.
    \item If $G$ is a $p$-adic Lie group over  $\qp$, then $\n{D}^{la}(G, K)$ is an idempotent algebra over the Iwasawa algebra $\n{K}_\sol[G]$, namely $\n{D}^{la}(G, K) \otimes_{\n{K}_\sol[G]}^L \n{D}^{la}(G, K) = \n{D}^{la}(G, K)$. This implies that the  category of $\n{D}^{la}(G, K)$-modules on $\n{K}_\sol$-vector spaces   embeds fully faithfully in the category of $\n{K}_\sol[G]$-modules on $\n{K}_\sol$-vector spaces (i.e. $\n{K}_{\sol}$-linear representations of $G$). In particular, $\Rep^{la}_{\n{K}_\sol}(G)$ is a full subcategory of $\Mod_{\n{K}_\sol}(\n{K}_\sol[G])$ and one can also define the locally analytic vectors of $\n{K}_\sol[G]$-modules as the right adjoint of this inclusion. This definition of locally analytic vectors coincides with the construction of Theorem \ref{TheoIntroAA} when restricted to $\n{D}^{la}(G, K)$-modules inside $\n{K}_\sol[G]$. Nevertheless, when the group is not defined over $\qp$, both constructions of locally analytic vectors differ, c.f. Remark \ref{RemarkLocAnvectorsSolidRep} for a detailed discussion.
    \item Not every solid $\n{D}^{la}(G, K)$-module is a locally analytic representation. The simplest example is $\n{D}^{la}(G, K)$ itself (see Corollary \ref{CoroLocAnDla}).
    \item The category $\Rep^{la}_{\n{K}_\sol}(G)$ is not stable under limits in $\Mod_{\n{K}_{\sol}}(\n{D}^{la}(G, K))$ (as witnessed, e.g., by $\n{D}^{la}(G, K)$). Nevertheless, being a right adjoint, the functor $(-)^{Rla}$ commutes with limits, i.e. $(\lim W_i)^{Rla} = \lim W_i^{Rla}$ (note though that the limit of the left hand side is a limit in the category $\Mod_{\n{K}_{\sol}}(\n{D}^{la}(G, K))$ while that of the right hand side is in $\Rep^{la}_{\n{K}_{\sol}}(G)$).
    \item We also give an analogue of Theorem \ref{TheoIntroAA} for solid smooth representations, c.f. \S \ref{SmoothQuasiCoherent}.
    \item As a corollary of Theorem \ref{TheoIntroAA}, we obtain a description of $\Rep^{la}_{\n{K}_\sol}(G)$ and $\Rep^{sm}_{\n{K}_\sol}(G)$ as quasi-coherent sheaves on the classifying stack $[* / G]$ of $G$, where $G$ is endowed with the sheaf of locally analytic or smooth functions, c.f. Theorem \ref{TheoremLocAnRepStack} and Proposition \ref{PropositionStackSmooth}.
    \end{enumerate}
\end{remark}

\subsubsection{Solid locally analytic representations via algebra}

It is standard in Representation Theory to view representations as modules over group rings. For instance, finite dimensional representations of a finite group $G$ are equivalent to finite dimensional modules over the group ring of $G$. More generally, complete continuous representations of a profinite group $G$ are equivalent to complete modules over the Iwasawa algebra of $G$ (i.e. the profinite completion of the group ring of $G$).
The case of locally analytic representations of a $p$-adic Lie group is subtler. If $G$ is compact,  the distribution algebra $\n{D}^{la}(G, K)$ is a Fr\'echet-Stein algebra in the sense of \cite{SchTeitDist}, and one can define the category of coherent $\n{D}^{la}(G, K)$-modules, which can be thought of as coherent sheaves over certain (non-commutative) Stein space associated with $\n{D}^{la}(G, K)$.
A coadmissible $\n{D}^{la}(G, K)$-module is a $\n{D}^{la}(G, K)$-module obtained as the global sections of a coherent $\n{D}^{la}(G, K)$-module. Moreover, the global section functor induces an equivalence of categories between coherent $\n{D}^{la}(G, K)$-modules and coadmissible $\n{D}^{la}(G, K)$-modules, c.f. Corollary 3.3. of \textit{loc. cit}.
One of the fundamental results of Schneider and Teitelbaum (\cite[and Theorem 6.3]{SchTeitDist}) is that the category of admissible locally analytic representations is anti-equivalent to the category of coherent $\n{D}^{la}(G, K)$-modules. Our next theorem is a vast generalization of this result, and states that Schneider and Teitelbaum's anti-equivalence upgrades to an equivalence of categories between the whole category of solid locally analytic representations of $G$ and \textit{solid quasi-coherent} sheaves of $\n{D}^{la}(G, K)$. \\

More precisely, for $h\in [0,\infty)$ a parameter depending on some choices,  there is a limit sequence of $h$-analytic distribution algebras $\{\n{D}^{h}(G,K)\}_{h\geq 0}$ such that $\n{D}^{la}(G,K)= \varprojlim_{h\to \infty} \n{D}^{h}(G,K)$. For example, if $G= \bb{Z}_p$ is the additive group of $p$-adic integers then, by the Amice transform, $\n{D}^{la}(\bb{Z}_p, K)$ is isomorphic to the global sections of an open unit disc $\mathring{\bb{D}}_{K}$ over $K$, and the algebras $\n{D}^{h}(\bb{Z}_p,K)$ are overconvergent algebras on closed discs of radius $p^{-\frac{p^{-h}}{p-1}}$.  In this way we can think of the sequence $\{\n{D}^{h}(G,K)\}_{h\geq 0}$ as a family of dagger  affinoid algebras defining closed subspaces of a non-commutative Stein space whose global functions are equal to $\n{D}^{la}(G,K)$. We define the category of \textit{solid quasi-coherent} $\n{D}^{la}(G,K)$-modules to be the limit $\infty$-category 
\[
\Mod^{qc}_{\n{K}_{\sol}}(\n{D}^{la}(G,K))= \varprojlim_{h\to \infty} \Mod_{\n{K}_{\sol}}(\n{D}^{h}(G,K))
\]
where the transition maps are given by the $K$-solid base change $\n{D}^{h}(G,K)\otimes_{\n{D}^{h'}(G,K)}^L-$ for $h'>h$. Concretely, an object in $\Mod^{qc}_{\n{K}_{\sol}}(\n{D}^{la}(G,K))$ is a sequence of objects $(\n{F}_h)_{h\geq 0}$ with $\n{F}_{h}\in \Mod_{\n{K}_{\sol}}(\n{D}^{h}(G,K))$, together with natural equivalences $\n{D}^{h}(G,K) \otimes^L_{\n{D}^{h'}(G,K)} \n{F}_{h'} \xrightarrow{\sim} \n{F}_{h}$ for $h'\geq h$, subject to higher coherences.  In the case where $G = \bb{Z}_p$, the category $\Mod^{qc}_{\n{K}_{\sol}}(\n{D}^{la}(G,K))$ is nothing but that of solid quasi-coherent sheaves on $\mathring{\bb{D}}_K$.  Our second main result gives an algebraic description of the full category of solid locally analytic representations of a compact $p$-adic Lie group.\\

\begin{theoremx} [{Theorem \ref{TheoremEquivalenceLocAn1}}] \label{TheoIntroA}
Let $G$ be a compact $p$-adic Lie group defined over $L$. Then there is an equivalence of (stable $\infty$-)categories
\[ \Mod^{qc}_{\n{K}_\sol}(\n{D}^{la}(G, K)) \xrightarrow{\sim} \Rep^{la}_{\n{K}_\sol}(G) \]
\[ (\n{F}_h)_{h} \mapsto j_! \n{F} := (\varprojlim_h \n{F}_h)^{Rla}. \]
\end{theoremx}

\begin{remark} \hfill
\begin{enumerate}
\item The functor $j_!$ giving the equivalence of categories can be thought of as taking cohomology with compact support of quasi-coherent sheaves. Indeed, if $G=\bb{Z}_p$ the functor  $j_!$ is the   cohomology  with compact support on  $\mathring{\bb{D}}_{K}$ of solid quasi-coherent sheaves as defined in \cite[Lecture XII]{CondensesComplex} for complex spaces. 

\item The functor $j_!$ of Theorem \ref{TheoIntroA} does not respect the natural $t$-structures on both sides and hence does not arise from a functor defined at the level of abelian categories. Indeed,  the module $\n{D}^{la}(G, K)$ defines a quasi-coherent sheaf  which is given by $\n{F}=(\n{D}^h(G, K))_{h \geq  0}$ and one has that $ j_! \n{F} = (\n{D}^{la}(G, K))^{Rla} = C^{la}(G,K)\otimes_{\n{K}_\sol} \chi[-d]$ where $d$ is the dimension of the group $G$ and $\chi=\det(\f{g})^{-1}$ denotes the one dimensional representation defined by the determinant of the dual adjoint representation of $G$ on its Lie algebra $\f{g}$, c.f. Corollary \ref{CoroLocAnDla}.

\item We also prove an analogous version of Theorem \ref{TheoIntroA} for solid smooth representations (Proposition \ref{PropositionQuasiCoherentSmooth}), where the category $\Mod^{qc}_{\n{K}_\sol}(\n{D}^{sm}(G, K))$ is defined as $\varprojlim_{H \subset G} \Mod_{\n{K}_{\sol}}(K[G / H])$ for $H$ running through all the open compact subgroups of $G$.
\end{enumerate}
\end{remark}

From Theorem \ref{TheoIntroA}, we can recover Schneider-Teitelbaum's anti-equivalence as follows.

\begin{theoremx} [Proposition \ref{PropCommutativeDiagramDual}]
Let $G$ be a compact $p$-adic Lie group defined over $L$.  There is a locally analytic contragradient functor on $\Rep_{\n{K}_{\sol}}^{la}(G)$ given by
\[
V\mapsto V^{\vee,Rla}= R\iHom_{K}(V,K)^{Rla},
\] and a duality functor $\bb{D}$ on $\Mod^{qc}_{\n{K}_{\sol}}(\n{D}^{la}(G,K))$, such that for $\n{F}\in \Mod^{qc}_{\n{K}_{\sol}}(\n{D}^{la}(G,K))$ one has 
\[
j_!(\bb{D}(\n{F}))= (j_!\n{F})^{\vee,Rla}.
\]
The functor $\n{F}\mapsto j_! \bb{D}(\n{F})= (j_!\n{F})^{\vee, Rla}$ restricts to Schneider and Teitelbaum's classical anti-equivalence between coadmissible $\D^{la}(G, K)$-modules and admissible locally analytic representations of $G$. 
\end{theoremx}

\begin{remark} \leavevmode
\begin{enumerate}
\item  For the duality functor $\bb{D}$ in the category  $\Mod^{qc}_{\n{K}_{\sol}}(\n{D}^{la}(G,K))$ we refer to Definition \ref{DefinitionDualityFunctor}. It is compatible with the duality functor  in the bigger category $\Mod_{\n{K}_\sol}(\n{D}^{la}(G, K))$ of all solid $\n{D}^{la}(G, K)$-modules, which is given by the formula 
\[
\bb{D}(V) = R \iHom_{\n{D}^{la}(G, K)}(V, \n{D}^{la}(G, K) \otimes_{\n{K}_\sol} \chi^{-1}[d]),
\]where $\chi$ and $d$ are as before. Note that this functor coincides (up to a shift in the cohomological degree, and after fixing a Haar measure) with the one considered  in \cite[Corollary  3.6]{SchTeitDuality} when $G$ is compact. Indeed, the character of \textit{loc. cit.} is the tensor of $\det (\f{g})=\chi^{-1}$ with the unimodular character $\delta_G$ of $G$ seen merely as a  profinite group. Since $G$ is compact and we are working in characteristic zero,  the choice of a Haar measure trivializes $\delta_G$. The duality in the noncompact case requires  the twist by $ \f{d}_G:=\det \f{g}\otimes \delta_G$ of \cite{SchTeitDuality}, but we  do not deal with this  in this article.

\item Even though the result is stated for a compact group $G$, one trivially recovers the anti-equivalence of Schneider-Teiltelbaum for non-compact groups, since the classical notions of admissible and coadmissble are local in $G$, i.e. they only depend on the restriction to an open compact subgroup.

\item  Along the way, the above proposition also answers a question raised in \cite[p. 26]{SchTeitDuality}, concerning the extension of the smooth contragradient functor from the category of admissible smooth representations to the cateory of admissible locally analytic representations.  We refer the reader to Proposition \ref{PropContragradientSmoothLocAn} for the precise answer to Schneider and Teitelbaum's question.
\end{enumerate}
\end{remark}

\subsubsection{Cohomology of representations}

We now explain our cohomology comparison results. Recall that, for a compact $p$-adic Lie group $G$ defined over $\Q_p$, and for any finite dimensional (and hence locally analytic) representation $V$ of $G$, Lazard (\cite{Lazard}) established an isomorphism between the continuous, locally analytic, and Lie algebra cohomology of $V$. In our previous work, we generalized this result for any solid locally analytic representation of $G$, and our proof was homological in nature. Moreover, we also showed an isomorphism between the solid group cohomology of a solid representation and of its locally analytic vectors. The following result generalizes all this to arbitraty $p$-adic Lie groups and also shows how all the comparison isomorphisms are consequence of certain adjunctions, which we think provide a better conceptual understanding of them. \\

There are natural functors\footnote{The first functor is the one that sends a solid $\n{K}_\sol$ module to itself equipped with the trivial action of $G$, the second one is given by viewing a smooth representation as a locally analytic representation, etc.}
\[    \Mod(\n{K}_\sol) \to \Rep^{sm}_{\n{K}_\sol}(G) \xrightarrow{F_1} \Rep^{la}_{\n{K}_\sol}(G_L) \xrightarrow{F_2} \Rep^{la}_{\n{K}_\sol}(G_\qp) \xrightarrow{F_3}  \Rep_{\n{K}_\sol}(G), \]
where we denote by $G_\qp$ the restriction of scalars of $G$ from $L$ to $\qp$, and $G_L = G$ to stress that the group is defined over $L$ in order to avoid confusion. All these functors commute with colimits and hence possess right adjoints. The main idea for our cohomological comparison results is to reinterpret them as formal identities coming from adjunctions and hence reduce them to calculating the right adjoints of the above arrows. Classically, there are many possible cohomology theories associated with $G$ that consider different possible structures of $G$, e.g., continuous, $\qp$ and $L$-locally analytic, smooth and Lie algebra cohomology.

\begin{definition}
We define
\begin{itemize}
\item Solid group cohomology $R \Gamma(G, -) : \Rep_{\n{K}_\sol}(G) \to \Mod(\n{K}_\sol)$,
\item ($\qp$-)Locally analytic group cohomology $R \Gamma^{la}(G_\qp, -) : \Rep^{la}_{\n{K}_\sol}(G_\qp) \to \Mod(\n{K}_\sol)$,
\item ($L$-)Locally analytic group cohomology $R \Gamma^{la}(G_L, -) : \Rep^{la}_{\n{K}_\sol}(G_L) \to \Mod(\n{K}_\sol)$,
\item Smooth group cohomology $R \Gamma^{sm}(G, -) : \Rep^{sm}_{\n{K}_\sol}(G) \to \Mod(\n{K}_\sol)$
\item Lie algebra cohomology $R \Gamma(\f{g}, -) : \Mod_{\n{K}_\sol}(U(\f{g})) \to \Mod(\n{K}_\sol)$,
\end{itemize}
as the right adjoint of the trivial representation functor from $\Mod(\n{K}_\sol)$ to the corresponding category.
\end{definition}

One can check (Proposition \ref{PropComparisonCochainsAndCoho}) that these definitions coincide with the usual definition of cohomology using (continuous, locally analytic, etc...) cochains. Our main key calculation is to show (Proposition \ref{PropositionAdjunctions}) that
\begin{enumerate}
    \item The right adjoint of $F_1$ is given by Lie algebra cohomology $R \Gamma(\f{g}_L, -) := R \iHom_{U(\f{g}_L)}(K, -)$.
    \item The right adjoint of $F_2$ is given by $R \Gamma(\f{k}, -) := R \iHom_{U(\f{k})}(K, -)$, where $\f{k} = \ker(\f{g}_\qp \otimes_\qp L \to \f{g}_L)$.
    \item The right adjoint of $F_3$ is given by the functor of locally analytic vectors $(-)^{Rla}$.
\end{enumerate}
Moreover, the right adjoint to the composition of $F_1 \circ \ldots \circ F_j$ ($j = 1, 2, 3$) can be interpreted as taking smooth vectors in the corresponding category. Analogously, the right adjoint of $F_2 \circ \ldots F_j$ ($ j = 2, 3$) can be interpreted as taking locally $L$-analytic vectors, and so on. Summarizing this, we obtain our third main result.

\begin{theoremx}[Theorem \ref{TheoremComparisonCohomologies}] \label{TheoremComparisonCohomologiesIntro}
We have the following commutative diagram:
\begin{small}
\[
\begin{tikzcd}
& \Rep^{la}_{\n{K}_\sol}(G_\qp) \ar[rr, bend left = 10, "R \Gamma(\f{k}{,} -)"] \arrow[dddr, "R \Gamma^{la}(G_\qp {,} -)" near start] & & \Rep^{la}_{\n{K}_\sol}(G_L) \ar[dr, bend left = 10, "R \Gamma(\f{g}{,} -)"] \ar[dddl, "R \Gamma^{la}(G_L{,}-)"' near start] & \\
\Rep_{\n{K}_\sol}(G) \ar[ur, bend left = 10, "(-)^{Rla}"] \arrow[ddrr, "R \Gamma(G{,} -)"'] & & & & \Rep^{sm}_{\n{K}_\sol}(G) \ar[ddll, "R \Gamma^{sm}(G{,}-)"] \\
& & & & \\
& & \Mod(\n{K}_\sol) & & 
\end{tikzcd}
\]
\end{small}
Moreover, since the embedding $\Rep^{la}_{\n{K}_\sol}(G_\qp)$ in $\Rep_{\n{K}_\sol}(G)$ is fully faithful, we have $R \Gamma(G, V) = R \Gamma(G_\qp, V^{Rla})$ for $V \in \Rep_{\n{K}_\sol}(G)$. In particular, we have that
\[ R \Gamma(G, V) = R \Gamma(G_\qp, V^{Rla}) = R \Gamma^{la}(G_L, V^{R L-la}) = R \Gamma^{sm}(G, R \Gamma(\f{g}, V^{R L-la})), \]
where $V^{R L-la} := R \Gamma(\f{k}, V^{Rla})$ are the $L$-locally analytic vectors of a solid representation of $G$.
\end{theoremx}

\begin{remark} \leavevmode
\begin{enumerate}
    \item When $G$ is compact and $V$ is a finite dimensional representation, the last two equivalences are a classical result of Lazard \cite{Lazard}. When $G$ is given by the $\qp$-points of an algebraic group and $V$ is finite dimensional, Casselman-Wigner generalized Lazard's result in \cite{CasselmanWigner}. The case of $G$  compact, defined over $\Q_p$, and any solid representation $V$ was obtained by the authors in \cite{SolidLocAnRR}.

    \item When $G$ is a $p$-adic reductive group over $\bb{Q}_p$ and $V$ is an admissible Banach representation of $G$, then $V^{Rla} = V^{la}$ and the isomorphism $R \Gamma(G, V) = R \Gamma(G, V^{la})$ was recently and independently shown by Fust in \cite{fust2023continuous} by reducing the problem to the compact case \cite[Theorem 5.3]{SolidLocAnRR} via a Bruhat-Tits building argument.
    
    \item We refer to Remark \ref{RemarkLocAnvectorsSolidRep} for comments on the functor of $L$-locally analytic vectors.
\end{enumerate}
\end{remark}

\subsubsection{Categorical $p$-adic locally analytic Langlands for $\GL_1$}

We conclude this introduction with an application of Theorem \ref{TheoIntroA} to the locally analytic $p$-adic Langlands correspondence for $\GL_1$. We heartily thank Eugen Hellmann for pointing out this application to us.

The $p$-adic local Langlands correspondence aims to establish a correspondence between categories of $p$-adic local Galois representations (or $(\varphi, \Gamma)$-modules) of dimension $n$ and certain category of $p$-adic representations of $\GL_n(\Q_p)$. Recently, inspired by  Fargues--Scholze's geometrization of the classical local Langlands correspondence, Emerton, Gee and Hellmann proposed a categorification of the $p$-adic local Langlands correspondence. In vague terms, they conjecture the existence of a functor from the category of Banach (resp. locally analytic) representations of $\GL_n(\Q_p)$ to certain category of quasi-coherent sheaves on a classifying stack of $(\varphi, \Gamma)$-modules over Fontaines's ring $\n{E}$ (resp. over the Robba ring), satisfying various properties. The Banach case for $\GL_1$ follows from local class field theory, but the locally analytic case already has some subtle difficulties.

Let $\n{X}_1$ be the classifying stack of rank $1$ $(\varphi, \Gamma)$-modules over the Robba ring on affinoid Tate algebras over $\n{K}=(K,K^+)$, c.f. \cite[\S 5]{EGH}. Since every such $(\varphi, \Gamma)$-module is given, up to a twist by a line bundle on the base, by a continuous (and hence locally analytic) character on $\qpe = \zpe \times p^\Z$, this stack is represented (c.f. \cite[\S 7.1]{EGH}) by the quotient 
\[ [(\widetilde{\n{W}} \times \bb{G}_m^{an}) / \bb{G}_m^{an}]\] with trivial action of $\bb{G}_m^{an}$, where $\widetilde{\n{W}}$ is the rigid analytic weight space of $\n{O}_L^\times$ whose points on an affinoid ring $A$ are given by continuous characters $\Hom(\n{O}_L^\times, A^\times)$, and where $\bb{G}_m^{an}$ denotes the rigid analytic multiplicative group. Let $\Mod^{qc}_{\n{K}_\sol}(\n{X}_1)$
be the category of solid quasi-coherent sheaves on $\n{X}_{1}$. In \cite{EGH}, the authors conjecture that the natural functor 
\begin{equation*}
\f{LL}_p^{la} : \Rep_{\n{K}_\sol}^{la}(L^{\times}_{\bb{Q}_p}) \to \Mod_\sol^{qc}(\n{X}_1) 
\end{equation*}
 given by $\f{LL}_p^{la}(V) = \n{O}_{\n{X}_1} \otimes^L_{\n{D}^{la}(L^\times_{\bb{Q}_p}, K)} V$ is fully faithful when restricted to a suitable category of ``tempered'' (or finite slope) locally analytic representations (c.f. \cite[Equation (7.1.3)]{EGH}). Here $L^{\times}_{\bb{Q}_p}$ is the restriction of scalars to $\bb{Q}_p$ of the $p$-adic Lie group $L^{\times}$. On the other hand, for the functor $\f{LL}_p^{la}$ to be fully faithful without restricting to a smaller subcategory of $\Rep^{la}_{\n{K}_\sol}(L^\times_{\bb{Q}_p})$, one can also modify the stack $\n{X}_1$, namely, we consider
\[ \n{X}_1^{mod}:=[\widetilde{\n{W}}\times \bb{G}_m^{alg}/\bb{G}_{m}^{alg}]\]   where $\bb{G}_m^{alg}$ is the analytic space, in the sense of \cite{ClauseScholzeanalyticspaces},  attached to the ring $(K[T^{\pm 1}], K^+)_{\sol}=\n{K}_{\sol}\otimes_{\bb{Z}} \bb{Z}[T^{\pm 1}]$.

In order to describe the category of solid quasi-coherent sheaves on the stacks $\n{X}_1^{mod}$ and $\n{X}_1$ in terms of representation theory, we need to introduce some notation. We let $\s{O}(\bb{G}_{m}^{an})= \varprojlim_{n\to \infty} K\langle  p^nT ,  \frac{p^{n}}{T}\rangle$ and  $\ell^{temp}_{\bb{Z},K}= \s{O}(\bb{G}_m^{an})^{\vee}$ be the Hopf algebras of functions of the group $\bb{G}_m^{an}$ and its dual. We let $\bb{Z}^{temp}$ denote the analytic space defined by the algebra $\ell^{temp}_{\bb{Z},K}$.
We also let $C^{temp}(L_{\bb{Q}_p}^{\times},K)= \s{O}(\widetilde{\n{W}}\times \bb{G}_{m}^{an})^{\vee}$ be the Hopf  algebra of tempered locally analytic functions on $L^\times$. Finally, we let $\Rep^{temp}_{\n{K}_\sol}(L^{\times}_{\bb{Q}_p}) := \Mod^{qc}_{\sol}([*/L^{\times,temp}_{\bb{Q}_p}])$ be the category of tempered (locally analytic) representations of $L^{\times}_{\bb{Q}_p}$.

\begin{theoremx} [Theorem \ref{TheoLLGL1Algebraic}]
 There are natural equivalences of stable $\infty$-categories
\begin{equation*}
 \Mod^{qc}_{\n{K}_\sol}([\bb{Z}/ L^{\times,la}_{\bb{Q}_p}]) \xrightarrow{\sim} \Mod^{qc}_{\n{K}_\sol}(\n{X}_1^{mod}), \;\;\;  \Mod^{qc}_{\n{K}_\sol}([\bb{Z}^{temp}/ L^{\times, temp}_{\bb{Q}_p}])  \xrightarrow{\sim} \Mod^{qc}_{\n{K}_\sol}(\n{X}_1)
\end{equation*}
Furthermore, the functor $\f{LL}_p^{la}$ induces equivalences
\begin{equation*}
\Rep^{la}_{\n{K}_{\sol}}(L^{\times}_{\bb{Q}_p}) \xrightarrow{\sim} \Mod_{\n{K}_\sol}^{qc}(\widetilde{\n{W}}\times \bb{G}_m^{alg}), \;\;\; \Rep^{temp}_{\n{K}_{\sol}}(L^{\times}_{\bb{Q}_p}) \xrightarrow{\sim} \Mod_{\n{K}_\sol}^{qc}(\widetilde{\n{W}}\times \bb{G}_m^{an}).
\end{equation*}  
\end{theoremx}

\subsection*{Acknowledgements}

We thank Eugen Hellmann for pointing out the application of Theorem \ref{TheoIntroA} to the categorical $p$-adic Langlands correspondence, and Arthur--C\'esar Le Bras for many inspiring conversations that concluded in the stacky interpretation of locally analytic representations. We thank Lucas Mann for several discussions on  six-functor formalisms and their connection with representation theory. We thank  Johannes Ansch\"utz,  Ko Aoki,  Yutaro Mikami, C\'edric Pepin, Vincent Pilloni,  Peter Scholze, Matthias Strauch and Otmar Venjakob for their comments and corrections. Finally, we would like to thank the anonymous referee for their numerous and detailed corrections, which have helped us to improve the exposition of this article considerably. The first author was supported by the project ANR-19-CE40-0015 COLOSS. The second author thanks the Max Planck Institute for Mathematics for its hospitality during the preparation and correction of this paper.

\subsection*{Notations}

Throughout this paper we use the language of $\infty$-categories of \cite{higherTopos}, and the techniques of higher algebra from \cite{HigherAlgebra}. We use Clausen and Scholze condensed approach to analytic geometry as presented in the lecture notes \cite{ClausenScholzeCondensed2019,ClauseScholzeanalyticspaces, CondensesComplex}. We refer to \cite{MannSix}  for complete and rigorous proofs of foundational results on the subject, particularly those regarding the set theoretical subtleties in condensed mathematics.  Nevertheless, throughout this paper we will fix an uncountable solid cutoff cardinal $\kappa$ as in \cite[Definition 2.9.11]{MannSix} and work with $\kappa$-small condensed sets, it will be clear from the definitions that  the functors and adjunctions constructed below are independent of $\kappa$, and therefore  that they extend naturally to the full condensed categories.

For $\n{C}$ an $\infty$-category with all small limits and colimits, we let $\Cond(\n{C})$ denote the $\infty$-category of condensed $\n{C}$-objects, see \cite[Definition 2.1.1]{MannSix}. Given $X\in \Cond(\n{C})$ and $S$ a profinite set, we let $\iCont(S,X)$  or  $C(S,X)$ be the object in $\Cond(\n{C})$ whose values at extremally disconnected $S'\in \Extdis$ are $X(S\times S')$. This is still a condensed object by \cite[Corollary 2.1.10]{MannSix} under a mild condition on $\n{C}$ (eg. if it is presentable). In particular, we shall write $\CondSet$,  $\CondAb$ and $\CondRing$ for the categories of condensed sets, abelian groups and commutative rings, respectively.

All the analytic rings considered in this document are assumed  to be animated and complete in the sense of  \cite[Definition 2.3.10]{MannSix}, unless otherwise specified.   Given $\n{A}=(\underline{\n{A}},\n{M})$ a commutative animated analytic ring we shall write $\Mod_{\n{A}}$ for the symmetric monoidal $\infty$-category of analytic $\n{A}$-modules and $\Mod^{\heartsuit}_{\n{A}}$  for the heart of its natural $t$-structure. Given $D$ an $\bb{E}_1$-algebra in $\Mod_{\n{A}}$, we let $\LMod_{\n{A}}(D)$ and $\RMod_{\n{A}}(D)$  be the $\infty$-category of left and right $D$-modules in $\Mod_{\n{A}}$, if it is clear from the context we will simply write $\Mod_{\n{A}}(D)=\LMod_{\n{A}}(D)$.  We say that an analytic ring $\n{A}$ is static if  for all extremally disconnected sets $S$, the object  $\underline{\n{A}}[S]$  is concentrated in cohomological degree $0$. We let $-\otimes^L_{\n{A}}-$ denote the complete tensor product of $\Mod_{\n{A}}$, and $R\iHom_{\n{A}}(-,-)$ the internal Hom space, right adjoint to the tensor.  By \cite[Warning 7.6]{ClausenScholzeCondensed2019}, the tensor $-\otimes^L_{\n{A}}-$ is the left derived functor of the tensor $-\otimes_{\n{A}}-$ if  $\n{A}[S\times T]$ sits in degree $0$ for all extremally disconnected sets. The analytic rings we will consider live over the solid base $\bb{Z}_{\sol}$, so this property is always true for them whenever they are static.

Recall that a map $f:N\to M$ of objects in $\Mod_{\n{A}}$ is called trace class (\cite[Definition 8.1]{CondensesComplex}) if there is a map $\n{A} \to N^{\vee} \otimes_{\n{A}} M$ with $N^{\vee} = R\iHom_{\n{A}}(N,\n{A})$, such that  $f$ factors as 
 \[
 N \to N\otimes_{\n{A}}^L N^{\vee} \otimes_{\n{A}}M \to M.
 \]
 An object $N\in \Mod_{\n{A}}$ is called nuclear (\cite[Definition 13.10]{ClauseScholzeanalyticspaces}) if for every extremally disconnected set $S$, the natural map 
 \[
 \n{A}[S]^{\vee} \otimes^L N (*) \to N(S)
 \]
 is an isomorphism. By \cite[Proposition 13.14]{ClauseScholzeanalyticspaces}, if $N \in \Mod_{\n{A}}$ is nuclear, then for every extremally disconnected  set $S$ and any $M\in \Mod_{\n{A}}$, the natural map
 \[
(R\iHom_{\n{A}}(\n{A}[S], M) \otimes_{\n{A}}^L N )(*)\to (M\otimes_{\n{A}}^L N ) (S)
 \]
is an isomorphism.

\section{Distribution algebras}

We record in this chapter basic properties of the several spaces of functions and algebras of distributions we will be working throughout the text. Most of the results are probably well known but we give statements and proofs for the sake of notation and completeness. Let $L$ be a finite extension of $\bb{Q}_p$ and $\varpi\in L$ a pseudo-uniformizer.  Let $G$ be a  $p$-adic Lie group over $L$. We normalize the $p$-adic absolute value of $L$ such that $|p|=p^{-1}$.

\subsection{Distribution algebras and  spaces of   functions}
\label{ss:SLADistributionAlgebras}

In the present section we will introduce and set up notations for all the distribution algebras and spaces of functions we will use in the article.

\subsubsection{Locally analytic distributions}
 
We start with the introduction of a family of locally analytic distribution algebras for the case of compact Lie groups.  Let $G$ be a compact $p$-adic Lie group of dimension $d$ over $L$. Let $\f{g}$ denote the Lie algebra of  $G$,  and let $\n{L}\subset \f{g}$ be an $\n{O}_L$-lattice such that $[\n{L},\n{L}]\subset p\n{L}$ (we call such a lattice \textit{good}). Let $L_{\sol}[G]$ be the Iwasawa algebra of $G$, i.e., $L_{\sol}[G]= (\varprojlim_{H\subset G}\n{O}_L[G/H] ) [\frac{1}{p}]$ where $H$ runs over all the compact open subgroups.
As it is explained in \cite[\S 5.2]{Emerton}, $\n{L}$ can be integrated to an analytic group $\bb{G}_{\n{L}}$ over $L$ whose underlying adic space can be identified with a polydisc of dimension $d$.  More precisely, let $\f{X}_1, \ldots, \f{X}_d$ be an $\n{O}_L$-basis of $\n{L}$, then the map 
\[
(T_1,\ldots, T_d)\mapsto \exp(T_1\f{X}_1)\hdots \exp(T_d\f{X}_d)
\]
induces an isomorphism of adic spaces between the polydisc $\bb{D}_L^{d}= \Spa(L\langle \underline{T} \rangle, \n{O}_L\langle \underline{T} \rangle)$ and $ \bb{G}_{\n{L}}$.  After shrinking $\n{L}$ if necessary we can assume that $\bb{G}_{\n{L}}(L)\subset G$ is a normal compact open subgroup which is moreover a uniform pro-$p$-group. In the following, we will always assume that $\n{L}$ is small enough such that this holds.

The previous construction can be slightly  generalized as follows. Let $\overline{L}$ be an algebraic closure of $L$, and let $\n{L}\subset \f{g}_{\overline{L}}$ be a free $\n{O}_{\overline{L}}$-lattice such that $[\n{L}, \n{L}]\subset p \n{L}$.  There exists a finite extension $F$ of $L$ such that $\n{L}$ is  defined over $F$,  and one can define an affinoid group $\bb{G}_{\n{L},F}$ over $F$ by integrating $\n{L}$. Furthermore, suppose that the action of $\Gal_{L}$ leaves $\n{L}$ stable, then $\bb{G}_{\n{L},F}$ can be obtained as the base change from $L$ of an affinoid group that we denote as $\bb{G}_{\n{L}}$.  A free lattice $\n{L}\subset \f{g}_{\overline{L}}$ is called \textit{good} if  it is $\Gal_{L}$-stable and  $[\n{L},\n{L}]\subset p \n{L}$. If $\n{L}$ is defined over $F$ we let $\n{L}_F$ denote the $\Gal_{F}$-invariants of $\n{L}$.

\begin{example}
\label{ExamplehGroups}
Let us fix a good  $\n{O}_L$-lattice $\n{L}_0 \subset \f{g}$ with group $\bb{G}_{0}$.   For $h>0$ a rational number,  the lattice $p^{h} \n{L}_0$ over $\f{g}_{\overline{L}}$ is good, and it defines an affinoid subgroup $\bb{G}_{h}\subset \bb{G}_0$ which is nothing but the polydisc of radius $p^{-h}$: 
\[
\bb{G}_h= \bb{G}_0 \big(\frac{\underline{T}}{p^{h}}\big). 
\]
\end{example}

Given a good lattice $\n{L}$ we can also define analytic groups which are Stein spaces, namely, we let $\mathring{\bb{G}}_{\n{L}} = \bigcup_{h>0} \bb{G}_{p^h \n{L}}$. If $\n{L}$ is already defined over $L$ then $\mathring{\bb{G}}_{\n{L}}$ is an open polydisc.

Finally, we can construct affinoid and Stein group neighbourhoods of $G$ by taking finitely many translates of the groups $\bb{G}_{\n{L}}$ and $\mathring{\bb{G}}_{\n{L}}$. Indeed, since $\bb{G}_{\n{L}}(L)$  and $\mathring{\bb{G}}_{\n{L}}(L)$ are normal subgroups of $G$, we can define 
\[
\bb{G}^{(\n{L})} := G \bb{G}_{\n{L}}= \bigsqcup_{g\in G/\bb{G}_{\n{L}}(L) } g \bb{G}_{\n{L}} \mbox{ and } \bb{G}^{(\n{L}^+)}:=G \mathring{\bb{G}}_{\n{L}}= \bigsqcup_{g\in G/\mathring{\bb{G}}_{\n{L}}(L)}   g \mathring{\bb{G}}_{\n{L}}.
\]
If $\n{L}_0$ is a fixed good lattice and $\n{L}=p^h\n{L}_0$  we will simply denote $\bb{G}^{(h)}= \bb{G}^{(\n{L})}$ and $\bb{G}^{(h^+)}= \bb{G}^{(\n{L}^+)}$.

With the previous notations we can now define the following locally analytic distribution algebras and analytic functions.

\begin{definition}
\label{Def:AlgebrasDist}
Let $\n{L} \subset \f{g}_{\overline{L}}$ be a good lattice defined over $F/L$.
\begin{enumerate}
    \item  Let $\bb{G}$ be one of the adic groups $\bb{G}_{\n{L}}$, $\mathring{\bb{G}}_{\n{L}}$, $\bb{G}^{(\n{L})}$ or $\bb{G}^{(\n{L}^+)}$. The space of analytic functions of $\bb{G}$ with values in $L$ is the space $C(\bb{G},L)=\s{O}(\bb{G})$ in the abelian category of solid $L$-vector spaces.  The algebra of distributions of $\bb{G}$ is the dual space $\n{D}(\bb{G},L)=\iHom_{L}(C(\bb{G},L),L )$. If $\n{L}_0$ is fixed as in Example \ref{ExamplehGroups} and  $\n{L}= p^h\n{L}_0$, we will simply denote $\n{D}^{h}(G,L) = \n{D}(\bb{G}^{(h^+)},L)$ and $C^h(G, L) = C(\bb{G}^{(h^+)}, L)$.

    \item  We let $\widehat{U}(\n{L})^+$ be the $\Gal_{F/L}$-invariants of the $p$-adic completion of the  enveloping algebra of $\n{L}_F$. We also denote $\widehat{U}(\n{L})= \widehat{U}(\n{L})^+[\frac{1}{p}]$.

\end{enumerate}

\end{definition}

\begin{remark}\label{rem02jsd}
Fix $h >0$ rational and denote $\mathring{\bb{G}}_{h} = \mathring{\bb{G}}_{p^h \n{L}_0}$, $G_h=\mathring{\bb{G}}_{h}(L)$. Note that, as $C(\bb{G}^{(h^+)}, L) = \prod_{g \in G / G_h} C(\mathring{\bb{G}}_h, L)$, by taking duals,  we have that 
\[
\n{D}^h(G,L)=\n{D}(\mathring{\bb{G}}_{h},L)\otimes_{L_{\sol}} L[G/G_h]
\]
as left $\n{D}(\mathring{\bb{G}}_{h},L)$-modules. A similar description holds as right module. 
\end{remark}

We finally define locally analytic functions and distribution algebras for general $p$-adic Lie groups. 

\begin{definition}\label{defijdasnj}
Let $G$ be a locally profinite $p$-adic Lie group over $L$ and $\f{g}$ its Lie algebra.  

\begin{enumerate}

\item We define the space of locally analytic functions of the Lie algebra $\mathfrak{g}$ as the colimit $C^{la}(\f{g},L):= \varinjlim_{\n{L}\subset \f{g}} C(\bb{G}_{\n{L}}, L)$. Its space of locally analytic distributions is defined as the dual $\n{D}^{la}(\f{g},L)= \iHom_{L}(C^{la}(\f{g},L),L)$ in the abelian category of solid $L$-vector spaces.

    \item   For $G$ profinite we define the space of $L$-analytic functions as $C^{la}(G,L)= \varinjlim_{\n{L}\subset \f{g}} C(\bb{G}^{(\n{L})},L)$. For $G$  locally profinite, its space of  $L$-analytic functions is defined as 
    \[
C^{la}(G,L):=\prod_{g\in G/G_0} C^{la}(gG_0, L)
    \]
where $G_0$ is an open compact subgroup.

  \item  We define the locally analytic distribution algebra of $G$ to be the dual   $\n{D}^{la}(G,L)= \iHom_{L}(C^{la}(G,L),L)$ in the abelian category of solid $L$-vector spaces.

\end{enumerate}
\end{definition}

\begin{remark} \leavevmode
\label{RemarkCofinalityAlgebrasAnalytic}
Let $G$ be a compact $p$-adic Lie group. 
 
\begin{enumerate}
    \item We note that, for $\bb{G} = \bb{G}_{\n{L}}$ or $\bb{G}^{(\n{L})}$ (resp. for $\bb{G} = \mathring{\bb{G}}_{\n{L}}$ or $\bb{G}^{(\n{L}^+)}$), the space $C(\bb{G}, L)$ is a Banach space (resp. a nuclear Fr\'echet space), and the distribution algebra $\n{D}(\bb{G}, L)$ is a Smith space (resp. an $LB$-space of compact type), c.f. \cite{SolidLocAnRR},  \cite{SchTeitDist} or \cite{SchneiderNFA}. In particular, by \cite[Theorem 3.40]{SolidLocAnRR} all the spaces of analytic functions and their corresponding distribution algebras of Definition\ref{Def:AlgebrasDist} are self dual in the abelian category of solid $L$-vector spaces.

    \item The algebra $C^h(G, L)$ is by definition the space of functions of $G$ that are analytic with radius $p^{-h'}$ for any $h' > h$ with respect to the coordinates of $\bb{G}_{\n{L}}$. The reason for considering analytic functions on open balls instead of affinoid balls comes from the fact that the algebras $\n{D}^{h}(G, L)$ are idempotent over $\n{D}^{la}(G, L)$, c.f. Theorem \ref{CoroIdempotentDistributions} below.  

    \item The filtered diagrams $\{C(\bb{G}_{\n{L}},L)\}_{\n{L}}$ and $\{ C(\mathring{\bb{G}}_{\n{L}},L)\}_{\n{L}}$ (resp. $\{C(\bb{G}^{(\n{L})},L)\}_{\n{L}}$,  $\{ C(\bb{G}^{(\n{L}^+)},L)\}_{\n{L}}$ and $\{C^{h}(G,L) \}_{h}$) are ind-isomorphic and their colimit is the space of locally analytic functions of $\f{g}$ (resp. of $G$).  Dually, the cofiltered diagrams of distribution algebras $\{\n{D}(\bb{G}_{\n{L}}, L) \}_{\n{L}}$, $\{\n{D}(\mathring{\bb{G}}_{\n{L}}, L)\}_{\n{L}}$  and $\{ \widehat{U}(\n{L})\}_{\n{L}}$ (resp. the cofiltered diagrams $\{\n{D}(\bb{G}^{(\n{L})}, L)\}_{\n{L}}$, $\{\n{D}(\bb{G}^{(\n{L}^+)}, L)\}_{\n{L}}$ and $\{\n{D}^h(G,L)\}_h$) are pro-isomorphic and their  limit   is equal to $\n{D}^{la}(\f{g}, L)$ (resp. $\n{D}^{la}(G, L)$). Observe moreover that the transition maps of all the previous projective diagrams of distribution algebras are trace class with dense image, which implies that they are also equivalent to diagrams of Banach spaces with trace class and dense image transition maps, see \cite[Corollary 3.38]{SolidLocAnRR}.
    
    \item For $h'>h\geq 0$ we have the inclusions
    \[
    \begin{gathered}
      \n{D}(\mathring{\bb{G}}_{p^{h'}\n{L}_0},L) \subset \n{D}(\bb{G}_{p^{h'}\n{L}_0},L) \subset  \n{D}(\mathring{\bb{G}}_{ p^{h}\n{L}_0},L),  \\
      \n{D}^{h'}(G,L) \subset \n{D}(\bb{G}^{(p^{h'}\n{L}_0)},L) \subset  \n{D}^h(G,L).
      \end{gathered}
    \]
On the other hand, we have that 
\begin{equation}
\label{eqColimitDist}
\n{D}(\mathring{\bb{G}}_{p^h\n{L}_0},L) = \varinjlim_{h'\to (h-\frac{1}{p-1})^{+}} \widehat{U}(p^{-h'} \n{L}_0)
\end{equation}
 for $h>\frac{1}{p-1}$, see  \cite[Proposition 5.2.6]{Emerton} and \cite[Corollary 4.18]{SolidLocAnRR}.

\item Let $\n{L}\subset \f{g}_{\overline{L}}$ be a good lattice defined over $F$ and let $\f{X}_1, \ldots, \f{X}_d$ be a base of $\n{L}_F$ over $\n{O}_F$. One has a power-series description
    \[
\widehat{U}(\n{L})\otimes_{L} F= \widehat{\bigoplus}_{\alpha \in \bb{N}^{d}} F \underline{\f{X}}^{\alpha}. 
    \]
\end{enumerate}
\end{remark}

\begin{remark}\label{RemarkCommentDualsOrNot}
Following \cite[Theorem 3.40]{SolidLocAnRR}, we have a perfect duality for some classical topological  $L$-vector spaces  when considered as objects in the abelian category of solid $L$-vector spaces. For example, we have a duality between Smith and Banach spaces, and between $LB$-spaces of compact type and nuclear Fr\'echet spaces. In some occasions we have a derived duality: if $V$ is a Smith space over $L$, then it is a  projective solid $L$-vector space (c.f. \cite[Lemma 3.8 (2)]{SolidLocAnRR}) and $\underline{\Hom}_L(V,L)=R\underline{\Hom}_L(V,L)$. Similarly, if $V$ is an $LB$-space of compact type, then we can write $V=\varinjlim_n V_n$ as a filtered colimit along injections of Smith spaces and trace class maps (c.f. \cite[Definition 3.34 (2) and Corollary 3.38 (2)]{SolidLocAnRR}), and we have that 
\[
R\underline{\Hom}_{L}(V, L)= R\varprojlim_{n} R\underline{\Hom}(V_n,L) = \varprojlim_n V_n^{\vee}= V^{\vee},
\]
where the first equivalence is formal, the second follows from the projectivity of Smith spaces, 
 and the third from the density of the maps at the level of duals (\cite[Lemma 3.19 (2)]{SolidLocAnRR}) and topological Mittag-Leffler  (\cite[Lemma 3.27]{SolidLocAnRR}).  
However, the derived duality for Banach spaces is known to depend on the axioms of set theory (see \cite{MR4984500}), and we currently do not know whether the same occurs for Fr\'echet spaces. 
\end{remark}

\subsubsection{Smooth distributions}

Next, we introduce smooth distribution algebras.

\begin{definition}\label{defjiasdj}
Let $G$ be a locally profinite group. 

\begin{enumerate}
    \item The space of $L$-valued smooth functions of $G$ is given by
    \[ C^{sm}(G,L)=\prod_{g\in G/G_0} C^{sm}(gG_0,L),\] where $G_0$ is a compact open subgroup of $G$ and $C^{sm}(gG_0,L)=C(gG_0,\mathbb{Z})\otimes_{\mathbb{Z}} L$ is the solid $L$-vector space of $L$-valued smooth functions of $gG_0$. 

	  \item The space $\n{D}^{sm}(G,L)$ of $L$-valued smooth distributions of $G$  is defined as \[ \n{D}^{sm}(G,L)=\iHom_L(C^{sm}(G,L),L). \] 
\end{enumerate}
\end{definition}

\begin{remark} \label{RemSmoothDist}
Let $G$ be a  profinite group. We can write  $C^{sm}(G,L)=\varinjlim_{H\subset G} C(G/H,L)$ where $H$ runs over all the compact open subgroups of $G$. Then $\n{D}^{sm}(G,L)=\varprojlim_{H} L[G/H]$. 
\end{remark}

\begin{lemma}
\label{CoroSmoothAlgebraMatrix}
    Let $G$ be a profinite group. Then 
    \[
    \n{D}^{sm}(G,L)= \prod_{\rho} \rho \otimes \iHom_{\n{D}^{sm}(G,L)}(\rho, \n{D}^{sm}(G,L))
    \]
     where $\rho$ runs over all the irreducible finite dimensional smooth representations of $G$.  In particular:
    \begin{enumerate}
        \item  $\rho$ is a direct summand of $\n{D}^{sm}(G,K)$,  and therefore internally projective. 
        
        \item $L[G/H]$ is an idempotent $\n{D}^{sm}(G,L)$-algebra for all $H\subset G$ normal  open subgroup. 
    \end{enumerate}
\end{lemma}
\begin{proof}
    Any group algebra of a finite group $G_0$ over a field of characteristic zero is isomorphic to the product of  $\rho \otimes \iHom_{L_{\sol}[G_0]}(\rho, L_{\sol}[G_0])$ where $\rho$ runs over all irreducible representations of $G_0$. Since $\n{D}^{sm}(G,L)= \varprojlim_{H} L[G/H]$ if $G$ is compact, the first part of the lemma follows. Assertion (1) is then clear.

    For the assertion (2),  notice that $L[G/H]$ is a direct summand of $\n{D}^{sm}(G,L)$, namely, the projection $\n{D}^{sm}(G,L)\to L[G/H]$ has a section given by the Haar measure of $H$. Writing $\n{D}^{sm}(G,L)= L[G/H]\oplus M$ as  $\n{D}^{sm}(G,L)$-modules, tensoring with $L[G/H]$ gives 
    \[
    L[G/H]= L[G/H]\otimes^L_{\n{D}^{sm}(G,L),\sol} L[G/H] \oplus M \otimes^L_{\n{D}^{sm}(G,L),\sol} L[G/H],
    \]
    but the image of $M$ in $L[G/H]$ is zero, this implies that $L[G/H]=L[G/H]\otimes^L_{\n{D}^{sm}(G,L),\sol} L[G/H]$ proving the corollary. 
\end{proof}

The locally analytic and smooth distribution algebras are self dual for arbitrary $p$-adic Lie groups:

\begin{lemma}\label{xjso2s}
Let $G$ be a locally profinite group. The following hold
\begin{enumerate}
\item We have
    \[ \n{D}^{sm}(G, L) = L_\sol[G] \otimes_{L_\sol[G_0]} \n{D}^{sm}(G_0, L)= \bigoplus_{g \in G / G_0} \n{D}^{sm}(g G_0, L) \]
    or any open compact subgroup $G_0 \subset G$. Moreover, we have $C^{sm}(G,L)=\iHom_L( \n{D}^{sm}(G,L),L)$; thus, the spaces of smooth functions and distributions are reflexive.

\item Suppose that $G$ is a $p$-adic Lie group.      We have 
    \[
    \n{D}^{la}(G,L)= L_\sol[G] \otimes_{L_\sol[G_0]} \n{D}^{la}(G_0, L) = \bigoplus_{g \in G/G_0}  \n{D}^{la}( gG_0, L)
    \]
   for any $G_0$ compact open subgroup of $G$. Moreover, $\n{D}^{la}(G,L)^{\vee}= C^{la}(G,L)$; thus, the spaces of locally analytic functions and distributions are reflexive.

    \item Let $C(G,L)$ be the space of continuous functions from $G$ to $L$.  We have  $\iHom_L(C(G,L),L)= L_{\sol}[G]$. In particular $L_{\sol}[G]$ and $C(G,L)$ are reflexive.

\end{enumerate}

The analogue statemetns of (1) and (2) also hold for right cosets $G_0\backslash G$.

\end{lemma}
\begin{proof}
    The claims for compact groups follows from  \cite[Theorem 3.40]{SolidLocAnRR}. The case of general groups follows  from Lemma \ref{xji92jw} down below and the description of the spaces of functions of $G$ as products of functions on cosets $gG_0$ for $G_0\subset G$ a compact open subgroup. 
\end{proof}

\begin{lemma}\label{xji92jw}
Let $(V_i)_{i\in I}$ be a family of $LB$ spaces over $L$, then 
\[
\iHom_L(\prod_i V_i, L)= \bigoplus_i V_i^{\vee}. 
\]
\end{lemma}
\begin{proof}
Write $V_i = \varinjlim_{j \in J_i} V_{i, j}$. We then have by (AB6) of \cite[Theorem 2.2]{ClausenScholzeCondensed2019}
\begin{equation*}
\prod_{i\in I} V_i = \prod_{i\in I} \varinjlim_{j \in J_i} V_{i, j}
=  \varinjlim_{\substack{\forall i\in  I \\
j_i\in J_i}} \prod_{i\in I} V_{i,j_i}.
\end{equation*}

Therefore, 
    \begin{equation*}
    \iHom_L(\prod_{i\in I} V_i, L) = \varprojlim_{\substack{\forall i\in  I \\
j_i\in J_i}}  \iHom_{L}(\prod_{i\in I} V_{i,j_i},L). 
    \end{equation*}

The product $\prod_{i\in I} V_{i,j_i}$ is a Fr\'echet space and \cite[Theorem 3.40]{SolidLocAnRR} implies that its dual is nothing but $\bigoplus_{i\in I}  V_{i,j_i}^{\vee}$ which yields 
\[
 \iHom_L(\prod_{i\in I} V_i, L)=\varprojlim_{\substack{\forall i\in  I \\
j_i\in J_i}} \bigg(\bigoplus_{i\in I}V_{i,j_i}^{\vee}\bigg). 
\]

By \cite[Theorem 3.40]{SolidLocAnRR} we have that $V_i^{\vee}=\varprojlim_{j\in I_i} V_{i,j}^{\vee}$. Therefore, we have inclusions
\begin{equation} \label{eqDualLocAnFunctions}
\bigoplus_{i\in I} V_i^{\vee}\hookrightarrow  \varprojlim_{\substack{\forall i\in  I \\
j_i\in J_i}} \bigg(\bigoplus_{i\in I}V_{i,j_i}^{\vee}\bigg) \hookrightarrow \prod_{i\in I} V_i^{\vee}.
\end{equation}
 We want to show that the left arrow of \eqref{eqDualLocAnFunctions} is an isomorphism. 

Let $(a_i)_{i\in I}$ be a sequence in the middle term of \eqref{eqDualLocAnFunctions}, we want to show that   all but finitely many elements $a_i\in V_i^{\vee}$ vanish. Suppose the opposite, then we can find infinitely many $i\in I$, and indices $j_i \in J_i$, such that the image $a_{i,j_i}$ of $a_{i}$ in $V_{i,j_i}^{\vee}$ is non-zero, but this contradicts the fact that $(a_{i,j_i})_{i\in I}$ defines an element in $\bigoplus_{i\in I} V_{i,j_i}^{\vee}$. Moreover, the same applies when evaluating at an arbitrary profinite set $S$. The lemma follows. 
\end{proof}

Locally analytic and smooth distribution algebras are related in the following way: 

\begin{lemma} \label{PropPropertiesSmoothDistributionAlgebrab}
Let $G$ be a locally profinite $p$-adic Lie group over $L$.  There is an isomorphism of left $\n{D}^{la}(\f{g},L)$-modules 
    \[ \n{D}^{la}(G, L) = \n{D}^{la}(\f{g},L) \otimes_{L_{\sol}} \n{D}^{sm}(G, L).
    \]
     In particular, the natural map 
    \[
    L\otimes_{\n{D}^{la}(\f{g},L)}^L \n{D}^{la}(G,L)\xrightarrow{\sim} \n{D}^{sm}(G,L)
    \]
    is an isomorphism. The same holds as right $\n{D}^{la}(\f{g},L)$-modules.
\end{lemma}
\begin{proof}
By Lemma \ref{xjso2s} (1) and (2) we can reduce to the case when $G$ is compact.  For any open compact subgroup $H$ of $G$ we have a $\n{D}^{la}(H,L)$-equivariant isomorphism
    \[ \n{D}^{la}(G, L) = \n{D}^{la}(H,L) \otimes_{L_\sol} L[G / H]. \]
    Taking limits for $H\subset G$,   using Remark \ref{RemSmoothDist} and \cite[Lemma 3.28]{SolidLocAnRR}, we get a $\n{D}^{la}(\f{g}, L)$-equivariant isomorphism
    \[
    \n{D}^{la}(G,L)= \n{D}^{la}(\f{g},L)\otimes_{L_{\sol}} \n{D}^{sm}(G,L). 
    \]
\end{proof}

\subsection{Finite projective resolutions and idempotency}

In this section we recollect some elementary algebraic computations of distribution algebras.  The main goal of the section is to show the idempotency of different algebras of locally analytic distributions. 

In the following section we let $G$ be a compact $p$-adic Lie group over $L$ unless otherwise specified. We also fix $\mathcal{L}\subset \f{g}_{\overline{L}}$ a good lattice as in Definition \ref{Def:AlgebrasDist}.  We use the notation $\bb{G}_{h}$, $\mathring{\bb{G}}_h$, $\bb{G}^{(h)}$ and $\bb{G}^{(h^+)}$ instead of $\bb{G}_{p^h\mathcal{L}}$, $\mathring{\bb{G}}_{p^h\n{L}}$, $\bb{G}^{(p^h\n{L})}$ and $\bb{G}^{(p^h\n{L}^+)}$ respectively.  Note in particular that the groups  $\bb{G}^{(h)}$ and $\bb{G}^{(h^+)}$  are not connected and their $L$-points are precisely $G$, while the groups $\bb{G}_{h}$ and $\mathring{\bb{G}}_h$ are geometrically connected and their $L$-points form a basis of compact open subgroups $G_h\subset G$.

\subsubsection{$p$-adic Lie groups over $\qp$}\label{sssIdempotentQp}

We begin by studying the case where $G$ is a compact $p$-adic Lie group over $\qp$. 

\begin{prop}[Lazard-Kohlhaase] \label{PropLazardKohlaase}
Let $G$ be a uniform pro-$p$-group over $\qp$. Then:
\begin{enumerate}
\item There exists a projective resolution of the trivial module $\zp$ of the form
\[ 0 \to \Z_{p, \sol}[G]^{{d \choose d}} \to \cdots \to \Z_{p, \sol}[G]^{{d \choose i}} \to \cdots \to \Z_{p, \sol}[G]^{{d \choose 0}} \to \zp \to 0. \]
\item For any $h > 0$, the above resolution extends by continuity to a resolution
\[ 0 \to \n{D}^h(G, \qp)^{{d \choose d}} \to \cdots \to  \n{D}^h(G, \qp)^{{d \choose i}} \to \cdots \to \n{D}^h(G, \qp)^{{d \choose 0}} \to \qp \to 0. \]
\item Moreover, the above resolution also extends to a resolution
\[ 0 \to \n{D}^{la}(G, \qp)^{{d \choose d}} \to \cdots \to  \n{D}^{la}(G, \qp)^{{d \choose i}} \to \cdots \to \n{D}^{la}(G, \qp)^{{d \choose 0}} \to \qp \to 0. \]
\end{enumerate}
\end{prop}

\begin{proof}
Part (1) is due to Lazard \cite[Lemme V.2.1.1]{Lazard}. Part (2) and (3) are essentially \cite[Theorem 4.4]{Kohlhaase}. Indeed, part (2)   follows from \cite[Theorem 4.4]{Kohlhaase} by taking filtered colimits of suitable distribution algebras $\n{D}_{(h')}(G,\qp)$, see  \cite[Theorem 5.8 and Corollary 5.11]{SolidLocAnRR}.  Part (3) follows from part (2)  after taking limits as $h\to \infty$; observe that taking limits is exact by topological Mittag-Leffler \cite[Lemma 3.27]{SolidLocAnRR}, or by observing that the complexes of $(2)$ admit a chain homotopy, being a colimit of complexes admitting compatible chain homotopies.
\end{proof}

A direct consequence of the proposition is the idempotency of the distribution algebras:

\begin{corollary}\label{coroIdemGQp}
    Let $G$ be a compact $p$-adic Lie group over $\qp$. The maps of associative solid $\qp$-algebras
    \[
    \mathbb{Q}_{p,\sol}[G]\to \n{D}^{la}(G,\qp)\to \n{D}^h(G,\qp)
    \]
    are idempotent. Furthermore, for $G$ a general locally profinite $p$-adic Lie group the map 
    \[
    \mathbb{Q}_{p,\sol}[G]\to \n{D}^{la}(G,\qp)
    \]
    is idempotent. 
\end{corollary}
\begin{proof}
    The claim for compact groups is \cite[Corollary 5.11]{SolidLocAnRR} where the key inputs are  Propositions \ref{PropLazardKohlaase}  and \ref{PropositionHopfAlgebras} (5). The claim for non-compact groups follows from the previous one and the fact that 
    \[
    \n{D}^{la}(G,\qp)=\mathbb{Q}_{p,\sol}[G]\otimes_{\mathbb{Q}_{p,\sol}[G_0]} \n{D}^{la}(G_0,\qp)
    \]
    for a compact open subgroup $G_0\subset G$. 
\end{proof}

\subsubsection{Lie algebras}  

Our next goal is to prove similar idempotency properties as those of Corollary \ref{coroIdemGQp} for the Lie algebras. We need in this case the following Koszul resolutions. We recall that we have fixed a good lattice $\mathcal{L}\subset \f{g}_{\overline{L}}$. For simplicity, and since it will suffice for us, let us assume that $\mathcal{L}$ is defined over $L$, though the following is true for general $\n{L}$ by Galois descent. 

\begin{lemma} 
\label{PropKoszulResolutionDist}
Let $\Kos(\f{g},  U(\f{g}))$ be the standard Koszul resolution of $L$ as $U(\f{g})$-module: 
\[
0 \to U(\f{g}) \otimes \bigwedge^d \f{g} \to \cdots \to U(\f{g})\otimes \f{g} \to U(\f{g})\to L \to 0,
\]
where the differentials are given by 
\begin{align*}
d(v \otimes Z_1 \wedge \ldots \wedge Z_k) &= \sum_{i = 1}^k (-1)^{i + 1} vZ_i \otimes  Z_1 \wedge \ldots \wedge \widehat{Z_i} \wedge \ldots \wedge Z_k \\
&+ \sum_{i<j} (-1)^{i+j} v \otimes [Z_i,Z_j] \wedge \cdots \wedge \widehat{Z}_i \wedge \cdots \wedge \widehat{Z}_j \cdots \wedge Z_k.
\end{align*}
 Let $\Dist$ denote $\Dist(\mathring{\bb{G}}_{\n{L}},L)$, $\widehat{U}(\n{L})$ or $\n{D}^{la}(\f{g},L)$. Then 
\[
\Kos(\f{g}, \Dist):= \Dist \otimes_{U(\f{g}), \sol} \Kos(\f{g}, U(\f{g}))
\]
is a resolution of $L$ as $\Dist$-module.  In particular,  $\Dist \otimes^L_{U(\f{g}), \sol} L = L$ and $\mathcal{D}$ is an idempotent $U(\mathfrak{g})$-algebra. 
\end{lemma} 

\begin{proof} 
Let $U(\n{L})^+$ be the enveloping algebra of $\n{L}$ over $\n{O}_L$.  Let $\Kos(\n{L},U(\n{L})^+)$ be the standard resolution of the trivial representation $\n{O}_L$  and $\epsilon: \Kos(\n{L}, U(\n{L})^+)\to \n{O}_{L}$  the augmentation map.  There is an $\n{O}_L$-linear homotopy $h_\bullet : U(\n{L})^+\otimes \bigwedge^{\bullet} \n{L} \to U(\n{L})^+ \otimes \bigwedge^{\bullet +1} \n{L}$ such that $d_{\bullet+1}h_{\bullet}+ h_{\bullet-1} d_{\bullet}= \id -\epsilon$ (\cite[Theorem 7.7.2]{WeibelHomologicalAlgebra}).  Taking a $p$-adic completion one obtains an homotopy $\widehat{h}_\bullet$ between $\id$ and $\epsilon$ for  $\Kos(\n{L}, \widehat{U}(\n{L})^+)$. Inverting $p$ we have an equivalence $\Kos(\f{g}, \widehat{U}(\n{L})) \xrightarrow{\epsilon} L$. Taking colimits of the Koszul resolutions for $p^{-h'}\n{L}$ as $h'\to (h-\frac{1}{p-1})^+$, one gets an equivalence $\Kos(\f{g}, \n{D}(\mathring{\bb{G}}_{\n{L}},L)) \xrightarrow{\epsilon} L$. Taking limits of $p^{h}\n{L}$ as $h\to \infty$, by topological Mittag-Leffler \cite[Lemma 3.27]{SolidLocAnRR} and Remark \ref{RemarkCofinalityAlgebrasAnalytic}, one gets an equivalence $\Kos(\f{g}, \n{D}^{la}(\f{g},L)) \xrightarrow{\epsilon} L$.  The idempotency of $\n{D}$ over $U(\mathfrak{g})$ follows from Proposition \ref{PropositionHopfAlgebras} (5). 
\end{proof}

The following lemma will be useful for reducing the study of distribution algebras of $p$-adic Lie groups over $L$ to those over $\qp$.

\begin{lemma}
\label{LemmaKoszulCauchiRiemann}
\label{CoroNonIdempotentDLan}
Let $\f{g}$ be a Lie algebra over $L$ and let $\f{h}\subset \f{g}$ be a subalgebra. Let $\n{L} \subset \f{g}$ be a good lattice and let $\n{L}_{\f{h}}= \n{L} \cap \f{h}$. Let $\n{D}(\n{L})$ denote $\widehat{U}(\n{L})$, $\n{D}(\mathring{\bb{G}}_{\n{L}},L)$ or $\n{D}^{la}(\f{g}, L)$ (resp. for $\n{L}_{\f{h}}$), and let $\n{D}(\n{L}/\n{L}_{\f{h}}):= \n{D}(\n{L})\otimes_{\n{D}(\n{L}_{\f{h}})} L$ (i.e. the non-derived tensor).
\begin{enumerate}
    \item  Let $\n{T}\subset \n{L}$ be a free complement of $\n{L}_{\f{h}}$ in  $\n{L}$ with basis $\f{Y}_1,\ldots, \f{Y}_s$ and $\f{t}=\n{T}[\frac{1}{p}]$. Let $\bb{G}_{\n{T}}\subset \bb{G}_{\n{L}}$ be the image by the exponential of the  ordered basis $\f{Y}_1,\ldots, \f{Y}_{s}$, and let $\mathring{\bb{G}}_{\n{T}}= \bigcup_{h>0} \bb{G}_{p^h\n{T}}$ be the open polydisc.  Then we have isomorphisms of solid $L$-vector spaces
    \[
    \n{D}(\n{L}/\n{L}_{\f{h}}) \cong \n{D}(\n{T}) 
    \]
where 
\[
\n{D}(\n{T}) = \begin{cases}
     \widehat{U}(\n{T}):= \widehat{\bigoplus}_{\alpha\in \bb{N}^{s}} L \underline{\f{Y}}^{\alpha}, \\ 
      \n{D}(\mathring{\bb{G}}_{\n{T}},L):= \iHom_L( \s{O} (\mathring{\bb{G}}_{\n{T}}),L) \\
  \n{D}^{la}(\f{t},L) := \varprojlim_{h\to \infty} \widehat{U}(p^h \n{T}).
\end{cases}
\]
    \item We have an isomorphism of right $\n{D}(\n{L}_{\f{h}})$-modules 
    \[
    \n{D}(\n{L})=  \n{D}(\n{T}) \otimes_{L_{\sol}} \n{D}(\n{L}_{\f{h}}).
    \]
    Furthermore, we have an equivalence of  left $\n{D}(\n{L})$-modules $\Kos(\f{h}, \n{D}(\n{L})) \xrightarrow{ \epsilon} \n{D}(\n{L}/\n{L}_{\f{h}})$ where $\Kos(\f{h}, \n{D}(\n{L}))$ is the Koszul complex
    \[
    \Kos(\f{h}, \n{D}(\n{L}))=[0 \to \n{D}(\n{L}) \otimes_L \bigwedge^{\dim \f{h}} \f{h}  \to \cdots \n{D}(\n{L}) \otimes_L \f{h}  \to \n{D}(\n{L})]. 
    \]
     In particular, $\n{D}(\n{L}/\n{L}_{\f{h}})= \n{D}(\n{L})\otimes^L_{\n{D}(\n{L}_{\f{h}}),\sol} L$, and taking $\f{h} = \f{g}$, one recovers Lemma \ref{PropKoszulResolutionDist}.
\end{enumerate}
\end{lemma}
\begin{proof}
The proof of $(1)$ is straightforward. The proof of $(2)$ follows the same lines as those of Lemma \ref{PropKoszulResolutionDist} and we give details. Since $\n{L}= \n{L}_{\f{h}} \oplus \n{T}$, we can write $\mathring{\bb{G}}_{\n{L}} = \mathring{\bb{G}}_{\n{T}} \times \mathring{\bb{G}}_{\n{L}_{\f{h}}}$. Taking global sections one finds that $\n{D}(\mathring{\bb{G}}_{\n{L}},L) = \n{D}(\mathring{\bb{G}}_{\n{T}},L) \otimes^L_{L_{\sol}} \n{D}(\mathring{\bb{G}}_{\n{L}_{\f{h}}},L) $.  We can  then take the  relative de Rham complex of the map 
\[
\mathring{\bb{G}}_{\n{L}} \to \mathring{\bb{G}}_{\n{L}}/ \mathring{\bb{G}}_{\n{L}_{\f{h}}}\cong \mathring{\bb{G}}_{\n{T}},
\]
and by taking duals we find the Koszul complex $\Kos (\f{h}, \n{D}(\mathring{\bb{G}}_{\n{L}} ,L))$, which is quasi-isomorphic to $\n{D}(\mathring{\bb{G}}_{\n{T}},L)$  by the Poincar\'e Lemma. Indeed, the relative de Rham complex is given by the Koszul complex of the action of the relative tangent space $T_{Y/X}$ of $Y=\mathring{\bb{G}}_{\n{L}} \to \mathring{\bb{G}}_{\n{T}}=X$ by derivations of the structure sheaf. Right derivations of $\mathfrak{h}$ on $\mathcal{O}_{Y}$ (induced by the right multiplication of $\mathring{\mathbb{G}}_{\mathcal{L}_{\mathfrak{h}}}$) give rise to an isomorphism of Lie algebroids $\mathcal{O}_{Y}\otimes_L \mathfrak{h}\xrightarrow{\sim} T_{Y/X}$ on $Y$; hence, the Koszul complex for the action of $T_{Y/X}$  is isomorphic to the Koszul complex of the Lie algebra $\mathfrak{h}$. The case for $\n{D}^{la}(\f{g},L)$ and $\n{D}^{la}(\f{h},L)$ is obtained by taking limits along all lattices $\n{L}$ in the previous construction. 

Finally, for $\widehat{U}(\n{L})$ and $\widehat{U}(\n{L}_{\f{h}})$, consider the Koszul complex of $U(\n{L})$ as $\n{L}_\f{h}$-module. Since $U(\n{L})= U(\n{L}_{\f{h}}) \otimes_{\n{O}_L} U(\n{T})$ where $U(\n{T})= \bigoplus_\alpha \underline{\f{Y}}^{\alpha}$, the same argument of \cite[Theorem 7.7.2]{WeibelHomologicalAlgebra} provides an homotopies between $\mathrm{id}$ and the augmentation map 
\[
\Kos(\n{L}_{\f{h}}, U(\n{L})) \xrightarrow{\epsilon } U(\n{T}). 
\]
Taking $p$-adic completions and inverting $p$ one gets the Koszul complex for the $\widehat{U}$-algebras, and  the equality $\widehat{U}(\n{L}) = \widehat{U}(\n{T}) \otimes^L_{L_{\sol}} \widehat{U}(\n{L}_{\f{h}})$. 
\end{proof}

\subsubsection{$p$-adic Lie groups over $L$}

Let now $G$ be a compact $p$-adic Lie group over $L$. In this section we prove that the relevant distribution algebras attached to $G$ are idempotent. This will follow formally from the case of $p$-adic Lie groups and Lie algebras dealt with before. 

 Let $\tilde{\f{g}}$ be the Lie algebra $\f{g}$ seen as a Lie algebra over $\bb{Q}_p$, similarly we let $\widetilde{G}$ be the restriction of $G$ to $\bb{Q}_p$. Take $\f{k} = \ker( \widetilde{\f{g}} \otimes_{\bb{Q}_p} L\to \f{g})$ the subalgebra of $\widetilde{\f{g}} \otimes_{\Q_p} L$.  Let $\n{L}\subset \f{g}_{\overline{L}}$ be a good lattice and let  $\widetilde{\n{L}}$ be its restriction to $\bb{Q}_p$, i.e. the lattice obtained by its $\Gal_{\bb{Q}_p}/ \Gal_{L}$-translates in 
 \[
 \widetilde{\f{g}}_{\overline{L}} = \f{g} \otimes_{\bb{Q}_p} \overline{L} = \prod_{\sigma:L\to \overline{L}} \f{g}_{\sigma, \overline{L}}. 
 \]

\begin{lemma}
\label{CoroCauchiRiemann}
\label{LemmmaBaseChangeDistAlgebras}
\label{LemmaTensorProducthansmoothDistributions}  The following holds:

 \begin{enumerate}
     \item  Let $\n{D}$ denote one of the algebras $\n{D}(\bb{G}^{(\n{L}^+)},L)$, $\widehat{U}(\n{L})$, $\n{D}(\mathring{\bb{G}}_{\n{L}},L)$,  $\n{D}^{la}(\f{g},L)$ or $\n{D}^{la}(G,L)$. Let $\widetilde{\n{D}}$ be the analogue algebra associated with $\widetilde{G}$ and $\widetilde{\n{L}}$ over $\mathbb{Q}_p$.   Then the natural map   
   \[
   (\widetilde{\n{D}}\otimes_{\bb{Q}_p} L) \otimes_{\n{D}^{la}(\f{k},L)}^LL\xrightarrow{\sim} \n{D}
   \]
   is an equivalence. The same holds as left $\mathcal{D}^{la}(\f{k},L)$-modules. 

  \item  The natural map   
    \[ 
    \n{D}(\bb{G}^{(\n{L}^+)}, L) \otimes^L_{\n{D}^{la}(G, L), \sol} \n{D}^{sm}(G, L) \xrightarrow{\sim} L_\sol[G / \mathring{\bb{G}}_{\n{L}}(L)], 
    \] 
    is an isomorphism.

     \item The morphism of algebras $
\n{D}^{la}(G,L)\to \n{D}(\bb{G}^{(\n{L}^+)},L)$
 is idempotent.
     
 \end{enumerate}
 \end{lemma}

\begin{proof} \leavevmode

\begin{enumerate}
 
    \item   For the algebras $\widehat{U}(\n{L})$, $\n{D}(\mathring{\bb{G}}_{\n{L}},L)$,  $\n{D}^{la}(\f{g},L)$  this follows from Lemma \ref{LemmaKoszulCauchiRiemann} (2) by taking $\widetilde{\f{g}} \otimes_{\bb{Q}_p} L$ for $\f{g}$,   $\f{k}$  for $\f{h}$ and $\widetilde{\mathcal{L}} \otimes_{\bb{Z}_p} \n{O}_L$ for $\mathcal{L}$. Indeed, in the notation of the lemma, $\widetilde{\n{D}} \otimes_{\bb{Q}_p} L$ corresponds to the algebra $\n{D}(\n{L})$, $\n{D}$ corresponds to $\n{D}(\n{L}_{\f{h}})$ and the base change $(\widetilde{\n{D}} \otimes_{\bb{Q}_p} L) \otimes_{\n{D}^{la}(\f{k}, L)}^L L$ corresponds to $\n{D}(\n{L} / \n{L}_{\f{h}})$; note that this base change is non-derived and coincides with the base change over $\n{D}(\n{L}_{\f{h}})$ by Lemma \ref{PropKoszulResolutionDist}.
    
For the algebra $\n{D}(G^{(\n{L}^+)},L)$ one reduces to the previous case by  writing the algebras as  finite products of translates of $\n{D}(\mathring{\bb{G}}_{\n{L}},L)$ as in Remark \ref{rem02jsd}. For  $\n{D}^{la}(G,L)$, using Lemma \ref{PropPropertiesSmoothDistributionAlgebrab} we get
\[ \widetilde{\n{D}} = \n{D}^{la}(\widetilde{G}, L) = \n{D}^{la}(\widetilde{\f{g}}, L) \otimes_{L_\sol} \n{D}^{sm}(\widetilde{G}, L).\]
Observing that $\n{D}^{sm}(\widetilde{G}, L) = \n{D}^{sm}(G, L)$, the result reduces to the case of $\n{D}^{la}(\widetilde{\f{g}}, L)$.
 
    \item By Lemma \ref{PropPropertiesSmoothDistributionAlgebrab} we have that 
      \[ 
    \n{D}(\bb{G}^{(\n{L}^+)}, L) \otimes^L_{\n{D}^{la}(G, L), \sol} \n{D}^{sm}(G, L) =\n{D}(\bb{G}^{(\n{L}^+)}, L) \otimes^L_{\n{D}^{la}(\f{g}, L), \sol} L.
    \] 
By Remark \ref{rem02jsd} we can write 
\[
\n{D}(\bb{G}^{(\n{L}^+)}, L) = L[G/\mathring{\bb{G}}_{\n{L}}(L)] \otimes_{L,\sol} \n{D}(\mathring{\bb{G}}_{\n{L}},L),
\]
thus the statement reduces to proving that the map 
\[
\n{D}(\mathring{\bb{G}}_{\n{L}},L)\otimes_{\n{D}^{la}(\f{g}, L),\square} L\to L
\]
is an equivalence which follows from idempotence (c.f. Lemma \ref{PropKoszulResolutionDist}). 
     
 \item      Finally, the claim for the idempotency follows from the idempotency of the algebras over $\Q_p$ of Corollary \ref{coroIdemGQp}, Lemma \ref{lemmaIdempotentBaseChange} and the relation 
     \[
     \n{D}=(\widetilde{\n{D}}\otimes_{\qp} L)\otimes_{\n{D}^{la}(\f{k}, L)} L. 
     \]

\end{enumerate}
 \end{proof}

\begin{lemma}
    Let $G$ be a compact $p$-adic Lie group over $L$. $L$ is a perfect $\mathcal{D}^{la}(G,L)$-complex of perfect amplitude $[- \dim_L G,0]$. In  particular $L$ is also a perfect  $\mathcal{D}(\bb{G}^{(\n{L}^+)},L)$-complex of perfect  amplitude  $[- \dim_L G,0]$.  
    
\end{lemma}
\begin{proof}
By Lemma \ref{PropPropertiesSmoothDistributionAlgebrab} we have that 
\[
\n{D}^{sm}(G,L)= \n{D}^{la}(G,L)\otimes_{\n{D}^{la}(\f{g},L)}^L L
\]
and, by the Koszul complex of Lemma \ref{PropKoszulResolutionDist} for $\n{D} := \n{D}^{la}(\f{g}, L)$, we know that $\n{D}^{sm}(G,L)$ is a perfect complex of finite  projective $\n{D}^{la}(G, L)$-modules of amplitude $\leq \dim_L G$. Finally, the trivial representation $L$ is a direct summand of $\n{D}^{sm}(G,L)$ as $\n{D}^{la}(G, L)$-module, this implies that it is a perfect complex with same bound for the perfect amplitude as wanted. The statement for $\mathcal{D}(\bb{G}^{(\n{L}^+)},L)$ follows from the idempotency over $\n{D}^{la}(G,L)$ and the fact that $L$ is a $\mathcal{D}(\bb{G}^{(\n{L}^+)},L)$-module. 
\end{proof}

For convenience for the reader we summarize the main results of this section in the following theorem.     

\begin{theorem}
\label{CoroIdempotentDistributions}
\label{TheoremIdempotentPefect}
Let $G$ be a compact $p$-adic Lie group over $L$ with Lie algebra $\f{g}$. We fix $\n{L}\subset \f{g}_{\overline{L}}$ a good lattice. Then the natural maps of solid distribution algebras 
\[
\n{D}^{la}(G,L)\to \n{D}(\bb{G}^{(\n{L}^+)}, L) \mbox{ and  }  U(\f{g})\to \n{D}^{la}(\f{g},L)\to  \n{D}(\mathring{\bb{G}}_{\n{L}},L)
\]
are idempotent. In particular, we have fully faithful embeddings of left modules 
\[
\LMod_{L_{\sol}}(\n{D}(\bb{G}^{(\n{L}^+)}, L))\subset  \LMod_{L_{\sol}}(\n{D}^{la}(G,L))
\]
\[ \LMod_{L_{\sol}}(\n{D}(\mathring{\bb{G}}_{\n{L}},L)) \subset \LMod_{L_{\sol}}(\n{D}^{la}(\mathfrak{g},L)) \subset  \LMod_{L_{\sol}}(U(\mathfrak{g}))\]
(resp. for right modules). 
Furthermore, as a module over any of the previous distribution algebras, the trivial representation $L$ is a perfect module with projective amplitude $[-\dim_L G,0]$. 
\end{theorem}
    
\begin{remark}
It is not true that the distribution algebra $\n{D}(\bb{G}^{(h)}, K)$ is an idempotent $\n{D}^{la}(G,K)$-algebra for general $G$. For example, if $G=\bb{Z}_p$, then $\n{D}(\bb{G}^{(h)},\bb{Q}_p)$ can be described as the generic fiber of the  formal complete PD-envelope of a polynomial algebra $\bb{Z}_p [X]$,  which is  not an idempotent $\bb{Z}_p[X]$-algebra, and hence neither as a $\n{D}^{la}(G, L)$-algebra.
\end{remark}

\section{Solid locally analytic representations}

\label{s:SolidLocAnRevisited}

In \cite{SolidLocAnRR} the authors introduced the concept of a solid locally analytic representation for  compact $p$-adic Lie groups over $\bb{Q}_p$.  The goal of this section is to extend the main results of \textit{loc. cit.} to the case where $G$ is a locally profinite  $p$-adic Lie group defined over a finite extension of $\bb{Q}_p$.

Let $L$ be a finite extension of $\bb{Q}_p$ and $\varpi\in L$ a pseudo-uniformizer. Let $\n{K}=(K,K^+)$ be a complete non-archimedean  field extension of $L$ and $\n{K}_{\sol}$ the associated analytic ring.   Let $G$ be a  $p$-adic Lie group over $L$.  In \S \ref{ss:SLALocallyAnalyticVectors}, motivated from the main theorems of  \cite{SolidLocAnRR},  we define the derived $L$-locally analytic vectors of $\n{K}_{\sol}$-complete solid $\n{D}^{la}(G,K)$-modules. We will show  that they can be recovered as the $\bb{Q}_p$-locally analytic vectors which are killed by some ``Cauchy-Riemann equations''.  In  \S \ref{ss:SLACategoriesofRepresentations} we define the $\infty$-category  of $\n{K}_{\sol}$-linear  locally analytic representations of $G$, which will be a full subcategory of the category of $\n{K}_{\sol}$-complete  solid $\n{D}^{la}(G,K)$-modules, where $\n{D}^{la}(G,K)$ is the locally analytic $K$-valued distribution algebra of $G$. If in addition $G$ is defined over $\bb{Q}_p$, the $\infty$-category of $\n{K}_{\sol}$-linear locally analytic representations is itself a full subcategory of $\n{K}_{\sol}$-linear solid $G$-representations, i.e. $\n{K}_{\sol}$-complete modules over the $\n{K}_{\sol}$-linear  group algebra $\n{K}_{\sol}[G]$.
Finally, in \S \ref{ss:SomeAdditionalResults}. we give sufficient conditions for a solid representation to be locally analytic.

\subsection{Locally analytic and smooth functions valued in solid vector spaces}\label{ss:FunctionsSolidV}

Let $G$ be a $p$-adic Lie group over a finite extension $L$ of $\qp$ and let $\n{K} = (K, K^+)$ be a complete non-archimedean extension of $L$. We denote $\n{K}_\sol$ the analytic ring associated with $\n{K}$. In the following we review the definition of locally $L$-analytic vectors of solid $G$-modules on $\n{K}_\sol$-vector spaces.   We shall fix a good lattice $\n{L}_0\subset \f{g}$ defined over $L$, and for $h>0$ we let $\bb{G}^{(h)}$ and $\bb{G}^{(h^+)}$ denote the analytic groups $\bb{G}^{(p^h\n{L}_0)}$ and $\bb{G}^{(p^h\n{L}^+_0)}$ containing $G$ (resp. we let  $\bb{G}_{h}$ and $\bb{G}_{h^+}$ denote $\bb{G}_{p^h\n{L}_0}$ and $\mathring{\bb{G}}_{p^h\n{L}_0}$).

We  define locally analytic functions on $G$ taking values in a solid vector space $V$. Recall from \cite{SolidLocAnRR} that we have defined analytic rings \[ C(\bb{G}^{(h)}, L)_\sol = (C(\bb{G}^{(h)}, L), C(\bb{G}^{(h)}, \n{O}_L))_\sol \] in order to define $h$-analytic and locally analytic vectors of a solid representation. The following Lemma says basically that, in the limit, the analytic structure becomes trivial.

\begin{lemma}
\label{LemmaIndependenceTensorColimit}
Let $h'>h$, we have natural maps of analytic rings 
\[(\s{O}(\bb{G}^{(h)}), \s{O}^+(\bb{G}^{(h)}))_{\sol} \to (\s{O}(\bb{G}^{(h')}), \n{O}_L)_{\sol} \to (\s{O}(\bb{G}^{(h')}), \s{O}^+(\bb{G}^{(h')}))_{\sol}.\] In particular for $V\in \Mod(L_{\sol})$ we have maps
\[
C(\bb{G}^{(h)},L)_{\sol} \otimes_{L_{\sol}}^L V\to C(\bb{G}^{(h')},L) \otimes_{L_{\sol}}^L V \to C(\bb{G}^{(h')}, L)_{\sol} \otimes^L_{L_\sol} V. 
\]
\end{lemma}

\begin{proof}
  By \cite[Lemma 3.31]{Andreychev} one has that for an affinoid ring $(A,A^+)$, $(A,\n{O}_L)_{\sol}= (A, A^{min,+})_{\sol}$ where $A^{min,+}$ is the integral closure of $\n{O}_L+ A^{00}$. The lemma   follows from \cite[Proposition 3.34]{Andreychev} and the fact that we have morphisms of Huber pairs $(\s{O}(\bb{G}^{(h)}), \s{O}^+(\bb{G}^{(h)}) ) \to ( \s{O}(\bb{G}^{(h')}), \s{O}(\bb{G}^{(h')})^{min,+} )\to ( \s{O}(\bb{G}^{(h')}), \s{O}^+(\bb{G}^{(h')}))$. Indeed, if $\frac{T}{p^h}$ denotes a variable of the group $\bb{G}^{(h)}$, one can write $\frac{T}{p^h}=  p^{h'-h}\frac{T}{p^{h'}}$, proving that the image of $\frac{T}{p^h}$ in $\s{O}(\bb{G}^{(h')})$ is topologically nilpotent. 
\end{proof}

\begin{definition} Let $V\in \Mod(L_{\sol})$, we define the following spaces of functions with values in $V$. 
\begin{enumerate}
    \item  For $G$ compact the space of $\bb{G}^{(h)}$-analytic functions 
    \[
    C(\bb{G}^{(h)},V):= C(\bb{G}^{(h)},L)_{\sol} \otimes_{L_{\sol}}^L V.
    \]

    \item  For $G$ compact the space of $\bb{G}^{(h^+)}$-analytic functions 
    \[
    C^h(G,V) = R\varprojlim_{h'>h} C(\bb{G}^{(h')},V) = R\varprojlim_{h'>h} (C(\bb{G}^{(h')},L)\otimes^L_{L_{\sol}} V) 
    \]
    where the second equality holds by Lemma \ref{LemmaIndependenceTensorColimit}. 

\item  For $G$ arbitrary the space of locally analytic functions 
\[
C^{la}(G, V):= \prod_{g\in G/G_0} (C^{la}(gG_0,L)\otimes^L_{L_{\sol}} V)
\]
with $G_0\subset G$ an open compact subgroup. 
\end{enumerate}
\end{definition}

\begin{remark}
Let $\n{K}=(K,K^+)$ be as above.  When $V = K$ is as above, the above definition gives the classical spaces of $K$-valued (locally) analytic functions. Since the spaces of functions $\n{C}$ are either Banach, Fr\'echet spaces or $LB$-type spaces, they are nuclear solid $L$-vector spaces by \cite[Proposition 3.29]{SolidLocAnRR}, and \cite[Proposition 2.3.22 (ii)]{MannSix} implies that, for $G$ compact, the base change along $\n{K}_{\sol}$ agrees with the solid product over $L$, that is, the natural map 
\[
\n{C}\otimes_{L_{\sol}}^L K \to \n{C}\otimes_{L_\sol}^L \n{K}_{\sol}
\]
is an equivalence. 
\end{remark}

\begin{remark} \label{RemarkClaColimit}
    Let $G$ be a compact $p$-adic Lie group and $V\in \Mod(L_{\sol})$. Then we have that 
    \[
    C^{la}(G,V)= \varinjlim_{h} C(\bb{G}^{(h)}, L)\otimes_{L_{\sol}}^L V = \varinjlim_{h} C(\bb{G}^{(h)}, V) = \varinjlim_{h} C^h(G,V)= \varinjlim_{h} C^{h}(G,L)\otimes^L_{L_{\sol}}V,
    \]
    where the first equality is by definition, the  others follow by Lemma \ref{LemmaIndependenceTensorColimit} and the cofinality of the algebras in Remark \ref{RemarkCofinalityAlgebrasAnalytic} (3). 
\end{remark}

\begin{definition} \leavevmode
We define $\n{D}(\bb{G}^{(h)}, K)$, $\n{D}^{h}(G, K)$ and $\n{D}^{la}(G, K)$ as the base change $- \otimes_{L_\sol} \n{K}_\sol$ of the corresponding distribution algebras over $L$, i.e., the $L$-linear duals of $C(\bb{G}^{(h)}, L)$, $C^h(G, L)$ and $C^{la}(G, L)$, respectively.
\end{definition}

\begin{remark}
We observe that, as the spaces $\n{D}^{h}(G, L)$ and $\n{D}^{la}(G, L)$ are nuclear $L_\sol$-vector spaces, the base change to $\n{K}_\sol$ coincides with the extension of scalars to $K$; in particular it is independent of $K^+$. However, the space $\n{D}(\bb{G}^{(h)}, L)$ is a Smith space and its base change to $\n{K}_\sol$ does depend on $K^+$. Since we will not be using this space that often, and since $\n{K}$ will remain fixed along the paper, we will allow ourselves this abuse of notation.
\end{remark}

\begin{remark}
    When $G$ is compact and $L=\bb{Q}_p$, the notation of \cite{SolidLocAnRR} and the one presented in this paper agree for the spaces of functions, i.e.  $C(\bb{G}^{(h)},K)$ and $C(\bb{G}^{(h^+)},K)$. Notice however that the distribution algebras $\n{D}(\bb{G}^{(h)}, K)$ and  $\n{D}(\bb{G}^{(h^+)}, K)$  are written, respectively,  as $\n{D}^{(h)}(G,K)$ and $\n{D}^{(h^+)}(G,K)$ in \textit{loc. cit.}.  In the current paper we are writing $\n{D}^{h}(G,K)= \n{D}(\bb{G}^{(h^+)},K)$ and $C^{h}(G,K)= C(\bb{G}^{(h^+)},K)$ instead since these are the spaces that we use most often, we apologize for the discrepancy in the notations. 
\end{remark}

\subsection{Locally analytic vectors}
\label{ss:SLALocallyAnalyticVectors}

We keep the same notations as before. In particular, $G$ denotes a $p$-adic Lie group over $L$. We will now define and study the functor of $L$-(locally) analytic vectors.

\begin{lemma}
\label{LemmaFunctorFunctionsAndAction}
\leavevmode
\begin{enumerate}
    \item Let $G$ be a compact group, then the functors $V\mapsto C(\bb{G}^{(h)},V)$ and $V\mapsto C^h(G,V)$ for $V\in \Mod({\n{K}_\sol})$ are naturally promoted to exact functors
    \[
    \Mod_{\n{K}_\sol}(\n{D}^{la}(G,K)) \to \Mod_{\n{K}_\sol}(\n{D}^{la}(G^3 ,K)). 
    \]
    
    \item  Let $G$ be arbitrary, then the functor $V\mapsto C^{la}(G,V)$ for $V\in \Mod({\n{K}_\sol})$ is naturally promoted to an exact functor 
    \[
     \Mod_{\n{K}_\sol}(\n{D}^{la}(G,K)) \to \Mod_{\n{K}_\sol}(\n{D}^{la}(G^3,K)).
    \]
    
\end{enumerate} 
Moreover, the functors $V \mapsto C(\bb{G}^{(h)},V)$ and $V\mapsto C^{la}(G,V)$ are exact in the abelian categories.

\end{lemma}
\begin{proof}
   For the compact case it suffices to prove the lemma for $C(\bb{G}^{(h)},-)$, namely the other functors are constructed as limits or colimits of this. But then by \cite[Corollary 2.19]{SolidLocAnRR} we have 
   \[
   C(\bb{G}^{(h)},V)= R\iHom_{K}(\n{D}(\bb{G}^{(h)},K), V), 
   \]
   as $\n{D}(\bb{G}^{(h)},K)$ is  a $\n{D}^{la}(G,K)$-algebra one has the desired  left and right natural actions of $\n{D}^{la}(G\times G,K)=\n{D}^{la}(G,K) \otimes^L_{\n{K}_{\sol}} \n{D}^{la}(G,K)$ on $C(\bb{G}^{(h)},V)$. On the other hand, $\n{D}(\bb{G}^{(h)},L)$ is a Smith space, so  projective as $L_{\sol}$-vector space by \cite[Lemma 3.8 (2)]{SolidLocAnRR}, hence its base change along $L_\sol \to \n{K}_\sol$ remains projective as $\n{K}_\sol$-vector space. This implies that $V\mapsto C(\bb{G}^{(h)},V)$ is exact in the abelian category. If in addition $V$ is a $\n{D}^{la}(G,K)$-module then one has the full action of $\n{D}^{la}(G^3,K)$   as wanted. 

   In the non-compact case, note that  we have natural equivalences
   \[
   C^{la}(G,V)= R\iHom_{\n{D}^{la}(G_0,K)} ( \n{D}^{la}(G,K), C^{la}(G_0, V))
   \]
   for both the left or right regular action of $\n{D}^{la}(G_0,K)$ on $C^{la}(G_0,V)$ and any compact open subgroup $G_0\subset G$.  This endows $C^{la}(G,V)$ with  commuting left and right regular action of $\n{D}^{la}(G,K)$, if in addition $V$ is a $\n{D}^{la}(G,K)$-module then we have the compatible action of \[
   \n{D}^{la}(G^3,K) = \n{D}^{la}(G,K)\otimes^L_{\n{K}_{\sol}}\n{D}^{la}(G,K)\otimes^L_{\n{K}_{\sol}}\n{D}^{la}(G,K)
   \]
   as desired. Finally, since $C^{la}(G,V)=\varinjlim_{h} C(\bb{G}^{(h)},V)$, the functor $V\mapsto C^{la}(G,V)$ is exact in the abelian category. 
\end{proof}
\begin{remark}
    The action of $G^3$ on a function $f$ in any of the three cases is heuristically given by $((g_1,g_2,g_3)\star  f) (h) =  g_3\cdot f(g_1^{-1}h g_2)$. If $V$ arises as the solid vector space attached to a locally convex vector space then the action of $G\times G\times G$ is  given precisely by these formulas.  
\end{remark}

Given $I\subset \{1,2,3\}$ a non-empty subset and $V\in \Mod_{\n{K}_{\sol}}(\n{D}^{la}(G^3, K))$ we let $V_{\star_I} \in \Mod_{\n{K}_{\sol}}(\n{D}^{la}(G, K))$ be the restriction of $V$ to the $I$-diagonal of $\n{D}^{la}(G^3,K)$, i.e. $V$ equipped with he $\D^{la}(G, K)$-module structure  induced by the embedding $\iota_I : G \to G^3$,  $\iota_I(g)_j = g$ if $j \in I$ and $\iota_I(g)_j = e_G$ if $j \notin I$, where $e_G \in G$ denotes the identity element. 

\begin{definition}
\label{DefinitionLocAnVectors}
Let $G$ be a $p$-adic Lie group over $L$.  
\begin{enumerate}
    \item For $G$ compact  the functor of (derived)  $\bb{G}^{(h)}$-analytic vectors $(-)^{Rh-an}: \Mod_{\n{K}_{\sol}}(\n{D}^{la}(G, K)) \to \Mod_{\n{K}_{\sol}}(\n{D}^{la}(G, K))$ is defined as 
    \[
    \begin{aligned}
V^{Rh-an} & := R\iHom_{\n{D}^{la}(G, K)}(K, (C(\bb{G}^{(h)},V)_{\star_{1,3}}),
\end{aligned}
    \]
    where the action of $\n{D}^{la}(G, K)$ on $V^{Rh-an}$ is induced by the $\star_{2}$-action (the right regular action). Similarly,  the  (derived)  $\bb{G}^{(h^+)}$-analytic vectors is the functor on solid $\n{D}^{la}(G,K)$-modules given by 
    \[
    V^{Rh^+-an}:= R\varprojlim_{h'>h} V^{Rh'-an}= R\iHom_{\n{D}^{la}(G,K)}(K, C^{h}(G,V)_{\star_{1,3}}).
    \]
    If $V\in \Mod^{\heartsuit}_{\n{K}_{\sol}}(\n{D}^{la}(G, K))$ we let $V^{h-an}$ and $V^{h^+-an}$ denote the $H^0$ of their derived analytic vectors.  

\item For $G$ compact, we say that an object $V\in \Mod_{\n{K}_{\sol}}(\n{D}^{la}(G,K))$ is $h$-analytic (resp. $h^+$-analytic) if the natural arrow $V^{Rh-an}\to V$ (resp. $V^{Rh^+-an}\to V$) is an equivalence. If $V\in \Mod_{\n{K}_{\sol}}^{\heartsuit}(\n{D}^{la}(G,K))$, we say that $V$ is non-derived $h$-analytic   if the map $V^{h-an}\to V$ is an equivalence (resp. for $h^+$).

\item  For $G$ arbitrary we define the functor of locally analytic vectors $(-)^{Rla}:\Mod_{\n{K}_{\sol}}(\n{D}^{la}( G,K)) \to \Mod_{\n{K}_{\sol}}(\n{D}^{la}( G,K))$ as 
\[
V^{Rla}=R\iHom_{\n{D}^{la}(G,K)}(K, C^{la}(G,V)_{\star_{1,3}})
\]
where we view $V^{Rla}$ endowed with the $\star_{2}$-action of $\n{D}^{la}(G, K)$. 

\item For $G$ arbitrary,  we say that an object $V\in \Mod_{\n{K}_{\sol}}(\n{D}^{la}(G,K))$ is a solid locally analytic representation of $G$ if the natural arrow $V^{Rla}\to V$ is an equivalence.  If $V\in \Mod^{\heartsuit}_{\n{K}_{\sol}}(\n{D}^{la}(G,K))$  we write $V^{la}:= H^0(V^{Rla})$. If $V\in \Mod_{\n{K}_{\sol}}^{\heartsuit}(\n{D}^{la}(G,K))$, we say that $V$ is non-derived locally analytic if $V^{la} \to V$ is an isomorphism. 
\end{enumerate}
\end{definition}

\begin{remark}
The distinction between derived and non derived locally analytic representations might look subtle at the beginning, we will see in Proposition \ref{PropPropertiesCategoryLocAn} that there is no actual difference. 
\end{remark}

\begin{remark}
\label{RemarkLocAnvectorsSolidRep}

\begin{enumerate}
\item The definition of locally analytic vectors might seem slightly strange since we are taking as an input a module over the distribution algebra instead of a solid representation of $G$ as it is usual. Note that, for any $V \in \Mod_{\n{K}_{\sol}}(\n{K}_\sol[G])$ one can define the $L$-analytic vectors of $V$ as
\[ V^{RL-la} := R \iHom_{\n{K}_\sol[G]}(K, C^{la}(G, V)_{\star 1, 3}). \] 
If $V^{R\bb{Q}_p-la}$ denotes the $\bb{Q}_p$-locally analytic vectors of $V$ seen as a  $G_{\bb{Q}_p}$-representation,  and $\f{k}:= \ker(\f{g}\otimes_{\bb{Q}_p} L \to \f{g})$, we have that 
\[
V^{RL-la}= R\Gamma(\f{k}, V^{R\bb{Q}_p-la}),
\]
see Theorem \ref{TheoremComparisonCohomologies}.  The non-derived variant of this functor is usually called the functor of $L$-analytic vectors.

\item If $G=G_{\bb{Q}_p}$ is defined over $\bb{Q}_p$, then $\n{D}^{la}(G, K)$ is an idempotent algebra over $\n{K}_\sol[G]$ and the inclusion of $\Mod_{\n{K}_\sol}(\n{D}^{la}(G, K))$ into $\Mod_{\n{K}_{\sol}}(\n{K}_{\sol}[G])$ is fully faithful. Then, for any $V\in \Mod_{\n{K}_{\sol}}(\n{D}^{la}(G,K))$ one has
\begin{gather*}
R\iHom_{\n{K}_{\sol}[G]}(K, C^{la}(G,V)_{\star_{1,3}}) \\ = R\iHom_{\n{D}^{la}(G,K)}(K, R\iHom_{\n{K}_{\sol}[G]}(\n{D}^{la}(G,K), C^{la}(G,V)_{\star_{1,3}}))  \\ 
  = R\iHom_{\n{D}^{la}(G,K)}(K, R\iHom_{\n{D}^{la}(G,K)}( \n{D}^{la}(G,K)\otimes^L_{\n{K}_{\sol}[G]}\n{D}^{la}(G,K), C^{la}(G,V)_{\star_{1,3}}))  \\
 = R\iHom_{\n{D}^{la}(G,K)}(K, C^{la}(G,V)_{\star_{1,3}}), \end{gather*}
proving  $\bb{Q}_p$-analytic vectors of $V$ seen as $\n{K}_{\sol}[G]$ or  $\n{D}^{la}(G,K)$-module agree.

\item  If $G$ is defined over $L\neq \bb{Q}_p$ and $V$ is a $\n{D}^{la}(G,K)$-module,   then we have a natural isomorphism 
\[
(V|_{\n{K}_{\sol}[G]})^{R\bb{Q}_p-la} \cong  V^{Rla}
\]
between the locally analytic  vectors of $V$ seen as $\n{D}^{la}(G,K)$-module, and the $\bb{Q}_p$-analytic vectors of $V$ seen as a  $\n{K}_{\sol}[G]$-module. Indeed, denote $\f{k}:=\ker(\f{g}\otimes_{\bb{Q}_p} L \to \f{g})$,  we have that  
\[
\begin{aligned}
C^{la}(G,K) & = R\iHom_{\n{D}^{la}(\f{k},K)}(K, C^{la}(G_{\bb{Q}_p}, K))   \\
&  = R\iHom_{\n{D}^{la}(G_{\bb{Q}_p},K)}(\n{D}^{la}(G,K) ,C^{la}(G_{\bb{Q}_p},K)),
\end{aligned}
\]
where in the second equivalence we use adjunction and the identity   $ \n{D}^{la}(G,K)=\n{D}^{la}(G_{\bb{Q}_p},K)\otimes_{\n{D}^{la}(\f{k}, K)}^L K $.

Then, we deduce that 
\begin{align*}
V^{Rla} &= R\iHom_{\n{D}^{la}(G,K)}(K, V\otimes_{\n{K}_{\sol}} C^{la}(G,K)) \\  
&=  R\iHom_{\n{D}^{la}(G,K)}(K, V\otimes_{\n{K}_{\sol}} R\iHom_{\n{D}^{la}(\f{k}, K)}(K, C^{la}(G_{\bb{Q}_p},K))) \\ 
&= R\iHom_{\n{D}^{la}(G,K)}(K,  R\iHom_{\n{D}^{la}(\f{k}, K)}(K,V\otimes_{\n{K}_{\sol}} C^{la}(G_{\bb{Q}_p},K))) \\
&= R\iHom_{\n{D}^{la}(G,K)}(K , R\iHom_{\n{D}^{la}(G_{\bb{Q}_p},K)}(\n{D}^{la}(G,K)   ,V\otimes_{\n{K}_{\sol}}C^{la}(G_{\bb{Q}_p},K) ))  \\ 
&= R\Hom_{\n{D}^{la}(G_{\bb{Q}_p},K)}(K, V\otimes_{\n{K}_{\sol}} C^{la}(G_{\bb{Q}_p},K) ) \\
&= (V|_{\n{K}_{\sol}[G]})^{R\bb{Q}_p-la},
\end{align*}
where the first and last equivalences are definitions, the second equivalence follows from the expresion  of $C^{la}(G,K)$ as $\f{k}$-Lie algebra cohomology of $C^{la}(G_{\bb{Q}_p},K)$ above,  the third follows form the fact that $K$ is a projective $\n{D}^{la}(\f{k},K)$-module and that its action on $V$ is trivial, the fourth equation follows from an adjunction and the formula $ \n{D}^{la}(G,K)=\n{D}^{la}(G_{\bb{Q}_p},K)\otimes_{\n{D}^{la}(\f{k}, K)}^L K $, and the second to last equivalence from another $\Hom-\otimes$ adjunction.   However, the   $L$-analytic vectors  of $V$ considered as a solid $G$-representation are given by
\[
(V|_{\n{K}_{\sol}[G]})^{RL-an} = R\Gamma(\f{k}, V^{R\bb{Q}_p-la}) = R\Gamma(\f{k}, V^{Rla}) \cong R\Gamma(\f{k}, K )\otimes_{\n{K}_{\sol}} V^{Rla}
\]
where in the last equivalence we use the fact that $\f{k}$ acts trivially on $V^{Rla}$, being a $\n{D}^{la}(G,K)$-module.
\end{enumerate} 
\end{remark}

Let us prove some basic properties of the functor of locally analytic vectors.

\begin{proposition}
\label{PropositionFirstPropertiesLocAn}
    The following assertions hold.
    \begin{enumerate}
        \item Let $G_0\subset G$ be any open subgroup and $V\in \Mod_{\n{K}_{\sol}}(\n{D}^{la}(G,K))$,  there is a natural equivalence $V^{Rla}|_{G_0}= (V|_{G_0})^{Rla}$ between the  restriction to $G_0$ of the  $G$-locally analytic vectors of $V$ and the $G_0$-locally analytic vectors of $V|_{G_0}$.
        
        \item The functor $(-)^{Rla}: \Mod_{\n{K}_{\sol}}(\n{D}^{la}(G,K))\to \Mod_{\n{K}_{\sol}}(\n{D}^{la}(G,K))$ is the right derived  functor of $W\mapsto W^{la}$ on the abelian category $\Mod^{\heartsuit}_{\n{K}_{\sol}}(\n{D}^{la}(G,K))$. 
        
        \item The functor $(-)^{Rla}: \Mod_{\n{K}_{\sol}}(\n{D}^{la}(G,K))\to \Mod_{\n{K}_{\sol}}(\n{D}^{la}(G,K))$ preserves small colimits.  In particular, small colimits of locally analytic representations are locally analytic. The same holds for $(-)^{Rh-an}$ and $G$ compact.

        \item If $G$ is compact, then $V^{Rla}= \varinjlim_{h} V^{Rh-an}= \varinjlim_{h} V^{Rh^+-an}$. 
    \end{enumerate}
\end{proposition}

\begin{proof} \leavevmode
\begin{enumerate}
    
\item By construction one has that 
\[
C^{la}(G,V)= R\iHom_{\n{D}^{la}(G_0,K)}(\n{D}^{la}(G,K),  C^{la}(G_0,V))
\]
where the $\n{D}^{la}(G_0,K)$ acts by left multiplication on $\n{D}^{la}(G,K)$ and by the left regular action on $C^{la}(G_0,V)$. One finds that 
\begin{align*}
V^{Rla}& = R\iHom_{\n{D}^{la}(G,K)}( K, (C^{la}(G,V))_{\star_{1,3}}) \\ 
        & = R\iHom_{\n{D}^{la}(G,K)}(K, R\iHom_{\n{D}^{la}(G_0,K)}(\n{D}^{la}(G,K), C^{la}(G_0,V)_{\star_{1,3}})) \\ 
        & = R\iHom_{\n{D}^{la}(G_0,K)}(K, C^{la}(G_0,V)_{\star_{1,3}}) \\
        & = (V|_{G_0})^{Rla}. 
\end{align*}

\item By Lemma \ref{LemmaFunctorFunctionsAndAction} the functor $V\mapsto C^{la}(G,V)$ is exact in the abelian category of solid $\n{D}^{la}(G,K)$-modules.  Then, one has that 
\[
V^{Rla}= R\iHom_{\n{D}^{la}(G,K)}(K, C^{la}(G,V)_{\star{1,3}})
\]
is a derived $\iHom$-functor,  which implies that it is the right derived functor of the invariants $V^{la}= C^{la}(G,K)^{G_{\star_{1,3}}}$. 

\item  By (1), we can assume that $G$ is compact. By definition of $(-)^{Rla}$ and $(-)^{Rh-an}$, since  $ V\mapsto C^{la}(G, V) = C^{la}(G,K) \otimes_{\n{K}_\sol}^L V$ and $V\mapsto C(\bb{G}^{(h)},V)= C(\bb{G}^{(h)},K)_{\sol}\otimes_{\n{K}_{\sol}}^L V$ commute with colimits, it suffices to show that $K$ is compact as $\n{D}^{la}(G,K)$-module, this follows from Theorem \ref{TheoremIdempotentPefect}.

\item Since taking locally analytic vectors commutes with colimits by (3), this is a consequence of the compactness of $K$ as a $\n{D}^{la}(G,K)$-module and Remark \ref{RemarkClaColimit}.        
\end{enumerate}
\end{proof}

The following proposition relates the functor of  analytic vectors  with the distribution algebras. 

\begin{proposition}
 \label{PropMainTheo1}
Let $G$ be a compact $p$-adic Lie group over $L$, and let $V \in \Mod_{\n{K}_{\sol}}(\n{D}^{la}(G,K))$. Then
\begin{gather*}
V^{Rh-an} = R\iHom_{\n{D}^{la}(G,K)}( \Dist(\bb{G}^{(h)},K), V)  \\ 
V^{Rh^+-an} = R\iHom_{\n{D}^{la}(G,K)}( \Dist^h(G,K), V). 
\end{gather*}
In particular, an object $V\in \Mod_{\n{K}_{\sol}}(\n{D}^{la}(G,K))$ is $h^+$-analytic if and only if it is a module over the idempotent $\n{D}^{la}(G,K)$-algebra $\n{D}^{h}(G,K)$, $h^+$-analytic implies $(h')^+$-analytic for any $h' > h$.
\end{proposition}

\begin{proof}
This follows from the same proof of Theorem 4.36 of \cite{SolidLocAnRR} using Corollary 2.19 of \textit{loc. cit.}
\end{proof}

For general groups, we have the following immediate consequence.

\begin{corollary}
\label{CoroLocAnColimitDmod}
Let $G$ be a $p$-adic Lie group over $L$ and let $V \in \Mod_{\n{K}_\sol}(\n{D}^{la}(G, K))$. Then, for any open compact subgroup $G_0$ of $G$, there is an equivalence of $G_0$-representations
\[ V^{Rla} = \varinjlim_{h} R \iHom_{\n{D}^{la}(G_0, K)}(\n{D}^{h}(G_0, K), V) = \varinjlim_{h} R \iHom_{\n{D}^{la}(G_0, K)}(\n{D}(\bb{G}^{(h)}_0, K), V). \]
In particular, if $V$ is $h$-analytic then it is locally analytic.
\end{corollary}

The following result verifies that taking locally analytic vectors defines an idempotent functor.

\begin{prop}
\label{PropRlaIdempotentGroup}
Suppose that $G$ is compact.  Let $V\in \Mod_{\n{K}_{\sol}}(\n{D}^{la}(G,K))$,  then the natural maps
\[ (V^{Rla})^{Rh-an} \to V^{Rh-an} \]
\[ (V^{Rla})^{Rh^+-an} \to V^{Rh^+-an} \]
induced by the natural map $V^{Rla} \to V$ are equivalences. In particular, for any group $G$, the natural map 
\[ (V^{Rla})^{Rla} \to V^{Rla} \]
is an equivalence, and the locally analytic vectors of a $\n{D}^{la}(G,K)$-module give rise to a solid locally analytic representation of $G$. 
\end{prop}
\begin{proof}
By Proposition \ref{PropositionFirstPropertiesLocAn} (1), we can assume that $G$ is compact in all the statements.  It suffices to prove that $(V^{Rla})^{Rh-an} \to V^{Rh-an}$ is an equivalence, as the other cases follow from this after taking limits or colimits.

    \begin{align*}
        (V^{Rla})^{Rh-an} & = \varinjlim_{h_1} (V^{R{h_1}^{+}-an})^{Rh-an} \\ 
        & = \varinjlim_{h_1}  R \iHom_{\n{D}^{la}(G, K)}( \n{D}(\bb{G}^{(h)}, K),  R \iHom_{\n{D}^{la}(G,K)}(\n{D}^{h_1}(G, K), V)) \\ 
        & = \varinjlim_{h_1} R\iHom_{\n{D}^{la}(G,K)} (\n{D}^{h_1}(G,K)\otimes^L_{\n{D}^{la}(G,K)} \n{D}(\bb{G}^{(h)},K), V) \\ 
        & = \varinjlim_{h_1} R\iHom_{\n{D}^{la}(G,K)}(\n{D}(\bb{G}^{(h)}, K), V )\\  
        & = V^{Rh-an},
    \end{align*}
where the first equality follows from Proposition \ref{PropositionFirstPropertiesLocAn} (3), the second equality follows from Proposition \ref{PropMainTheo1}, the third equality is a $\otimes$-$\iHom$ adjunction, the fourth equality follows from the fact that $\n{D}^{h_1}(G,K)$ is an idempotent $\n{D}^{la}(G,K)$-algebra and that $\n{D}(\bb{G}^{(h)},K)$ is a $\n{D}^{h_1}(G,K)$-module for all $h_1$ big enough, and the last equality is   Proposition \ref{PropMainTheo1} again. 
\end{proof}

The following proposition provides a different way to compute locally analytic vectors as a relative tensor product of $\n{D}^{la}(G,K)$-modules. We start with a lemma.

\begin{lemma} \label{LemmatrivialD}
Let $G$ be a compact $p$-adic Lie group and let $\n{D}$ denote $\n{D}^{la}(\f{g}, K)$ or $\n{D}^{la}(G, K)$. Then the choice of a Haar measure on $G$ produces an equivalence 
\begin{equation*}
R\iHom_{\n{D}}(K,\n{D})= K(\chi).
\end{equation*}
\end{lemma}

\begin{proof}
For $\n{D}=\n{D}^{la}(\f{g}, K)$ this follows by an explicit computation using the Koszul resolution of Lemma \ref{PropKoszulResolutionDist}. For $\n{D}=\n{D}^{la}(G,K)$ one argues as follows:  $K$ is a $\n{D}^{sm}(G,K)$-module and $\n{D}^{sm}(G,K)= K \otimes^L_{\n{D}^{la}(\f{g},K)} \n{D}^{la}(G,K)$ by Lemma \ref{PropPropertiesSmoothDistributionAlgebrab}. Then
\begin{align*}
    R\iHom_{\n{D}^{la}(G,K)}(K,  \n{D}^{la}(G,K)) & = R\iHom_{\n{D}^{la}(G,K)}(\n{D}^{sm}(G,K) \otimes^{L}_{\n{D}^{sm}(G,K)} K,  \n{D}^{la}(G,K)) \\ 
    & = R\iHom_{\n{D}^{sm}(G,K)}(K, R\iHom_{\n{D}^{la}(G,K)}(\n{D}^{sm}(G,K), \n{D}^{la}(G,K))) \\ 
    & = R\iHom_{\n{D}^{sm}(G,K)}(K, R\iHom_{\n{D}^{la}(\f{g},K)}(K,\n{D}^{la}(G,K))) \\ 
    & = R\iHom_{\n{D}^{sm}(G,K)}(K, K(\chi)\otimes_{\n{D}^{la}(\f{g},K)}^L \n{D}^{la}(G,K)) \\ 
    & = R\iHom_{\n{D}^{sm}(G,K)}(K,  K(\chi) \otimes_{\n{K}_\sol} \n{D}^{sm}(G,K)) \\
    & = \iHom_{\n{D}^{sm}(G,K)}(K,  K(\chi) \otimes_{\n{K}_\sol} \n{D}^{sm}(G,K)) \\ 
    & = K(\chi),
\end{align*}
where the first two equivalences are obvious, the third one follows from Lemma \ref{PropPropertiesSmoothDistributionAlgebrab}, the fourth one follows from the case of of $\n{D}^{la}(\f{g}, K)$, the fifth one follows again using Lemma \ref{PropPropertiesSmoothDistributionAlgebrab}, the sixth one follows from  Lemma \ref{CoroSmoothAlgebraMatrix} (1), and the last one follows from the fact that the trivial representation has multiplicity one in $\n{D}^{sm}(G,K)$ (given by the Haar measure of the  compact group).
\end{proof}

\begin{remark}
Lemma \ref{LemmatrivialD} is a special case of our cohomological comparison isomorphisms that will be shown in \S \ref{s:AdjunctionsCoho}.
\end{remark}

\begin{prop}
\label{PropCohomolocyAndHomology}
    Let $G$ be a compact $p$-adic Lie group. The following assertions hold.
    \begin{enumerate}
\item  Let  $V,W\in \Mod_{\n{K}_{\sol}}(\n{K}_{\sol}[G])$. Let $V\otimes_{\n{K}_\sol}^L W$ be endowed with the diagonal action. Then there is a natural equivalence
\[
R\iHom_{\n{K}_{\sol}[G]}(K, V\otimes^L_{\n{K}_\sol} W ) = ( K(\chi_{\bb{Q}_p}) \otimes_{\n{K}_{\sol}}^L \iota(V)) \otimes_{\n{K}_{\sol}[G]}^L W [-d]
\]
where $\iota(V)$ is the right $G$-module induced by $V$ under the natural involution $\iota: \n{K}_{\sol}[G] \to \n{K}_{\sol}[G]$,  $\chi_{\bb{Q}_p}= \det(\f{g}_{\bb{Q}_p})^{-1}$, and $d=\dim_{\bb{Q}_p} G$. 

\item Let $\n{D}$ denote $\n{D}^{la}(G,K)$ or $\n{D}^{la}(\f{g},K)$. Let $V,W \in \Mod_{\n{K}_{\sol}}(\n{D})$, then there is a natural equivalence  
 \[
    R\iHom_{\n{D}}(K, V\otimes_{\n{K}_{\sol}}^L W) = (K(\chi)\otimes^L_{\n{K}_\sol} \iota(V)) \otimes^L_{\n{D}} W [-d]
    \]
    where $\iota(V)$ is the right $\n{D}$-module obtained by the involution of $\n{D}$,  $\chi=\det(\f{g})^{-1}$ and $d=\dim_L G$.

    \end{enumerate}
\end{prop}

\begin{proof}
Without loss of generality we can take $K$ to be a finite extension of $\bb{Q}_p$, the general case is deduced by  base change.  By Theorem 5.19 of \cite{SolidLocAnRR} one has that 
\[
R\iHom_{\n{K}_{\sol}[G]}(K, V\otimes^L_{\n{K}_{\sol}} W ) = K(\chi_{\bb{Q}_p})\otimes_{\n{K}_{\sol}[G]}^L ( V\otimes_{\n{K}_{\sol}}^{L} W)[-d].
\]
where we see $K(\chi_{\bb{Q}_p})$ as a right representation. By Proposition \ref{PropositionHopfAlgebras} (4), we have natural equivalences
\begin{align*}
K(\chi_{\bb{Q}_p})\otimes_{\n{K}_{\sol}[G]}^L ( V\otimes_{\n{K}_{\sol}}^{L} W)[-d] &= 1 \otimes_{\n{K}_{\sol}[G]}^L ( \iota(K(\chi_{\bb{Q}_p}))\otimes_{\n{K}_\sol}^L V \otimes_{\n{K}_\sol}^L W) [-d] \\
&= (K(\chi_{\bb{Q}_p}) \otimes_{\n{K}_\sol}^L \iota(V)) \otimes^L_{\n{K}_{\sol}[G]} W[-d], 
\end{align*}
this shows (1). By Theorem \ref{TheoremIdempotentPefect},  the trivial representation is a perfect $\n{D}$-module, in particular dualizable,  this implies that the natural functor 
\[
 R\iHom_{\n{D}}(K,\n{D})\otimes^L_{\n{D}} W \to R\iHom_{\n{D}}(K, W)
\]
for any $W\in \Mod_{\n{K}_{\sol}}(\n{D})$ is an equivalence.  Then, for $V,W\in \Mod_{\n{K}_{\sol}}(\n{D})$,  by Proposition \ref{PropositionHopfAlgebras} (4) we have natural equivalences
\[
R\iHom_{\n{D}}(K, V\otimes^L_{\n{K}} W )= R\iHom_{\n{D}}(K,\n{D})\otimes^L_{\n{D}} (V\otimes^L_{\n{K}_{\sol}} W) = (R\iHom_{\n{D}}(K, \n{D}) \otimes_{\n{K}_{\sol}}^L \iota (V)) \otimes_{\n{D}}^L W.
\]
The result follows from Lemma \ref{LemmatrivialD}.
\end{proof}

\begin{remark}\label{RemarkActionsContragradientwas}
    In Proposition \ref{PropCohomolocyAndHomology} we see $\chi$ as a right $\n{D}^{la}(G,K)$-module. It arises as the determinant of the right action of $G$ on $\f{g}^{\vee}$ given by 
    \[
   ( H\cdot g)(v)= H(gvg^{-1})
    \]
    for $H\in \f{g}^{\vee}$, $v\in \f{g}$ and $g\in G$. We will often consider $\chi$ as a left representation as well, in this case, it arises as the determinant of the contragradient representation $\f{g}^{\vee}$ with action 
    \[
    (g\cdot H)(v)=H(g^{-1}vg). 
    \]
\end{remark}

\begin{corollary}
\label{CorollaryFiniteCohoDim}
Let $G$ be an arbitrary $p$-adic Lie group over $L$ of dimension $d$. The following assertions hold.
\begin{enumerate}
    \item  Let $V\in \Mod_{\n{K}_{\sol}}(\n{D}^{la}(G,K))$, then for any open compact sugroup $G_0\subset G$ one has an equivalence of $\n{D}^{la}(G_0, K)$-modules
    \[
    V^{Rla}=( \iota(C^{la}(G_0,K)_{\star_{1}})\otimes_{\n{K}_\sol}  K(\chi)[-d] ) \otimes_{\n{D}^{la}(G_0, K)}^L V. 
    \]
    In particular, the functor $(-)^{Rla}$ has cohomological dimension $d$. 

    \item Suppose that $G$  is defined over $\bb{Q}_p$ and let $V\in \Mod_{\n{K}_{\sol}}(\n{K}_{\sol}[G])$.  Then  for any open compact sugroup $G_0\subset G$ one has an equivalence of $\n{D}^{la}(G_0, K)$-modules
    \[
    V^{Rla}=( \iota(C^{la}(G_0,K)_{\star_{1}})\otimes_{\n{K}_\sol}  K(\chi)[-d] ) \otimes_{\n{K}_{\sol}[G_0]}^L V
    \]
where the locally analytic vectors are as in Remark \ref{RemarkLocAnvectorsSolidRep}. 
    \item Let $V,W\in \Mod_{\n{K}_{\sol}}(\n{D}^{la}(G,K))$. Then $V^{Rla} \otimes_{\n{K}_\sol} W^{Rla}$ is a locally analytic representation and the natural map
    \[
    V^{Rla}\otimes_{\n{K}_{\sol}} W^{Rla} \to   
    (V\otimes_{\n{K}_{\sol}}^L W^{Rla})^{Rla}
    \]
is an equivalence of $\n{D}^{la}(G, K)$-modules.
\end{enumerate}
\end{corollary}
\begin{proof}
    \begin{enumerate}
        \item  By  Proposition \ref{PropositionFirstPropertiesLocAn} (1) the locally analytic vectors are independent of $G_0\subset G$ compact open, so we can assume without loss of generality that $G$ is compact. Then, by Proposition \ref{PropCohomolocyAndHomology} (2) one has 
        \begin{align*}
            V^{Rla}& = R\iHom_{\n{D}^{la}(G, K)}(K, C^{la}(G,K)_{\star_{1}}\otimes_{\n{K}_{\sol}}^L V) \\ 
            & =  (\iota( C^{la}(G,K)_{\star_1})\otimes_{\n{K}_\sol} K(\chi)[-d] ) \otimes^L_{\n{D}^{la}(G,K)} V. 
        \end{align*} 

        \item This follows from the same argument of the previous point using Proposition \ref{PropCohomolocyAndHomology} (1) instead. 

        \item Again by Proposition \ref{PropositionFirstPropertiesLocAn} (1) we can assume that $G$ is compact.  The fact that $V^{Rla} \otimes_{\n{K}_\sol} W^{Rla}$ is locally analytic follows by observing that $V^{Rla} \otimes_{\n{K}_\sol} W^{Rla} = \varinjlim_{h > 0} V^{Rh^+ -an} \otimes_{\n{K}_\sol} W^{Rh^+-an}$, that $V^{Rh^+ -an} \otimes_{\n{K}_\sol} W^{Rh^+-an}$ is a $\n{D}^h(G, K)$-module, hence $h^+$-analytic and hence locally analytic by Corollary \ref{CoroLocAnColimitDmod}, and that colimits of locally analytic representations are locally analytic by Proposition \ref{PropositionFirstPropertiesLocAn} (3).

We now prove the final equivalence. The orbit map $\s{O}_{W}: W^{Rla}\to C^{la}(G,K)\otimes_{\n{K}_{\sol}}^L W^{Rla}$ induces a natural equivalence 
    \begin{equation}
    \label{eqIsoChangeVariables}
    C^{la}(G,W^{Rla})_{\star_1} \xrightarrow{\sim } C^{la}(G,W^{Rla})_{\star_{1,3}},
    \end{equation}
    at the level of functions this maps sends $f:G\to W$ to the function $\widetilde{f}:G \to W$ given by $\widetilde{f}(g)= g\cdot f(g)$. 
    Then, one computes 
    \begin{align*}
    (V\otimes_{\n{K}_\sol}^L  W^{Rla})^{Rla} & = ( \iota(C^{la}(G,K)_{\star_1})\otimes_{\n{K}_\sol} K(\chi)[-d]) \otimes_{\n{D}^{la}(G,K)}^L (V\otimes_{\n{K}_{\sol}}^L W^{Rla}) \\ 
    & = ( \iota(C^{la}(G,K)_{\star_1} \otimes_{\n{K}_{\sol}}^L W^{Rla})\otimes_{\n{K}_\sol} K(\chi)[-d]) \otimes_{\n{D}^{la}(G,K)}^L V \\
    & = (\iota ( C^{la}(G,W^{Rla})_{\star_{1,3}}) \otimes_{\n{K}_\sol} K(\chi)[-d]) \otimes^L_{\n{D}^{la}(G,K)} V \\ 
    & = (\iota( C^{la}(G,W^{Rla})_{\star_{1}})\otimes_{\n{K}_\sol} K(\chi)[-d]) \otimes^L_{\n{D}^{la}(G,K)} V \\ 
    & = \bigg((\iota(C^{la}(G,K)_{\star_{1}})\otimes_{\n{K}_\sol} K(\chi)[-d] ) \otimes^L_{\n{D}^{la}(G,K)} V \bigg) \otimes^L_{\n{K}_{\sol}}W^{Rla} \\ 
    & = V^{Rla}\otimes_{\n{K}_{\sol}}^L W^{Rla}.
    \end{align*}
    In the first equality we use part (1). In the second equality we move $W$ to the left part of the tensor using Proposition  \ref{PropositionHopfAlgebras} (4). The third equality is the definition  $C^{la}(G,W)= C^{la}(G,K)\otimes^L_{\n{K}_{\sol}} W$. The fourth equality  uses  the natural equivalence \eqref{eqIsoChangeVariables}. In the  fifth equality we take out the tensor with $W^{Rla}$ since $\n{D}^{la}(G,K)$ is acting trivially on it. In the last equality we use part (1) again. 
    \end{enumerate}
\end{proof}

The previous computation implies that there are  representations with higher locally analytic vectors. 

\begin{corollary}
\label{CoroLocAnDla}
Let $G$ be a $p$-adic Lie group  over $L$ of dimension $d$. Then for any profinite set $S$ we have 
\[
(\n{D}^{la}(G,K)\otimes_{\n{K}_{\sol}} \n{K}_{\sol}[S] )^{Rla}= (C^{la}_{c}(G,K) \otimes_{\n{K}_\sol} K(\chi) [-d])\otimes_{\n{K}_{\sol}}\n{K}_{\sol}[S] 
\]
where $C^{la}_c(G,K)= \n{D}^{la}(G,K) \otimes_{\n{D}^{la}(G_0,K)} C^{la}(G_0,K)$ is the space of compactly supported locally analytic functions of $G$. If $G$ is defined over $\bb{Q}_p$ we also have 
\[
(\n{K}_{\sol}[G\times S])^{Rla}= (C^{la}_{c}(G,K) \otimes_{\n{K}_\sol} K(\chi) [-d])\otimes_{\n{K}_{\sol}}\n{K}_{\sol}[S].
\]
\end{corollary}
\begin{proof}
 By Corollary \ref{CorollaryFiniteCohoDim} (1) we have that 
\begin{gather*}
     (\n{D}^{la}(G,K)\otimes_{\n{K}_{\sol}} \n{K}_{\sol}[S] )^{Rla} \\
     = (\iota(C^{la}(G_0,K))_{\star_{1}}\otimes_{\n{K}_\sol} K(\chi)[-d]) \otimes_{\n{D}^{la}(G_0,K)} (\n{D}^{la}(G,K)\otimes_{\n{K}_{\sol}} \n{K}_{\sol}[S]) \\
      =\bigg(  (\iota(C^{la}(G_0,K))_{\star_{1}}\otimes_{\n{K}_\sol} K(\chi)[-d]) \otimes_{\n{D}^{la}(G_0,K)} \n{D}^{la}(G,K)  \bigg)  \otimes_{\n{K}_{\sol}} \n{K}_{\sol}[S] \\
     = (C^{la}_{c}(G,K)\otimes_{\n{K}_\sol}  K(\chi)[-d]) \otimes_{\n{K}_{\sol}} \n{K}_{\sol}[S].
    \end{gather*}
The second claim follows by the same argument using Corollary \ref{CorollaryFiniteCohoDim} (2) instead. 
\end{proof}

\subsection{The category of locally analytic representations}
\label{ss:SLACategoriesofRepresentations}

Let $L$ be a finite extension of $\bb{Q}_p$. Our next goal is to define the $\infty$-category of locally analytic representations and discuss some general properties of it.

\begin{definition}
    We define    the $\infty$-category of locally analytic representations, denoted as $\Rep^{la}_{\n{K}_{\sol}}(G)$,   to be the full subcategory of $\Mod_{\n{K}_{\sol}}(\n{D}^{la}(G,K))$ whose objects are locally analytic representations of $G$. In other words, $\Rep^{la}_{\n{K}_{\sol}}(G)$ is the full subcategory of solid $\n{D}^{la}(G,K)$-modules whose objects are the $V$ such that $V^{Rla}=V$. 
\end{definition}

 Our next task is to show that the derived category of locally analytic representations has a natural $t$-structure and that it is the derived category of its heart.

\begin{lemma}
\label{LemmaLocAnBounded}
Given $V\in \Mod_{\n{K}_{\sol}}(\n{D}^{la}(G,K))$, one has that  
\[
V^{Rla} = \varinjlim_{b\to +\infty} \varprojlim_{a\to -\infty} (\tau^{[a,b]} V)^{Rla},
\]
in the homotopy category, where $a,b\in \bb{Z}$ with $a\leq b$ and $\tau^{[a,b]}= \tau^{\geq a}\circ \tau^{\leq b} $ is the canonical truncation  in the interval $[a,b]$ in cohomological notation. 
\end{lemma}
\begin{proof}
This follows from the fact that $(-)^{Rla}$ has finite cohomological dimension, see Corollary \ref{CorollaryFiniteCohoDim} (1). 
\end{proof}

We now prove some basic and fundamental properties of the category of solid locally analytic representations.

\begin{prop}
\label{PropositionStableColimitsAndTensor}
$\Rep^{la}_{\n{K}_\sol}(G)$ is stable under all  small colimits of $\Mod_{\n{K}_{\sol}}(\n{D}^{la}(G, K))$ and tensor products over $\n{K}_{\sol}$.
\end{prop}
\begin{proof}
This follows from the fact that taking locally analytic vectors preserves colimits, c.f. Proposition \ref{PropositionFirstPropertiesLocAn}, and the projection formula of  Corollary \ref{CorollaryFiniteCohoDim} (3).  
\end{proof}

\begin{lemma}
\label{LemmaLocAnColimitDh}
Let $G$ be compact. Then  $\Rep^{la}_{\n{K}_\sol}(G)$ is the full subcategory of $\Mod_{\n{K}_{\sol}}(\n{D}^{la}(G, K))$ stable under  all small colimits containing the categories $\Mod_{\n{K}_\sol}(\n{D}^h(G, K))$ for all $h\geq 0$.
\end{lemma}

\begin{proof}
 By Proposition \ref{PropMainTheo1} the category  $\Rep^{la}_{\n{K}_\sol}(G)$ contains $\Mod_{\n{K}_\sol}(\n{D}^{h}(G, K))$ for any $h \geq 0$ and. Then, by Proposition \ref{PropositionStableColimitsAndTensor},   $\Rep^{la}_{\n{K}_\sol}(G)$ is stable under colimits, hence the category generated by  $\n{D}^h(G,K)$-modules under colimits is contained in $\Rep^{la}_{\n{K}_\sol}(G)$. Conversely,    by Corollary \ref{CoroLocAnColimitDmod}, for any $V \in \Mod_{\n{K}_\sol}(\n{D}^{la}(G, K))$, we have an equivalence
	\[ V^{Rla} = \varinjlim_{h} R \iHom_{\n{D}^{la}(G_0, K)}(\n{D}^{h}(G_0, K), V) \]
	in $\Mod_{\n{K}_\sol}(\n{D}^{la}(G, K))$, which shows the reverse inclusion and finishes the proof.
\end{proof}

\begin{prop} \leavevmode
\label{PropPropertiesCategoryLocAn}
An object $V\in \Mod_{\n{K}_{\sol}}(\n{D}^{la}(G,K))$ is locally analytic if and only if $H^i(V)$ is non-derived locally analytic for all $i\in \bb{Z}$.  In particular, the $t$-structure of $\Mod_{\n{K}_{\sol}}(\n{D}^{la}(G, K))$ induces a $t$-structure on $\Rep_{\n{K}_{\sol}}(G)$.
\end{prop}
\begin{proof}
We can assume without loss of generality that $G$ is compact. If $V$ is locally analytic then $V^{Rla}=\varinjlim_{h} V^{Rh^+-an}$ and $H^i(V)= \varinjlim_h H^i( V^{Rh^+-an})$, but $V^{Rh^+-an}$ is a $\n{D}^{h}(G,K)$-module. This shows that the cohomology groups $H^i(V)$ are colimits of $\n{D}^{h}(G,K)$-modules for $h\to \infty$ and are locally analytic by Lemma \ref{LemmaLocAnColimitDh} (so a fortriori non-derived locally analytic). Conversely,  suppose that $H^i(V)$ is non-derived locally analytic for all $i\in \bb{Z}$. We want to show that $V$ is locally analytic. By Lemma \ref{LemmaLocAnBounded} we can assume that $V$ is bounded with support in cohomological degrees $[0,k]$. By an inductive argument one  the length of the support of $V$, we can assume that $\tau^{\geq 1} V$ is locally analytic, then $V$ is an extension 
\[
H^0(V) \to V \to \tau^{\geq 1} V \to H^0(V)[1].
\] 
Since $H^0(V)$ is non-derived locally analytic, it can be written as a filtered colimit of $\n{D}^h(G,K)$-modules, then it is actually locally analytic by Lemma \ref{LemmaLocAnColimitDh}.  This exhibits $V$ as the fiber of $ \tau^{\geq 1} V \to H^0(V)[1]$  which is a locally analytic representation by Proposition \ref{PropositionStableColimitsAndTensor}. 
\end{proof}

\begin{prop} \label{PropReplaGrothendieck}
The category $\Rep^{la, \heartsuit}_{\n{K}_\sol}(G)$  is a Grothendieck abelian category.   Moreover  $\Rep^{la}_{\n{K}_\sol}(G)$ is the  $\infty$-derived category of $\Rep^{la, \heartsuit}_{\n{K}_\sol}(G)$.
\end{prop}

\begin{proof}
To show that $\Rep^{la, \heartsuit}_{\n{K}_\sol}(G)$ is a Grothendieck category, by the above results, it is enough to see that it has a small family of generators. Let $G_0\subset G$ be a compact open subgroup. Recall that, since we are working with $\kappa$-small condensed sets, the category $\Mod_{\n{K}_\sol}(\n{D}^{h}(G_0, K))$ admits a system of compact projective generators given by the modules  $\n{D}^h(G_0, K) \otimes_{\n{K}_\sol} \n{K}_{\sol}[S]$ for $S$ a $\kappa$-small profinite set. Hence, by Lemma \ref{LemmaLocAnColimitDh} the category $\Rep^{la}_{\n{K}_\sol}(G)$ is generated by $\{ \n{D}^{la}(G,K)\otimes_{\n{D}^{la}(G_0,K)} \n{D}^h(G_0, K) \otimes_{\n{K}_\sol} \n{K}_{\sol}[S]\}_{h,S}$ where $h>0$ and $S$ runs over all the $\kappa$-small profinite sets.

Let us first prove that the right adjoint of the fully faithful inclusion $\Rep^{la, \heartsuit}_{\n{K}_\sol}(G)\to \Mod^{\heartsuit}(\n{D}^{la}(G, K))$ of abelian categories is given by the  (non-derived) locally analytic vectors. Let $V\in \Rep^{la,\heartsuit}_{\n{K}_{\sol}}(G)$ and $W\in \Mod_{\n{K}_{\sol}}^{\heartsuit}(\n{D}^{la}(G,K))$, we want to prove that the natural map 
\[
\iHom_{\n{D}^{la}(G,K)}(V,W^{la}) \to \iHom_{\n{D}^{la}(G,K)}(V,W)
\]
is an equivalence.  Then, it suffices to take $V= \n{D}^{la}(G,K)\otimes_{\n{D}^{la}(G_0,K)} \n{D}^{h}(G_0,K) \otimes_{\n{K}_{\sol}} \n{K}_{\sol}[S]$   and show that the natural map
\[
R\iHom_{\n{D}^{la}(G,K)}(V,W^{Rla}) \to R\iHom_{\n{D}^{la}(G,K)}(V,W)
\]
is an equivalence. Indeed, one can find a resolution $P^\bullet \to V$ of $V$ where each term is a direct sum of elements in $\{ \n{D}^{la}(G,K)\otimes_{\n{D}^{la}(G_0,K)} \n{D}^h(G, K) \otimes_{\n{K}_\sol} K_{\sol}[S]\}_{h,S}$ and calculate $R\iHom_{\n{D}^{la}(G, K)}(V ,W)$ in terms of this resolution. Let $V=\n{D}^{la}(G,K)\otimes_{\n{D}^{la}(G_0,K)} \n{D}^h(G_0, K) \otimes_{\n{K}_\sol} K_{\sol}[S]$,  since we are taking internal $\iHom$, we can assume that $S=*$. By Proposition \ref{PropRlaIdempotentGroup}  we have that
\begin{align*}
    R\iHom_{\n{D}^{la}(G,K)}(V,W)& = R\iHom_{\n{D}^{la}(G_0,K)}(\n{D}^h(G_0,K), W) \\ 
    & = W^{Rh^+-an} \\ 
    & = (W^{Rla})^{Rh^+-an} \\ 
    & = R\iHom_{\n{D}^{la}(G,K)}(V,W^{Rla}),
\end{align*}
proving the claim.

Now, let $I$ be a $\kappa$-small injective $\n{D}^{la}(G, K)$-module. Since the functor of derived locally analytic vectors is a right derived functor by Proposition \ref{PropositionFirstPropertiesLocAn} (2), one has $I^{la}=I^{Rla}$. By  \cite[\href{https://stacks.math.columbia.edu/tag/015Z}{Tag 015Z}]{stacks-project}  the object $I^{la}$ is injective in  $\Rep^{la, \heartsuit}_{\n{K}_\sol}(G)$.

  Then, if $V\in \Rep^{la,\heartsuit}_{\n{K}_{\sol}}(G)$ and $I^{\bullet}$ is an injective resolution of $V$ as $\n{D}^{la}(G,K)$-module, we have 
\[
V= V^{Rla}= I^{\bullet, Rla}= I^{\bullet,la},
\]
so that $I^{\bullet,la}$ is an injective resolution of $V$ in $\Rep^{la,\heartsuit}_{\n{K}_{\sol}}(G)$. The previous implies that the $R\iHom$ in the derived category of $\Rep^{la, \heartsuit}_{\n{K}_\sol}(G)$ can be computed as the $R\iHom$ in $\Mod_{\n{K}_{\sol}}(\n{D}^{la}(G, K))$. Since $\Rep^{la}_{\n{K}_{\sol}}(G)$ is left complete by Lemma \ref{LemmaLocAnBounded},  one deduces that it is the derived $\infty$-category  of $\Rep^{la, \heartsuit}_{\n{K}_\sol}(G)$. 
\end{proof}

As a byproduct of the proof of Proposition \ref{PropReplaGrothendieck}, we have the following result.

\begin{corollary}
\label{CoroRightAdjointLocAnGroup}
The fully faithful inclusion $\Rep^{la}_{\n{K}_\sol}(G) \to \Mod_{\n{K}_\sol}(\n{D}^{la}(G, K))$ has for  right adjoint the functor of locally analytic vectors $V \mapsto V^{Rla}$. 
\end{corollary}

We end this section by  briefly discussing  some functorial properties of the categories of locally  analytic representations.  Let $H \to G$ be a morphism of $p$-adic Lie groups over $L$ and denote by $\f{h} \to \f{g}$ the corresponding map between their Lie algebras. We have a natural morphism of projective systems of good lattices $\{\n{M}\}_{\n{M}\subset \f{h}_{\overline{L}}} \to \{\n{L}\}_{\n{L}\subset \f{g}_{\overline{L}}}$. In particular, if $\n{M}$ maps to $\n{L}$, we have a morphism of completed enveloping algebras  $\widehat{U}(\n{M})\to \widehat{U}(\n{L})$. On the other hand, the forgetful functor $F: \Mod_{\n{K}_{\sol}}(\n{D}^{la}(G,K)) \to \Mod_{\n{K}_{\sol}}(\n{D}^{la}(H,K))$ restricts to a forgetful functor $\Rep^{la}_{\n{K}_{\sol}}(G) \to \Rep^{la}_{\n{K}_{\sol}}(H)$. It has a right adjoint which is given by the locally analytic induction
\[
F: \Rep^{la}_{\n{K}_\sol}(G) \rightleftharpoons \Rep^{la}_{\n{K}_\sol}(H): \mbox{la-}\Ind^{G}_H(-)
\]
where 
\[
\mbox{la-}\Ind_H^{G}(V) := R\iHom_{\n{D}^{la}(H,K)}(\n{D}^{la}(G,K), V)^{RG-la}.  
\]
If $H \subset G$ is an open subgroup, then the forgetful functor commutes with limits in the category of locally analytic representations (computed as the locally analytic vectors of the limit in $\Mod_{\n{K}_{\sol}}(\n{D}^{la}(G,K))$). Then it has a left adjoint called the  \textit{compactly supported induction} and is given by 
\[
\mbox{la-}\mathrm{cInd}_{H}^{G}(V)= \n{D}^{la}(G,K)\otimes^L_{\n{D}^{la}(H,K)} V. 
\]

\subsection{Detecting locally analyticity}

\label{ss:SomeAdditionalResults}

We finish this section with an additional result that can come in handy when proving that a  solid representation is locally analytic.  In the following we let $G$ be an uniform pro-$p$-group over $\bb{Q}_p$.  

\begin{prop}
\label{PropTestLocAn}
Let  $V\in \Mod_{\n{K}_{\sol},\geq 0}(\n{K}_{\sol}[G])$ a  connective solid $\n{K}_{\sol}[G]$-module. Suppose that the following holds: 
\begin{enumerate}
    \item There exists a $p$-adically complete object $V^+\in \Mod_{K^+_{\sol},\geq 0}( K^+_{\sol}[G])$ with $V^+\otimes_{K^+_{\sol}}^L K = V$. 
  
    \item The action of $G$ on $V^+/p$ factors through a finite quotient, i.e. there exists an open subgroup $G_0\subset G$ such that the restriction of $V^+/p$ to  $G_0$  belongs to the image of  $\Mod((K^+/p)_{\sol})$ into $\Mod_{K^+_{\sol}}(K^+_{\sol}[G_0])$ via the trivial representation.  
\end{enumerate}
Then $V$ is locally analytic.
\end{prop}

\begin{proof}
We can assume without loss of generality that $K=\bb{Q}_p$, that is the property of being locally analytic is independent of the base field. Let $h>0$, let $C^{h}(G,\bb{Q}_p)$ be the space of $h$-analytic functions of $G$, and let $C^{h,+}_G\subset C^{h}(G,\bb{Q}_p)$ be the subspace of power bounded functions. The cohomology  $R\Gamma(G, C^{h,+}_G \otimes^L_{\bb{Z}_{p},\sol} V^+) $ has a natural structure of $R\Gamma(G, C^{h,+}_G )$-module, that is, $G$-cohomology is right adjoint to the symmetric monoidal functor $\mathrm{Mod}_{\bb{Z}_{p,\sol}}\to \mathrm{Mod}_{\bb{Z}_{p,\sol}[G]}$ given by the trivial representation, and hence it is lax symmetric monoidal. We claim that the natural map (where the base is considered as an $\mathbb{E}_{\infty}$-algebra in $\Mod_{\bb{Z}_{p,\sol}}$)
\begin{equation}\label{eqojpamfpawfraw}
R\Gamma(G,  C^{h,+}_G\otimes^L_{\bb{Z}_{p},\sol} V^+)\otimes^L_{R\Gamma(G, C^{h,+}_G)} \bb{Z}_p\xrightarrow{\sim} V^+
\end{equation}
is an isomorphism. Suppose this holds, by inverting $p$ and taking colimits as $h\to \infty$, we have that $\varinjlim_h R\Gamma(G, C^{h,+}_G)[\frac{1}{p}]=\bb{Q}_p^{Rla}=\bb{Q}_p$ and $\varinjlim_h R\Gamma(G,  C^{h,+}_G\otimes^L_{\bb{Z}_{p},\sol} V^+)[\frac{1}{p}]=V^{Rla}$ so that 
\[
V^{Rla}=\varinjlim_h \bigg(R\Gamma(G,  C^{h,+}_G\otimes^L_{\bb{Z}_{p},\sol} V^+)\otimes^L_{R\Gamma(G, C^{h,+}_G)} \bb{Z}_p \bigg)[\frac{1}{p}] \xrightarrow{\sim} V
\]
is an isomorphism as wanted.

We finish by proving \eqref{eqojpamfpawfraw}. Since $V^+$ is connective, \cite[Proposition 2.12.10 (i)]{MannSix} implies that all the previous tensor products are derived $p$-complete. By derived Nakayama's lemma \cite[\href{https://stacks.math.columbia.edu/tag/0G1U}{Tag 0G1U}]{stacks-project}, it suffices to prove the claim modulo $p$. In that case, we have to prove that the natural map 
\[
R\Gamma(G,  (C^{h,+}_G/p)\otimes^L_{\bb{F}_{p},\sol} V^+/p)\otimes^L_{R\Gamma(G, C^{h,+}_G/p)} \bb{F}_p\xrightarrow{\sim} V^+/p
\]
is an isomorphism. Since $G$ is uniform pro-$p$, $\bb{F}_p$ admits a Lazard resolution as perfect $\bb{F}_{p,\sol}[G]$-modules (see \cite[Theorem 5.7]{SolidLocAnRR}), which implies the projection formula for $G$-cohomology, that is, if $W$ is a trivial $\bb{F}_{p,\sol}$-solid $G$-representation  and $M$ is an $\bb{F}_{p,\sol}$-linear $G$-representation,  then the natural map
\[
R\Gamma(G,M)\otimes^L_{\bb{F}_{p,\sol}} W \xrightarrow{\sim} R\Gamma(G, M\otimes^L_{\mathbb{F}_{p,\sol} }W)
\]
is an isomorphism. More precisely, this projection formula follows from the analogue formula of \cite[Theorem 5.19]{SolidLocAnRR} with $\bb{F}_p$-coefficients. 

Applying the previous to $V^+/p$, we formally deduce that 
\[
\begin{aligned}
R\Gamma(G,  (C^{h,+}_G/p)\otimes^L_{\bb{F}_{p},\sol} V^+/p)\otimes^L_{R\Gamma(G, C^{h,+}_G/p)} \bb{F}_p & = \bigg(R\Gamma(G,  C^{h,+}_G/p)\otimes^L_{\bb{F}_{p},\sol} V^+/p\bigg) \otimes^L_{R\Gamma(G, C^{h,+}_G/p)} \bb{F}_p \\ 
& = \bigg(R\Gamma(G,  C^{h,+}_G/p)  \otimes^L_{R\Gamma(G, C^{h,+}_G/p)} \bb{F}_p  \bigg) \otimes^L_{{\bb{F}_{p},\sol} }V^+/p \\ 
& = V^{+}/p,	
\end{aligned}
\]
proving what we wanted. 
\end{proof}

\begin{remark}
    The same proof of Proposition \ref{PropTestLocAn} holds for a quotient $V^+/p^{\epsilon}$ for any $\epsilon>0$, namely, it is enough to suppose that $V^+/p^{\epsilon}$ arises as a trivial $G_0$-representation.  
\end{remark}

\section{Geometric interpretation  of locally analytic representations}

Let $G$ be a  $p$-adic Lie group over a finite extension $L$ of $\bb{Q}_p$.  The purpose of this section is to give two different algebro-geometric identifications of the category of locally analytic representations, and to apply this to the $p$-adic Langlands correspondence.

For $G$ compact, in \S \ref{SectLocAnModDla} and \S \ref{SectLocAnModDla2}, we identify the category of locally analytic representations inside the category $\Mod_{\n{K}_\sol}(\n{D}^{la}(G, K))$. The algebra $\n{D}^{la}(G,K)$ can be thought of as the global sections of a non-commutative Stein space. Global  sections of sheaves over this space will give objects of $\Mod_{\n{K}_{\sol}}(\n{D}^{la}(G, K))$,  and we will prove that the functor of ``global sections with compact support'' induces an equivalence of stable $\infty$-categories between quasi-coherent sheaves of this space and  $\Rep^{la}_{\n{K}_\sol}(G)$. 

In \S \ref{SectLocAnCoModCla}, we give a second interpretation, for general $G$, of the category of solid locally analytic representations of $G$ as the derived category of comodules of the coalgebra $C^{la}(G,K)$ of $L$-analytic functions. Heuristically, if $G^{la}$ denotes the ``analytic spectrum of $C^{la}(G,K)$'', the previous description provides a natural equivalence  between $\Rep^{la}_{\n{K}_{\sol}}(G)$ and solid quasi-coherent sheaves of the classifying stack $[*/ G^{la}]$. 

Finally, in \ref{SectpadicLL} we conclude with an application to the categorical locally analytic $p$-adic Langlands correspondence for $\GL_1$.

\subsection{Locally analytic representations as quasi-coherent $\n{D}^{la}(G,K)$-modules} \label{SectLocAnModDla}

In the following we let $G$ be a compact $p$-adic Lie group over $L$.

\begin{definition}
Let us write $\n{D}^{la}(G,K)= \varprojlim_{h\to \infty} \n{D}^{h}(G,K)$ as a limit of $h$-analytic distribution algebras.  We define the category $\Mod^{qc}_{\n{K}_\sol}(\D^{la}(G, K))$ of solid quasi-coherent modules over $\n{D}^{la}(G,K)$ as the $\infty$-category
\[ \Mod^{qc}_{\n{K}_\sol}(\n{D}^{la}(G, K)) := \varprojlim_{h > 0} \Mod_{\n{K}_\sol}(\n{D}^{h}(G, K)), \]
where the transition maps in the limit are given by base change. 
\end{definition}

Objects in the category $\n{C}=\Mod_{\n{K}_{\sol}}^{qc}(\n{D}^{la}(G,K))$ are sequences of modules $(V_h)_{h}$ with $V_h\in \Mod_{\n{K}_{\sol}}(\n{D}^{h}(G,K))$, and  for $h'>h$  the datum of an isomorphism $\n{D}^{h}(G,K)\otimes^L_{\n{D}^{h'}(G,K)} V_{h'}\xrightarrow{\sim}V_{h}$, compatible with the $h$'s and higher coherences. Given two objects $(V_h)_{h}$ and $(W_h)_{h}$ in $\Mod_{\n{K}_{\sol}}^{qc}(\n{D}^{la}(G,K))$, the spectrum of morphisms is given by 
\[
R\Hom_{\n{C}}((V_h)_{h}, (W_{h})_{h})= \varprojlim_{h} R\Hom_{\n{D}^{h}(G,K)}(V_h,W_h).
\]
The following lemma will give a sufficient condition for a morphism of objects in $\n{C}$ to be an equivalence. 

\begin{lemma}
\label{LemmaCofinalDist}
    Let $(R_n)_{n\in \bb{N}}$ be a limit sequence of $\bb{E}_1$-$\n{K}_{\sol}$-algebras and let $\n{C}= \varprojlim_{n} \Mod_{\n{K}_{\sol}}(R_n)$ be the limit category along base change. Let $f_{\bullet}:(X_n)_{n}\to (Y_{n})_{n}$ be a morphism of objects in $\n{C}$, and suppose that there are arrows $h_{n+1}:Y_{n+1}\to X_{n}$ of $R_{n+1}$-modules making the following diagram commute:
    \[
    \begin{tikzcd}
        X_{n+1} \ar[r] \ar[d, "f_{n+1}"] & X_{n}  \ar[d, "f_{n}"] \\ 
        Y_{n+1} \ar[r] \ar[ur, "h_{n+1}"'] & Y_{n}.
    \end{tikzcd}
    \]
    Then $f_{\bullet}$ is an equivalence in $\n{C}$. 
\end{lemma}
\begin{proof}
    We have to prove that each $f_{n+1}: X_{n+1}\to Y_{n+1}$ is an equivalence. We have a commutative diagram by extension of scalars
    \[
    \begin{tikzcd}
        X_{n+1} \ar[r]  \ar[d, "f_{n+1}"]&  R_n \otimes^L_{R_{n+1}} X_{n+1} \ar[d, "1\otimes f_{n+1}"']  \ar[r, "\sim"]& X_{n}  \ar[d, "f_{n}"]\\ 
        Y_{n+1} \ar[r]  & R_n\otimes^L_{R_{n+1}} Y_{n+1}  \ar[r, "\sim"] \ar[ur, "1\otimes h_{n+1}"'] & Y_{n}.
    \end{tikzcd}
    \]
  A diagram chase shows that the map $Y_n\xrightarrow{\sim} R_{n}\otimes_{R_{n+1}}^L Y_{n+1} \xrightarrow{1\otimes h_{n+1}} X_n$ defines a homotopy inverse of $f_n$ proving that $f_{\bullet}$ is an equivalence. 
\end{proof}

\begin{remark}\label{LemmaCofinalDistApplication}
Keep the notation of Lemma  \ref{LemmaCofinalDist}, and denote $R:=\varprojlim_{n} R_{n}$. Suppose that all the maps $R\to R_n$ are idempotent making $\mathrm{Mod}_{\n{K}_{\sol}}(R_n)\subset \mathrm{Mod}_{\n{K}_{\sol}}(R)$ a full subcategory along the forgetful functor. Then, Lemma \ref{LemmaCofinalDist} implies that if $f_{\bullet}\colon (X_n)\to (Y_n)$ is a morphism in $\n{C}=\varprojlim_{n} \mathrm{Mod}_{\n{K}_{\sol}}(R_n)$ such that $f_{\bullet}$ is a pro-isomorphism of objects in $\mathrm{Mod}_{\n{K}_{\sol}}(R)$, then it is an isomorphism in the limit category $\n{C}$. 
\end{remark}

Next, we define natural functors between the category of modules over $\n{D}^{la}(G,K)$ and $\Mod^{qc}_{\n{K}_{\sol}}(\n{D}^{la}(G,K))$. 

\begin{lemma} \label{LemmaAdjunction}
    Let $j^*: \Mod_{\n{K}_{\sol}}(\n{D}^{la}(G,K))\to \Mod^{qc}_{\n{K}_{\sol}}(\n{D}^{la}(G,K))$ be the localization functor sending a $\n{D}^{la}(G,K)$-module $V$ to the sequence $(V_h)_{h}$ with $V_h=\n{D}^{h}(G,K) \otimes_{\n{D}^{la}(G,K)}^L V$. Then $j^*$ has a right adjoint $j_*$  given by 
    \[
    j_*(V_{h})_{h} = R\varprojlim_{h} V_{h}. 
    \]
\end{lemma}
\begin{proof}
    Let us denote $\n{C}=\Mod_{\n{K}_{\sol}}^{qc}(\n{D}^{la}(G,K))$, let $V=(V_h)_h \in \n{C}$ and $W\in \Mod_{\n{K}_{\sol}}(\n{D}^{la}(G,K))$. We have a natural map $W\to R\varprojlim_{h}(\n{D}^{h}(G,K)\otimes^L_{\n{D}^{la}(G,K)} W)$, and by construction we have 
    \begin{align*}
    R\Hom_{\n{C}}(j^* W, V) &= R\varprojlim_{h} R\Hom_{\n{D}^{h}(G,K)}(\n{D}^{h}(G,K)\otimes^L_{\n{D}^{la}(G,K)} W, V_h) \\              & = R\varprojlim_{h} R\Hom_{\n{D}^{la}(G,K)}(W, V_h) \\ 
    & = R\Hom_{\n{D}^{la}(G,K)}(W, R\varprojlim_h V_{h}),
    \end{align*}
    proving that the right adjoint of $j^*$ is $j_*$ as wanted. 
\end{proof}

The main goal of this section is to prove the following theorem. 

\begin{theorem} \label{TheoremEquivalenceLocAn1}
Let $G$ be a compact $p$-adic Lie group over $L$.  Then the map $j_*$ is  fully faithful and has essential image the $\n{D}^{la}(G,K)$-modules $W$ such that the map
\[
 W\to \varprojlim_h(\n{D}^{h}(G,K)\otimes^L_{\n{D}^{la}(G,K)} W)
\]
is an equivalence.   Furthermore, the map $j^*$ has a left adjoint $j_!$ given by
\[ j_! : \Mod^{qc}_{\n{K}_{\sol}}(\n{D}^{la}(G, K)) \to \Mod_{\n{K}_{\sol}}(\n{D}^{la}(G, K))  \]
\[ j_! (V_h)_{h}= (R\varprojlim_{h} V_h)^{Rla}. \] The functor $j_!$ is fully faithful, and  $j_!j^*W=W^{Rla}$ for all $W\in \Mod_{\n{K}_{\sol}}(\n{D}^{la}(G,K))$, so that   the essential image of $j_!$ is the category $\Rep^{la}_{\n{K}_\sol}(G)$. In particular, the functor $j^*$ induces an equivalence of (stable $\infty$)-categories
\[  \Rep^{la}_{\n{K}_{\sol}}(G)  \xrightarrow{\sim} \Mod^{qc}_{\n{K}_{\sol}}(\n{D}^{la}(G, K)). \]
\end{theorem}

To construct the functor   $j_!$ we shall exploit the fact that the maps $\n{D}^{h'}(G,K)\to \n{D}^{h}(G,K)$ and $C^h(G,K)\to C^{h'}(G,K)$ are of trace class for $h'>h$ (this follows essentially from the standard fact that the transition map between a strict inclusion of rigid analytic balls is compact, see the discussion of \cite[Remark 4.9]{SolidLocAnRR}). Moreover, they factor through $\n{D}^{la}$-modules (see \cite[Remark 4.9]{SolidLocAnRR})
\[
\n{D}^{h'}(G,K)\to \overline{\n{D}}^{h'}(G,K)\to  \n{D}^{h}(G,K)
\]
and 
\[
C^{h}(G,K) \to \overline{C}^{h}(G,K) \to C^{h'}(G,K)
\]
with $\overline{\n{D}}^{h'}(G,K)$ and $\overline{C}^{h}(G,K)$ being compact projective as $\n{K}_{\sol}$-vector spaces.  We will write $C^{h',B}(G,K)$ and $\n{D}^{h,B}(G,K)$ for the duals of $\overline{\n{D}}^{h'}(G,K)$ and $\overline{C}^{h}(G,K)$ respectively, these are $K$-Banach spaces.  Before discussing the proof of Theorem \ref{TheoremEquivalenceLocAn1} we need a series of technical lemmas.

\begin{lemma}
\label{LemmaTraceClassFactor}
Let $f:V\to W$ be a trace class  map of $\n{K}_{\sol}$-vector spaces. Then  there is a morphism $R\iHom_{K}(W,-)\to V^{\vee} \otimes_{\n{K}_{\sol}}^L-$ making the following diagram commutative
\[
\begin{tikzcd}
W^{\vee}\otimes_{\n{K}_{\sol}}^L - \ar[r] \ar[d] &  R\iHom_{K}(W,-)  \ar[ld] \ar[d]\\ 
V^{\vee}\otimes_{\n{K}_{\sol}}^L-  \ar[r]& R\iHom_{K}(V,-).
\end{tikzcd}
\]
Similarly, there is a morphism $R\iHom(V^{\vee},-)\to W\otimes^L_{\n{K}_{\sol}}-$ making the following diagram commutative
\[
\begin{tikzcd}
V\otimes^L_{\n{K}_{\sol}} -  \ar[r] \ar[d] & R\iHom_{K}(V^{\vee},-)  \ar[ld] \ar[d]\\ 
 W\otimes^L_{\n{K}_{\sol}} - \ar[r]& R\iHom_K(W^{\vee},-).   
\end{tikzcd}
\]

\end{lemma}
\begin{proof}
    This is analogue to \cite[Lemma 8.2]{CondensesComplex}. By definition the map $f$ arises from a morphism $K\to V^{\vee}\otimes_{\n{K}_{\sol}}^L W$. For the first diagram we note that, for any $P\in \Mod(\n{K}_{\sol})$, we have  functorial morphisms on $P$
    \begin{align*}
       R\iHom_{K}(W, P) & \to R\iHom_{ K}(V^{\vee}\otimes_{\n{K}_{\sol}}^L W, V^{\vee}\otimes_{\n{K}_{\sol}}^L P) \\ 
                & \to R\iHom_{K}(K, V^{\vee}\otimes_{\n{K}_{\sol}}^L P)\\ 
                &  = V^{\vee}\otimes_{\n{K}_{\sol}}^L P \\ 
                & \to R\iHom_{K}(V,P). 
    \end{align*}
    For the second diagram, we use instead the following functorial morphisms on $P$
    \begin{align*}
    R\iHom_K(V^{\vee},P) & \to  R\iHom_{K} (W \otimes^L_{\n{K}_{\sol}}  V^{\vee}, W \otimes^L_{\n{K}_{\sol}}   P) \\ 
    & \to R\iHom_{K} (K,  W \otimes^L_{\n{K}_{\sol}}   P)  \\ 
    & =   W \otimes^L_{\n{K}_{\sol}}   P \\ 
    & \to R\iHom_K(W^{\vee},P).  
    \end{align*}
\end{proof}

\begin{lemma}\label{LemmaEquivariantVersionTraceClass}
Let $f\colon V \to W$ be a morphism of $\n{D}^{la}(G,K)$-modules which is trace class as $\n{K}_{\sol}$-vector space. Suppose in addition that the map $\widetilde{f}\colon K\to V^{\vee}\otimes^L_{\n{K}_{\sol}} W$  defining $f$ is $\n{D}^{la}(G,K)$-equivariant. Then $\widetilde{f}$ induces natural transformations 
\[
V\otimes_{\n{K}_{\sol}}^L(-)\to    R\iHom_{K}(V^{\vee}, -) \to W\otimes_{\n{K}_{\sol}}^L (-)
\]
and 
\[
V^{\vee} \otimes_{\n{K}_{\sol}}^L(-)\to R\iHom_{K}(V,-)\to W^{\vee}\otimes_{\n{K}_{\sol}}^L(-)
\]
of functors from $\mathrm{Mod}_{\n{K}_{\sol}}(\n{D}^{la}(G,K))\to \mathrm{Mod}_{\n{K}_{\sol}}(\n{D}^{la}(G,K))$.
\end{lemma}
\begin{proof}
 This follows from   Lemma  \ref{LemmaTraceClassFactor} where the only additional observation is that we use the $\n{D}^{la}(G,K)$-equivariance of the map $\widetilde{f}$, together with the natural  Hopf algebra structure of $\n{D}^{la}(G,K)$, to endow  the Hom-spaces and tensor products with the structure of $\n{D}^{la}(G,K)$-modules, and to make the morphisms $\n{D}^{la}(G,K)$-linear. 
\end{proof}

From Lemma \ref{LemmaEquivariantVersionTraceClass} we have the following control on functors of analytic vectors. 

\begin{corollary}\label{CorollaryFactorizationAnFunctions}
For any $h_1>h_2>h_3$, one has the following factorization of functors $\Mod_{\n{K}_{\sol}}(\n{D}^{la}(G,K))\to \Mod_{\n{K}_{\sol}}(\n{D}^{la}(G,K))$.
\begin{enumerate}

\item  
\begin{align*}
R\iHom_{\n{D}^{la}(G,K)}(K, (-)\otimes^L_{\n{K}} C^{h_3}(G,K))   & \to R\iHom_{\n{D}^{la}(G,K)}(\n{D}^{h_2}(G,K), -)  \\
&\to  R\iHom_{\n{D}^{la}(G,K)}(K, (-)\otimes^L_{\n{K}} C^{h_1}(G,K)).
\end{align*}

\item 

\begin{align*}
  \n{D}^{h_1}(G,K) \otimes^L_{\n{D}^{la}(G,K)} (-) &\to  R\iHom_{\n{D}^{la}(G,K)}(\chi [-d]\otimes_{\n{K}_{\sol}} C^{h_2}(G,K) , - )  \\
& \to  \n{D}^{h_3}(G,K) \otimes^L_{\n{D}^{la}(G,K)} (-) 
\end{align*}

\end{enumerate}

\end{corollary}
\begin{proof}
 Recall that the spaces $\n{D}^{h}(G,K)$ are $LB$ space of compact type, so its derived dual agrees with $C^{h}(G,K)$ by Remark \ref{RemarkCommentDualsOrNot}.  Part (1)  follows formally from  Lemma \ref{LemmaEquivariantVersionTraceClass} after taking $\n{D}^{la}(G,K)$-linear Hom spaces from the trivial representation to the factorization  of $\n{D}^{la}(G,K)$-modules
 \[
 (-)\otimes^L_{\n{K}_{\sol}} C^{h_3}(G,K)\to R\iHom_{K}(\n{D}^{h_2}(G,K),-)\to (-)\otimes^L_{\n{K}_{\sol}} C^{h_1}(G,K).
 \]

For part (2), consider instead the factorization of $\n{D}^{la}(G,K)$-modules
\[
( \chi^{-1}[d]\otimes^L_{\n{K}_{\sol}} \n{D}^{h_1}(G,K)) \otimes^L_{\n{K}_{\sol}}(-) \to R\iHom_{K}( \chi[-d]\otimes_{\n{K}_{\sol}} C^{h_2}(G,K), -) \to  ( \chi^{-1}[d]\otimes^L_{\n{K}_{\sol}} \n{D}^{h_3}(G,K)) \otimes^L_{\n{K}_{\sol}}(-). 
\]
obtained from the $\n{D}^{la}(G,K)$-linear  factorizations of Lemma \ref{LemmaEquivariantVersionTraceClass} applied to the trace class maps of distrubution algebras and twisting with $\chi^{-1}[d]$.  The corollary follows by taking $\n{D}^{la}(G,K)$-linear Hom spaces from the trivial representation and Proposition \ref{PropCohomolocyAndHomology} (2) to identify the left and right terms with the base changes along $\n{D}^{la}(G,K)\to \n{D}^{h_i}(G,K)$, as well as Proposition \ref{PropositionHopfAlgebras} (4) for the middle term. 
\end{proof}

\begin{proof}[Proof of Theorem \ref{TheoremEquivalenceLocAn1}]

\textit{Step 1.} Let us first show that the functor
\[
j_*\colon \Mod_{\n{K}_{\sol}}^{qc}(\n{D}^{la}(G,K))\to  \Mod_{\n{K}_{\sol}}(\n{D}^{la}(G,K))
\]
sending $(V_h)_{h}$ to $\varprojlim_{h} V_h$ is fully faithful, and has by essential image those $\n{D}^{la}(G,K)$-modules $W$ such that the natural map 
\[
W\to \varprojlim_{h} (W\otimes_{\n{D}^{la}(G,K)} \n{D}^{h}(G,K)) 
\]is an equivalence. To lighten notations, we will denote $\n{D}^{la} = \n{D}^{la}(G, K)$, $\n{D}^{h} = \n{D}^h(G, K)$ and $C^{h}=C^{h}(G,K)$ for any $h > 0$, and we omit the decoration for derived limits and tensor products. 

Let $(V_h)_{h}\in \Mod_{\n{K}_{\sol}}^{qc}(\n{D}^{la}(G,K))$ and let $V := j_* (V_h)_h = \varprojlim_h V_h$. By Corollary   \ref{CorollaryFactorizationAnFunctions} (2)  we have a pro-equivalence of functors $ \Mod_{\n{K}_{\sol}}(\n{D}^{la}(G,K))\to  \Mod_{\n{K}_{\sol}}(\n{D}^{la}(G,K))$
\[
j^*(-)= \left( \n{D}^h\otimes_{\n{D}^{la}}  (-) \right)_h \xrightarrow{\sim}  \left( \iHom_{ \n{D}^{la}}( \chi[-d]\otimes_{\n{K}_{\sol}}\n{C}^{h}, -) \right)_h. 
\]
Since the right projective diagram of functors preserves pointwise limits,  we deduce an equivalence of pro-objects
\begin{equation}\label{ProEquivqcModules}
j^*V=\left(  \n{D}^h \otimes_{\n{D}^{la}} V \right)_h= \left( \n{D}^h \otimes_{\n{D}^{la}} (\varprojlim_{h'} V_{h'})  \right)_h \xrightarrow{\sim} \left(\varprojlim_{h'} ( \n{D}^h \otimes_{\n{D}^{la}} V_{h'}) \right)_h = ( V_{h} )_{h}.
\end{equation}
The previous pro-equivalence together with Lemma \ref{LemmaCofinalDist} and Remark  \ref{LemmaCofinalDistApplication} imply that the counit map $j^*j_* (V_h)_{h}\to (V_{h})_h$ is an equivalence, proving that $j_*$ is fully faithful as wanted.

Next, we want to identify the essential image of $j_*$ as claimed.
From
Equation \eqref{ProEquivqcModules}, we know that $V=\varprojlim_h (\n{D}^h\otimes_{\n{D}^{la}} V)$. Conversely, if $W=\varprojlim_{h}(\n{D}^h\otimes_{\n{D}^{la}} W)$, then $W=j_*(\n{D}^h\otimes_{\n{D}^{la}} W)_h$ is in the essential image of $j_*$, proving what we wanted.\\

\textit{Step 2.}
We claim that, for $W\in \Mod_{\n{K}_{\sol}}(\n{D}^{la}(G,K))$, the natural map 
\[
j^* (W^{Rla})\to j^* W
\]
is an equivalence.

 Write $W^{Rh-an}:=\iHom_{\n{D}^{la}}(K,W\otimes_{\n{K}_{\sol}} \n{C}^h)$.  By Corollary  \ref{CorollaryFactorizationAnFunctions} (1), we have a natural equivalence of ind-systems
\begin{equation}\label{eqIndEquivSystem}
\left(  W^{Rh-an} \right)_h \xrightarrow{\sim} \left( \iHom_{\n{D}^{la}}(\n{D}^h,W)\right)_h.
\end{equation}
Hence, we deduce that 
\[
W^{Rla}=\varinjlim_h  \iHom_{\n{D}^{la}}(\n{D}^h,W).
\]
Now, using the pro-equivalences of \eqref{eqIndEquivSystem} and Corollary \ref{CorollaryFactorizationAnFunctions} (2) we see that $j^* W^{Rla} = (\varinjlim_{h'>h} \n{D}^{h}\otimes_{\n{D}^{la}} W^{Rh'-an} )_{h}$ is pro-equivalent to the pro-system 
\begin{align*}
\left( \varinjlim_{h'>h} \iHom_{\n{D}^{la}}( \chi[-d]\otimes_{\n{K}_{\sol}} \n{C}^{h}, \iHom_{\n{D}^{la}}(\n{D}^{h'}, W))\right)_{h} & 
= \left( \varinjlim_{h'>h}   \iHom_{\n{D}^{la}} (  \n{D}^{h'} \otimes_{\n{D}^{la}} (\chi[-d]\otimes_{\n{K}_{\sol}} \n{C}^h), W)  \right)_{h} \\ 
&= \left( \iHom_{\n{D}^{la}}(\chi[-d]\otimes_{\n{K}_{\sol}} \n{C}^h, W) \right)_h,
\end{align*}
where in the last equivalence we use the fact that, for $h'>h$, the space $\n{C}^h$ is a $\n{D}^{h}$-module and $\n{D}^{h'}\to \n{D}^h$ is idempotent. We deduce pro-equivalences
\[
j^* W^{Rla} \xrightarrow{\sim}   \bigg( \iHom_{\n{D}^{la}}(\chi[-d]\otimes_{\n{K}_{\sol}} \n{C}^h, W) \bigg)_h \cong j^* W,
\]
proving what we wanted. \\

\textit{Step 3.} Let $W\in \Mod_{\n{K}_{\sol}}(\n{D}^{la}(G,K))$, we claim that the natural map 
\[
W^{Rla} \to j_! j^* W = (j_*j^* W)^{Rla}
\]
obtained by taking locally analytic vectors of the unit map $W\to j_*j^* W $ is an equivalence. In particular, this shows that the essential image of $j_!$ is $\Rep^{la}_{\n{K}_{\sol}}(G)$.

 By Step 2, we know that  $j^* W^{Rla}=j^* W$, so that 
\[
(j_*j^* W^{Rla})^{Rla}= (j_*j^* W)^{Rla}.
\]
Hence, since passing to locally analytic vectors is an idempotent operation by Proposition  \ref{PropRlaIdempotentGroup}, we can assume without loss of generality that $W=W^{Rla}$ is locally analytic. To prove the claim, it suffices to show that we have an equivalence of ind-systems
\[
\left( W^{Rh-an}\bigg)_h\xrightarrow{\sim} \bigg( (j_*j^*W)^{Rh-an} \right)_h. 
\]
The right hand side term in the equivalence of ind-systems of \eqref{eqIndEquivSystem} commutes with pointwise limits, hence we have  equivalences of ind-systems 
\begin{align*}
\left( (j_*j^*W)^{Rh-an} \right)_h &= \bigg( \varprojlim_{h'>h}( \n{D}^{h'}\otimes_{\n{D}^{la}} W)^{Rh-an} \bigg)_h  \\ 
&= \bigg( \varprojlim_{h'>h} \big( \iHom_{\n{D}^{la}}(K,( \n{D}^{h'}\otimes_{\n{D}^{la}} W) \otimes_{\n{K}_{\sol}} \n{C}^{h}) \big) \bigg)_h  \\ 
&=  \bigg( \varprojlim_{h'>h}\big(  ( \chi[-d]\otimes_{\n{K}_\sol} \n{D}^{h'}\otimes_{\n{D}^{la}} \iota(W)) \otimes_{\n{D}^{la}} \n{C}^{h}) \big) \bigg)_h   \\ 
&=    \bigg( \varprojlim_{h'>h} \big(   (\chi[-d]\otimes_{\n{K}_\sol}  \iota(W) ) \otimes_{\n{D}^{la}} \n{C}^{h}) \big) \bigg)_h  \\  
&=  \bigg(  \big( \chi[-d]\otimes_{\n{K}_\sol}  \iota(W) \big) \otimes_{\n{D}^{la}} \n{C}^{h} \bigg)_h    \\
&= \bigg(  \iHom_{\n{D}^{la}}( K,  W  \otimes_{\n{K}_{\sol}} \n{C}^{h}) \bigg)_h  \\
&= (W^{Rh-an})_h
\end{align*}
where the first equivalence follows from the definition of the functors $j^*$ and $j_*$, the second and seventh  equivalences from the definition of $(-)^{Rh-an}$, the third and sixth equivalences use Proposition \ref{PropCohomolocyAndHomology} (2),  the fourth equivalence uses the fact that $\n{C}^{h}$ is a $\n{D}^h$-module and that for $h'>h$ the map $\n{D}^{h'}\to \n{D}^h$ is idempotent, and the fifth equivalence is clear. This proves the claim. \\

\textit{Step 4.} Finally, we prove that the functor $j_! (V_h)_h:= (\varprojlim_h V_h )^{Rla}$ is the left adjoint to $j^*$. If this holds, then a standard argument in category theory shows that $j_!$ is fully faithfull since its double right adjoint $j_*$ is fully faithful. Indeed, this follows from the fact that for $V,W$ in $\Mod^{qc}_{\n{K}_{\sol}}(\n{D}^{la}(G,K))$ one has 
\[
\mathrm{Map}_{\Mod^{qc}_{\n{K}_{\sol}}(\n{D}^{la}(G,K))}(j^*j_!V,W)=  \mathrm{Map}_{\Mod^{qc}_{\n{K}_{\sol}}(\n{D}^{la}(G,K))}(V,j^*j_*W) = \mathrm{Map}_{\Mod^{qc}_{\n{K}_{\sol}}(\n{D}^{la}(G,K))}(V,W)
\]
proving that $j^*j_!V\cong V$ naturally on $V$. 
 To show the claim, let $V=(V_h)_h\in \Mod^{qc}_{\n{K}_{\sol}}(\n{D}^{la}(G,K))$ and $W\in \Mod_{\n{K}_{\sol}}(\n{D}^{la}(G,K))$.  By Step 2 the natural map
\begin{equation}\label{eqKeyCOmputationMainLocAnqcTheo}
j^* j_! V = j^*(j_*V)^{Rla}\to j^*j_*V = V
\end{equation}
is an equivalence. Thus, we deduce   natural equivalences  of spectra
\begin{align*}
\Hom_{\n{D}^{la}}(j_!V, W)   &\xleftarrow{\sim}  \Hom_{\n{D}^{la}}(j_!V, W^{Rla})  \\
&\xrightarrow{\sim}\Hom_{\n{D}^{la}}(j_!V, (j_*j^*W)^{Rla})   \\
&\xrightarrow{\sim} \Hom_{\n{D}^{la}}(j_!V, j_*j^*W) \\
&\cong  \Hom_{\n{D}^{la}}(j^*j_!V, j^*W) \\
&\cong \Hom_{\n{D}^{la}}(V, j^*W),
\end{align*}
where the first arrow arises from the fact that $j_! V$ is locally analytic and the right adjoint property of locally anlaytic vectors (see Corollary  \ref{CoroRightAdjointLocAnGroup}), the second equivalence follows from Step 3,  the third equivalence also follows from Corollary  \ref{CoroRightAdjointLocAnGroup}, the fourth equivalence from the adjunction between $j^*$ and $j_*$, and the last equivalence from \eqref{eqKeyCOmputationMainLocAnqcTheo}.  This finishes the proof of the theorem. 
\end{proof}

We now give some examples showing how the equivalence  of Theorem \ref{TheoremEquivalenceLocAn1} behaves. In particular, it does not preserve the natural $t$-structures on both sides and hence does not induce at all an equivalence of abelian categories.

\begin{example}
\label{ExampleComputationjDla}
We have
\begin{enumerate}
\item $j^* \n{D}^{la}(G, K) = (\n{D}^{h}(G,K))_{h}$.

\item $j_! j^* \n{D}^{la}(G, K) = C^{la}(G, K) \otimes_{\n{K}_\sol} \chi[-d]$.
\item $j^* C^{la}(G, K) = (\n{D}^{h}(G,K)\otimes_{\n{K}_\sol} \chi^{-1}[d])_{h}$.

\item If $V$ is a $\n{D}^{h}(G,K)$-module then the sequence $(V)_{h'>h}$ defines an element in $\Mod_{\n{K}_{\sol}}^{qc}(\n{D}^{la}(G,K))$ and one has $j_!(V)_{h}= j_* (V)_{h} = V$. In particular, for each $h > 0$, Theorem \ref{TheoremEquivalenceLocAn1} restricts to the equivalences of \cite[Theorem 4.36]{SolidLocAnRR}.
\end{enumerate}
\end{example}

\begin{proof}
The first point follows by definition. Part
(2) follows from (1) and Corollary  \ref{CoroLocAnDla}. Indeed, we have
\[ j_! j^* \n{D}^{la}(G, K) = (\varprojlim_h \n{D}^{h}(G, K))^{Rla} = \n{D}^{la}(G, K)^{Rla}. \] Applying $j^*$ to the second example, we obtain
\[ j^*C^{la}(G,K) = j^* j_! j^* \n{D}^{la}(G, K) \otimes_{\n{K}_\sol} \chi^{-1}[d] = j^* \n{D}^{la}(G, K) \otimes_{\n{K}_\sol} \chi^{-1}[d] = (\n{D}^{h}(G, K) \otimes_{\n{K}_\sol} \chi^{-1}[d])_h, \]
where for the second equality we used the equivalence of $j^* j_! \to \mathrm{id}$ of Theorem \ref{TheoremEquivalenceLocAn1}. The last point follows directly from the definitions. Indeed, if $V \in \Mod_{\n{K}_\sol}(\n{D}^{la}(G, K))$ is in fact a $\n{D}^{h}(G, K)$-module, then $j^* V = (\n{D}^{h'}(G, K) \otimes_{\n{D}^{la}(G,K)} V)_{h'} = (V)_{h'\geq h}$, which is a constant sequence, and we have $j_*(V)_h = \varprojlim_h V=  V$ and $j_!(V)_h = (j_* V)^{Rla} = V^{Rla} = V$. 
\end{proof}

\begin{remark}
    As the notation suggests,  the functors $j^*$, $j_*$ and $j^!$ should come from a  $6$-functor formalism  of ``non-commutative spaces'' which at the moment is not available.  When $G= \bb{Z}_p$, nevertheless, the functors $j^*$, $j_*$ and $j_!$ can be interpreted as part of the six functors of the open rigid ball of radius one.
\end{remark}

\begin{definition} \label{DefinitionDualityFunctor}
    We define a duality functor on $\n{C}=\Mod^{qc}_{\n{K}_{\sol}}(\n{D}^{la}(G,K))$ by mapping an object $V=(V_h)_{h}$ to  
    \[
    \bb{D}(V) :=j^* \big( \varprojlim_{h} R\iHom_{\n{D}^{la}}(V_h, \n{D}^{h}(G,K)\otimes_{\n{K}_\sol} \chi^{-1} [d])\big)=\varprojlim_{h} j^* R\iHom_{\n{D}^{la}}(V_h, \n{D}^{h}(G,K)\otimes_{\n{K}_\sol} \chi^{-1}[d]). 
    \]
\end{definition}

\begin{lemma} \label{LemmaRhomjstar}
Let $V \in \Mod_{\n{K}_\sol}(\n{D}^{la}(G, K))$, then
\[ j^* R \iHom_{\n{D}^{la}}(V, \n{D}^{la}(G, K)\otimes_{\n{K}_\sol} \chi^{-1}[d]) = \bb{D}(j^*V).
\]
\end{lemma}

\begin{proof}
    We compute
\[ j^* R \iHom_{\n{D}^{la}}(V, \n{D}^{la}) = \bigg(\n{D}^h \otimes_{\n{D}^{la}}^L R \iHom_{\n{D}^{la}}(V, \n{D}^{la})\bigg)_h. \]
By Corollary \ref{CorollaryFactorizationAnFunctions}, this pro-system  is pro-isomorphic to the system
\[ \bigg( R \iHom_{\n{D}^{la}}(C^{h} \otimes_{\n{K}_\sol} \chi[-d], R\iHom_{\n{D}^{la}}(V, \n{D}^{la}) \bigg)_h. \]
But
\begin{align*}
R \iHom_{\n{D}^{la}}(C^{h} \otimes_{\n{K}_\sol} \chi[-d], R\iHom_{\n{D}^{la}}(V, \n{D}^{la})) &=  R \iHom_{\n{D}^{la}}(C^{h} \otimes_{\n{K}_\sol} \chi[-d] \otimes_{\n{D}^{la}} V, \n{D}^{la}) \\
&=  R \iHom_{\n{D}^{la}}(V, R\iHom_{\n{D}^{la}}(C^{h} \otimes_{\n{K}_\sol} \chi[-d], \n{D}^{la})). \\
\end{align*}
Moreover, applying again Corollary \ref{CorollaryFactorizationAnFunctions}, we get that the pro-system $\bigg(R \iHom_{\n{D}^{la}}(C^h \otimes_{\n{K}_\sol} \chi[-d], \n{D}^{la})\bigg)_h$ is pro-isomorphic to the pro-system $\bigg(\n{D}^{h} \otimes_{\n{D}^{la}}  \n{D}^{la}\bigg)_h = (\n{D}^h)_h$. We deduce  a  natural equivalence of pro-systems 
\[
j^* R\iHom_{\n{D}^{la}}(V,\n{D}^{la}) \cong  \bigg(R\iHom_{\n{D}^{la}}(V, \n{D}^{h})\bigg)_{h} = \bigg(R\iHom_{\n{D}^{la}}(V_h, \n{D}^{h})\bigg)_{h}
\]
with $V_h=\n{D}^h\otimes_{\n{D}^{la}} V$.  One deduces the lemma after twisting by $\chi^{-1}[d]$. 
\end{proof}

\begin{proposition} \label{PropDualityjstar}
Let $V\in \Mod_{\n{K}_{\sol}}(\n{D}^{la}(G,K))$, then 
\[
j^* R\iHom_{K}(V,K)= \bb{D}(j^*V)
\]
where we use the involution of $\n{D}^{la}(G,K)$ to see both modules as left $\n{D}^{la}(G,K)$-modules. In other words, the duality functors as $K$-vector spaces or quasi-coherent $\n{D}^{la}(G,K)$-modules are intertwined (modulo a twist) by the localization functor $j^*$. 

\end{proposition}
\begin{proof}
    By definition we have that $j^* R\iHom_{K}(V,K)= \bigg(\n{D}^{h}\otimes_{\n{D}^{la}}^L R\iHom_{K}(K,V)\bigg)_{h}$. By Corollary \ref{CorollaryFactorizationAnFunctions}  the pro-system  $j^*(R\iHom_{K}(V,K))$ is pro-isomorphic to the pro-system
    \[
\bigg( R\iHom_{\n{D}^{la}}(\n{C}^{h}\otimes_{\n{K}_\sol} \chi[-d], R\iHom_{K}(V,K))\bigg)_{h}.    
    \]
Using adjunction twice, we also have that  
  \begin{align*}
        R\iHom_{\n{D}^{la}}(\n{C}^{h}\otimes_{\n{K}_\sol} \chi[-d], R\iHom_{K}(V,K))  & = R\iHom_{K}( \n{C}^{h}\otimes_{\n{K}_\sol} \chi[-d] \otimes_{\n{D}^{la}}^LV, K ) \\ 
        & = R\iHom_{\n{D}^{la}}(V, R\iHom_{K}(C^{h}\otimes_{\n{K}_\sol} \chi [-d], K)).
\end{align*}
Using Corollary \ref{CorollaryFactorizationAnFunctions} we see that the pro-system  $\bigg(R\iHom_{K}(\n{C}^{h}\otimes_{\n{K}_\sol} \chi[-d],K)\bigg)_{h}$ is pro-equivalent  to $\bigg(\n{D}^{h}\otimes_{\n{K}_\sol} \chi^{-1} [d]\bigg)_{h}$.  One deduces that $j^*(R\iHom_{K}(V,K))$ is pro-equivalent to the pro-system $\bigg(R\iHom_{\n{D}^{la}}(V, \n{D}^{h}\otimes_{\n{K}_\sol} \chi^{-1}[d])\bigg)_{h} $. Hence 
\[
j_* j^*(R\iHom_{K}(V,K)) =  \varprojlim_h R\iHom_{\n{D}^{la}}(V, \n{D}^{h}\otimes_{\n{K}_\sol} \chi^{-1}[d]) = R\iHom_{\n{D}^{la}}(V, \n{D}^{la}\otimes_{\n{K}_\sol} \chi^{-1}[d]).
\]   One concludes thanks to  Lemma  \ref{LemmaRhomjstar}.
\end{proof}

\subsection{Admissible and coadmissible representations} \label{SectLocAnModDla2}

Let $G$ be a compact $p$-adic Lie group over $L$.

\begin{definition}
    We define the derived category of coherent  $\n{D}^{la}(G,K)$-modules  to be the inverse limit 
    \[
    \Mod^{coh}_{\n{K}_\sol}(\n{D}^{la}(G,K))= \varprojlim_{h} \Mod^{perf}_{\n{K}_\sol}(\n{D}^{h}(G,K))\subset \Mod^{qc}_{\n{K}_{\sol}}(\n{D}^{la}(G,K))
    \] of perfect $\n{D}^{h}(G,K)$-modules.\footnote{It would be more ``coherent''  to call these objects \textit{perfect} quasi-coherent  $\n{D}^{la}(G,K)$-modules, but  this terminology would conflict with that of just perfect $\n{D}^{la}(G,K)$-modules. However, because of the noetherian and Auslander properties of the rings $\n{D}^{(h)}(G,K)$ introduced below, coherent and perfect $\n{D}^{(h)}(G,K)$-modules agree. }
     Under the fully faithful embedding $j_*: \Mod^{qc}_{\n{K}_{\sol}}(\n{D}^{la}(G,K))\to \Mod_{\n{K}_\sol}(\D^{la}(G, K))$, we denote by $\Mod^{coad}_{\n{K}_\sol}(\n{D}^{la}(G,K))$ the essential image of $\Mod^{coh}_{\n{K}_\sol}(\n{D}^{la}(G,K))$ and call it the derived category of coadmissible $\n{D}^{la}(G, K)$-modules. Analogously, under the equivalence $j_!: \Mod^{qc}_{\n{K}_{\sol}}(\n{D}^{la}(G,K))\to \Rep^{la}_{\n{K}_{\sol}}(G)$, we denote by $\Rep^{ad}_{\n{K}_\sol}(G)$ the essential image of $\Mod^{coh}_{\n{K}_\sol}(\n{D}^{la}(G,K))$ and call it the derived category of admissible locally $L$-analytic representations of $G$. 
\end{definition}

Let us relate $\Rep^{ad}_{\n{K}_\sol}(G)$ with a more classical definition of the category of admissible representations. We first need to recall some properties of the distribution algebras.

\begin{proposition}[{\cite{SchTeitDist}}] \leavevmode
\begin{enumerate}
    \item  There are Banach distribution algebras $\n{D}^{(h)}(G,K)$ with dense and trace class transition maps $\n{D}^{(h')}(G,K)\to \n{D}^{(h)}(G,K)$ for $h'>h$, such that  $\n{D}^{la}(G,K)= \varprojlim_{h} \n{D}^{(h)}(G,L)$ is  presented as a Fr\'echet-Stein algebra. In particular the rings $\n{D}^{(h)}(G,K)$ are noetherian so any finite $\n{D}^{(h)}(G,K)$-module is naturally a Banach space,  and the morphisms of algebras $\n{D}^{la}(G,K)\to \n{D}^{(h)}(G,K)$ and $\n{D}^{(h')}(G,K)\to \n{D}^{(h)}(G,K)$ for $h'>h$ are flat. 

    \item The rings $\n{D}^{(h)}(G,K)$ are Auslander of dimension $d=\dim_L G$. In particular, any $\n{D}^{(h)}(G,K)$-module of finite type has a finite projective resolution of length at most $d$. 
\end{enumerate}
\end{proposition}

\begin{remark} \label{remarkqcDlamods}
The algebras $\n{D}^{(h)}(G,K)$ used by Schneider and Teitelbaum (denoted by $D_r(G,K)$ in \textit{loc. cit.}) are different from those $\n{D}^{h}(G,K)$ used in this paper. It should  be true  that algebras $\n{D}^{h}(G,K)$ are noetherian and Auslander of dimension $d$, and that the transition maps $\n{D}^{la}(G,K)\to \n{D}^{h}(G,K)$ and $\n{D}^{h'}(G,K)\to \n{D}^{h}(G,K)$ are flat  for $h'>h$, see \cite[Theorem 10.5]{CondensesComplex}. On the other hand, the systems $(\n{D}^{(h)}(G,K))_{h}$ and  $(\n{D}^h(G,K))_h$ are cofinal, this implies that we can also write 
\[
\Mod^{qc}_{\n{K}_{\sol}}(\n{D}^{la}(G,K)) = \varprojlim_{h} \Mod_{\n{K}_{\sol}}(\n{D}^{(h)}(G,K)). 
\]
\end{remark}

\begin{corollary}
    The category $\Mod^{coh}_{\n{K}_\sol}(\n{D}^{la}(G,K))$ has a natural $t$-structure with heart given by the abelian category of coadmissible $\n{D}^{la}(G,K)$-modules, i.e. $\D^{la}(G, K)$-modules of the form $V = \varprojlim_h (V_{h})_{h}$, where the $V_h$'s are $\n{D}^{(h)}(G,K)$-modules of finite type such that $\n{D}^{(h)}(G,K)\otimes^{L}_{\n{D}^{(h')}}V_{h'}=V_h$ for $h'>h$. 
\end{corollary}
\begin{proof}
    The flatness of the rings of distribution algebras implies that the $t$-structures on the categories $\Mod^{perf}_{\n{K}_\sol}(\n{D}^{(h)}(G,K))$ are preserved under base change, this shows that $\Mod^{coh}_{\n{K}_\sol}(\n{D}^{la}(G,K))$ has a natural $t$-structure and that the heart is, by definition, the abelian category of coadmissible $\n{D}^{la}(G, K)$-modules of \cite{SchTeitDist}.
\end{proof}

\begin{remark}
    One can ask for the relation of the (triangulated) bounded derived category of the abelian category of  coadmissible $\n{D}^{la}(G, K)$-modules and the homotopy category of the bounded objects in $\Mod^{coh}_{\n{K}_\sol}(\n{D}^{la}(G,K))$. We do not have an answer to this question, however the first could be poorly behaved as the abelian category of coadmissible $\n{D}^{la}(G,K)$-modules might not have enough injectives or projectives. 
\end{remark}

\begin{lemma} \label{LemmaRlaCoadmissible}
    Let $V\in \Mod^{coh,\heartsuit}_{\n{K}_{\sol}}(\n{D}^{la}(G))$ be a coherent $\n{D}^{la}(G, K)$-module in the heart. Then $(j_*V)^{\vee, Rla}$ is a locally analytic representation concentrated in degree $0$.   Furthermore, it agrees with the  $K$-linear dual $\iHom_{K}(j_*V,K)$ in the abelian category of solid $\n{K}_{\sol}$-vector spaces.  
\end{lemma}
\begin{proof}
   Let $V \in \Mod^{coh,\heartsuit}_{\n{K}_{\sol}}(\n{D}^{la}(G))$ be a coherent $\n{D}^{la}(G, K)$-module in the heart. By Remark \ref{remarkqcDlamods}, we can write $V=(V_h)_h$, where $V_h$ is a $\n{D}^{(h)}(G,K)$-module.  In the following we shall denote $\n{D}^{la}(G,K)$ and $\n{D}^{(h)}(G,K)$ simply by $\n{D}^{la}$ and $\n{D}^{(h)}$-respectively.  
   
     By definition of the functor locally anlaytic vectors and Corollary \ref{CorollaryFactorizationAnFunctions} (1),  we have 
    \[
    (j_*V)^{\vee, Rla}= \varinjlim_{h} R\iHom_{\n{D}^{la}}(\n{D}^{(h)}, R\iHom_{K}(j_*V,K))= \varinjlim_h R\iHom_{K}(\n{D}^{(h)} \otimes^L_{\n{D}^{la}} j_*V,K).
    \]
    By Theorem \ref{TheoremEquivalenceLocAn1} we have $j^*j_* V = V$, so that $\n{D}^{(h)} \otimes_{\n{D}^{la}} j_*V= V_h$. Therefore
    \[
    (j_* V)^{\vee,Rla} = \varinjlim_{h} R\iHom_{K}( V_h, K),
    \]
    but $V_{h'}$ is a $\n{D}^{(h')}$-module of finite presentation, and $V_{h}= \n{D}^{(h)} \otimes_{\n{D}^{(h')}} V_{h'}$. One deduces that $V_{h'}\to V_h$ is a trace class map, defined by a trace map $K\to H^0(V_{h'}^{\vee}) \otimes_{\n{K}_{\sol}} V_{h'}$. Let $W_{h'}:= H^0(V_{h'}^{\vee})$, one then has a factorization 
    \begin{align*}
        V_{h}^{\vee} & \to R\iHom_{K} (W_{h'}\otimes^L_{\n{K}_{\sol}} V_{h}, W_{h'})\\ 
        & \to  W_{h'} \\ 
        & \to V_{h'}^{\vee}
    \end{align*}
    where the first map is the obvious one, the second follows from the trace map $K\to H^0(V_{h'}^{\vee}) \otimes_{\n{K}_{\sol}} V_{h'}$, and the last from the natural map $W_{h'} = H^0(V_{h'}^{\vee})  \to  V_{h'}^{\vee}$. One concludes that 
    \[
    \varinjlim_{h} V_{h}^{\vee} = \varinjlim_{h} W_{h}
    \]
    sits in degree $0$ which proves the lemma. 
\end{proof}

The reader might ask  about the relation between the equivalence provided by Theorem \ref{TheoremEquivalenceLocAn1} and the classical anti-equivalence between the (abelian) category of admissible locally analytic representations and the (abelian) category of coherent $\n{D}^{la}(G, K)$-modules of Schneider and Teitelbaum (we refer the reader to \cite[Theorem 6.3]{SchTeitDist} for the details). In \cite[Proposition 4.42]{SolidLocAnRR} we have shown how one can recover their anti-equivalence from our previous work. The following proposition, which is a summary of many of the previous results of this section, shows how  Schneider and Teitelbaum's anti-equivalence is an instance of the equivalence of Theorem \ref{TheoremEquivalenceLocAn1} after applying duality.

\begin{proposition} \label{PropCommutativeDiagramDual}
We have a commutative diagram
    \[
    \begin{tikzcd}
        \Mod_{\n{K}_{\sol}}(\n{D}^{la}(G,K))  \ar[r, "j^*"] \ar[d, "(-)^{\vee, Rla}"] & \Mod^{qc}_{\n{K}_\sol}(\n{D}^{la}(G, K)) \ar[d, "\bb{D}(-)"] \\
        \Rep^{la}_{\n{K}_{\sol}}(G)  & \Mod^{qc}_{\n{K}_\sol}(\n{D}^{la}(G, K)) \ar[l, "j_!"'],
    \end{tikzcd}
    \]
    where the right vertical arrow is given by the dualizing functor of Definition \ref{DefinitionDualityFunctor}. Moreover, when restricted to the  abelian category of coadmissible $\n{D}^{la}(G,K)$-modules, the composition $j_! \circ \bb{D} \circ j^*$ restricts to the anti-equivalence of \cite[Theorem 6.3]{SchTeitDist}.
\end{proposition}    
    
\begin{proof}
We first prove that the diagram is commutative. By Proposition \ref{PropDualityjstar}, we know that $\bb{D} \circ j^* V = j^* R \iHom_{K}(V, K)$, so that 
\[ j_! \circ \bb{D} \circ j^* = (j_! j^* V^\vee) = (V^\vee)^{Rla} \]
by the third step of the proof of Theorem \ref{TheoremEquivalenceLocAn1}. Lemma \ref{LemmaRlaCoadmissible} shows that, when we restrict to the subcategory $\Mod^{coh, \heartsuit}_{\n{K}_\sol}(\n{D}^{la}(G, K))$, this composition of functors is concentrated in degree $0$ and  coincides with $V \mapsto \iHom_K(j_*V, K)$ which is an admissible  locally analytic representation  being the $K$-linear dual of a coadmissible $\n{D}^{la}(G,K)$-module.
\end{proof}

\begin{remark}\label{CoroDualityLocAn}
  We define a duality functor in $\Mod_{\n{K}_{\sol}}(\n{D}^{la}(G,K))$ by 
\[
\bb{D}(W)=R\iHom_{\n{D}^{la}(G,K)}(W, \overline{\n{D}}^{la}(G,K)\otimes_{\n{K}_\sol} \chi^{-1} [d]).
\]
  This  duality functor  is compatible with that on $\Mod^{qc}_{\n{K}_\sol}(\n{D}^{la}(G, K))$ of Definition \ref{DefinitionDualityFunctor}, namely if $W\in  \Mod_{\n{K}_\sol}(\n{D}^{la}(G, K))$, then by Lemma \ref{LemmaRhomjstar} one has $\bb{D}(j^* W) = j^* \bb{D}(W)$. In particular, combining with Proposition \ref{PropCommutativeDiagramDual} we have a commutative diagram 
    \[
\begin{tikzcd}
\Rep^{la}_{\n{K}_{\sol}}(G)  \ar[d, " (-)^{\vee,Rla}"'] & \Mod_{\n{K}_{\sol}}(\n{D}^{la}(G)) \ar[l, "(-)^{Rla}"'] \ar[d, " \bb{D}(-)"] \\  
\Rep^{la}_{\n{K}_{\sol}}(G) & \Mod_{\n{K}_{\sol}}(\n{D}^{la}(G)) \ar[l, "(-)^{Rla}"].
\end{tikzcd}
    \]
\end{remark}

\subsection{Locally analytic representations as comodules of $C^{la}(G,K)$} \label{SectLocAnCoModCla}

Let $G$ be a (not necessarily compact) $p$-adic Lie group over $L$.  In this section we show that the category of locally $L$-analytic representations of $G$ can be undestood as the derived category of quasi-coherent sheaves over a suitable ``classifying stack'' $[*/G^{la}]$ of $G$. Throughout this paper we will only see this stack as a formal object for which the category of quasi-coherent sheaves can be defined by hand as a limit of a cosimplicial diagram; an honest definition as a stack will require a notion of stack on analytic rings that we will not explore in this work. 

\begin{definition} \label{DefinitionCosimplicial} \leavevmode
\begin{enumerate}
  \item Let $G$ be a group acting on a space $X$.  We define the  simplicial  diagram $(G^{n}\times X)_{[n]\in \Delta^{op}}$ with boundary maps  $d_{n}^i: G^{n}\times X\to G^{n-1}\times X$ for $0\leq i\leq n$ given by 
   \[
    d_{n}^i(g_n,\cdots , g_1,x) = \begin{cases}
 (g_n,\ldots, g_2,g_1x) & \mbox{ if } i=0 \\
(g_n,\ldots, g_{i+1}g_i, \ldots, g_1,x) & \mbox{ if } 0<i<n \\ 
 (g_{n-1},\ldots, g_1,x) &  \mbox{ if } i=n 
    \end{cases}
 \]
and degeneracy maps $s_{n}^{i}: G^{n}\times X\to G^{n+1}\times X$ for $0\leq i\leq n$ given by sending the tuple $(g_n,\ldots, g_1,x)$ to $(g_n,\ldots,1,\ldots, g_1,x)$ with $1$ in the $i+1$-th coordinate. 

\item  Let $G_0\subset G$ be an open compact subgroup. We define the category of quasi-coherent sheaves on $G^{la}$ to be 
\[
\Mod^{qc}_{\n{K}_{\sol}}(G^{la}):= \prod_{g\in G/G_0} \Mod_{\n{K}_{\sol}}(C^{la}(gG_0,K)) 
\]
(and similarly for $G^{n,la}$).

\item We define the category of quasi-coherent sheaves on the classifying stack $[*/G^{la}]$ to be the limit 
\[
\Mod^{qc}_{\n{K}_{\sol}}([*/G^{la}])= \varprojlim_{[n]\in \Delta}  \Mod_{\n{K}_{\sol}}^{qc}(G^{n,la}).
\]
 
\end{enumerate}
\end{definition}

\begin{remark}
    The definition of $\Mod^{qc}_{\n{K}_{\sol}}(G^{la})$ is made in such a way that for $G$ compact we can see $G^{la}$ as the analytic spectrum of $C^{la}(G,K)$, and that for $G$ arbitrary   $G^{la}= \bigsqcup_{g\in G/G_0} g G_0^{la}$.  Then,  the definition of $\Mod^{qc}_{\n{K}_{\sol}}([*/G^{la}])$ follows the intuition that $[*/G^{la}]$ is the geometric realization of the simplicial space $(G^{n,la})_{n\in \Delta^{op}}$. 
\end{remark}
 
\begin{theorem}
\label{TheoremLocAnRepStack}
    There is a natural equivalence of symmetric monoidal stable $\infty$-categories 
    \[
    \Rep^{la}_{\n{K}_{\sol}}(G) = \Mod^{qc}_{\n{K}_{\sol}}([*/G^{la}]),
    \]
    where the  tensor product in the LHS is the tensor product over $\n{K}_{\sol}$. 
\end{theorem}

We need some tecnical lemmas.

\begin{lemma}
\label{LemmaComoduleAndRepLa}
    There is a natural symmetric monoidal equivalence between  the abelian category $\Rep^{la,\heartsuit}_{\n{K}_{\sol}}(G)$ of solid locally analytic representations of $G$, and the abelian category of comodules of the functor $C^{la}(G,-)$ mapping $V\in \Mod^{\heartsuit}(\n{K}_{\sol})$ to $C^{la}(G,V)=\prod_{g\in G/G_0} (C^{la}(gG_0,K)\otimes_{\n{K}_{\sol}} V)$. 
\end{lemma}
\begin{proof}
We highlight   that the category of comodules for the  functor $C^{la}(G,-)$ is abelian thanks to the $t$-exactness properties stablished in Lemma \ref{LemmaFunctorFunctionsAndAction}.   Given a map $\s{O}:V\to C^{la}(G,V)$ we have a morphism  $V\to C^{la}(G,V)\to \iHom_{K}(\n{D}^{la}(G,K),V)$ which by adjunction gives rise a map $\rho:\n{D}^{la}(G,K)\otimes_{K}V\to V$. If $\s{O}$ is a comodule structure then $\rho$ is a module structure and $V$ defines an object in $\Mod_{\n{K}_{\sol}}^{\heartsuit}(\n{D}^{la}(G,K))$. Restricting the comodule structure to $G_0$ one finds that the morphism $\s{O}|_{G_0}:V\to C^{la}(G_0,K)\otimes_{\n{K}_{\sol}} V$ lands in the invariants of the $\star_{1,3}$-action of $\n{D}^{la}(G_0,K)$ in right term. Thus, by taking invariants one finds that $V$ is a direct summand of $V^{la}$ which implies that $V$ is locally analytic itself, i.e. $V\in \Rep^{la,\heartsuit}_{\n{K}_{\sol}}(G)$. Conversely, given $V\in \Rep^{la,\heartsuit}_{\n{K}_{\sol}}(G)$ one has an orbit map $\s{O}: V \to  C^{la}(G,V)$ which is clearly a comodule for the functor $C^{la}(G,-)$. It is easy to check that these constructions are inverse each other. 
\end{proof}

\begin{lemma}\label{LemmaDerivedCategoriesHeart}
Let $\n{A}$ be a Grothendieck abelian category and let $\s{D}^+:=\s{D}^+(\n{A})$ be its left bounded derived $\infty$-category obtained by inverting quasi-isomorphisms on left bounded chain-complexes, see \cite[Section 1.3]{HigherAlgebra}. The following hold:

\begin{enumerate}

\item Let $\End^{\mathrm{exact}}(\s{A})$ denote the category of exact endofunctors of $\n{A}$, that is, endofunctors of $\n{A}$ preserving short exact sequences. Then the natural morphism from \cite[Theorem 1.3.3.2]{HigherAlgebra}
\[
\End^{\mathrm{exact}}(\s{A})\hookrightarrow \End(\s{D}^+)
\]
obtained by passing to derived categories is fully faithful and has essential image given by the $t$-exact endofunctors of $\s{D}^+$.

\item Let $T\in \mathrm{Alg}(\End(\s{D}^+))$ be a $t$-exact colimit preserving monad, and let $\mathrm{LMod}_T(\s{D}^+)$ denote its category of left  modules.  Then $\mathrm{LMod}_T(\s{D}^+)$ carries a natural $t$-structure making the forgetful  map $\mathrm{LMod}_{T}(\s{D}^+)\to \s{D}^+$ $t$-exact.  Its heart is the  Grothendieck abelian category given by  $\mathrm{LMod}_T(\s{A})$,  and there is a natural  $t$-exact equivalence of stable $\infty$-categories
\[
\s{D}^{+}(\mathrm{LMod}_T(\s{A})) \xrightarrow{\sim} \mathrm{LMod}_T(\s{D}^+).
\]
Furthermore, this equivalence passes to left completions
\[
\widehat{\s{D}}(\mathrm{LMod}_T(\s{A}))\xrightarrow{\sim} \mathrm{LMod}_{T}(\widehat{\s{D}}). 
\]

\item   Let $U\in \mathrm{coAlg}(\End(\s{A}))$ be a $t$-exact colimit preserving comonad on $\s{D}^+$.   Then $\mathrm{coLMod}_{U}(\s{D}^{+})$  has a natural $t$-structure making the natural forgetful functor $\mathrm{coLMod}_{U}(\s{D}^+)\to \s{D}^+$ $t$-exact. Its heart is the Grothendieck abelian category  given  by $\mathrm{coLMod}_{U}(\s{A})$. Moreover, the natural map  of left completions
\[
\widehat{\s{D}}(\mathrm{coLMod}_{U}(\s{A}))\xrightarrow{\sim} \mathrm{coLMod}_{U}(\widehat{\s{D}})
\]  
is a $t$-exact equivalence of categories.

\end{enumerate}
\end{lemma}
\begin{proof}
\begin{enumerate}

\item  Part (1) is a direct consequence of \cite[Theorem 1.3.3.2]{HigherAlgebra}: the theorem says that right exact endofunctors of $\s{D}$ that carry injective objects of $\s{A}$ into $\s{A}$ identifies with right $t$-exact endofunctors of $\s{A}$. This equivalence identifies $t$-exact endofunctors of $\s{D}$ with exact endofunctors of $\s{A}$.

\item  Since $T$ is a $t$-exact monad, the category of  $T$-modules $\mathrm{LMod}_T(\s{A})$ on $\s{A}$ form a Grothendieck abelian category, and produces a $t$-exact monadic adjunction
\[
  F\colon \s{A} \rightleftharpoons \mathrm{LMod}_T(\s{A}) \colon G
\]
with monad $T=GF$.  Indeed, the $t$-exactness of $T$ guarantees that if $f\colon M\to N$ is a morphism of $T$-modules in $\s{A}$ both $\mathrm{ker}(f)$ and $\mathrm{coker}(f)$ have a natural structure of $T$-module, making $T$-modules in $\s{A}$ an abelian category, and $G$ exact. To see that  $\mathrm{LMod}_T(\s{A})$ satisfies the Grothendieck axioms, one has to see that it is presentable, that is, that it admits filtered colimits and has a set of generators. For filtered colimits, this follows from the fact that $T$ commutes with filtered colimits, so that if $\{M_i\}_i$ is a filtered diagram of $T$-modules in $\s{A}$, their colimit $\varinjlim_i M_i$ in $\s{A}$ has a natural structure of $T$-module. For constructing  a set of generators, let $\{N_{i}\}_i$ be a set of generators of $\s{A}$, then the free $T$-modules $\{T(N_i)\}_i$ form a set of generators in $\mathrm{LMod}_T(\s{A})$.

 By the equivalence of (1), pasing to derived categories we find  a  $t$-exact adjunction 
\[
  F\colon \s{D}^+ \rightleftharpoons \s{D}^+(\mathrm{LMod}_T(\s{A}) )\colon G
\]
with monad $T=GF$.  We claim that $G$ is monadic on derived categories,  by \cite[Theorem 4.7.3.5]{HigherAlgebra} it suffices to show that $G$ is conservative and commutes with $G$-split geometric realizations. For conservativity,   since  $G$ is a stable functor, it suffices to show that if $G(M)=0$ then $M=0$. By $t$-exactness of $G$ we have 
\[
0=H^*(G(M))=G(H^*(M)).
\]
Therefore, as $G$ is monadic on the underlying abelian categories, we know that $G(H^*(M))=0$  yields $H^*(M)=0$ and therefore $M=0$ as wanted.  For commutation with  $G$-split geometric realizations, we note that $G$ is already a colimit preserving functor, namely, it commutes with finite colimits being stable, and  as the original functor on abelian categories commutes with filtered colimits by the assumption on $T$, the functor $G$ on the derived categories also commutes with filtered colimits.  The last claim about left completions is clear as the monad $T$ is $t$-exact so that the left completion of $ \mathrm{LMod}_T(\s{D}^+)$ is precisely $ \mathrm{LMod}_T(\widehat{\s{D}})$.

\item  Let $U$ be a $t$-exact colimit preserving comonad of $\s{A}$. Since $U$ is a $t$-exact comonad, the category $\mathrm{coLMod}_{U}(\s{A})$ of $U$-comodules on $\s{A}$ is abelian, that is, if $f\colon M\to N$ is a morphism of $U$-comodules, the $t$-exactness of $U$ guarantees that both $\mathrm{ker}(f)$ and $\mathrm{coker}(f)$ have $U$-comodule structures. We show that  $\mathrm{coLMod}_{U}(\s{A})$ satisfies the  Grothendieck axioms, i.e. that it is presentable.  Since $U$ commutes with colimits, one deduces that the filtered colimit in $\s{A}$ of a filtered diagram of $U$-comodules has a natural structure of $U$-comodule. It is left to exhibit a set of generators for $\mathrm{coLMod}_{U}(\s{A})$, for that, let $\{N_i\}_{i}$ be a set of generators of $\s{A}$ stable under kernels, cokernels and extensions,  and let $\{U(N_i)\}_{i\in I}$ be the set of $U$-cofree comodules. We let $\{\widetilde{N}_j\}_j$ denote the family of subquotients of $U$-comodules  of the form $U(N_i)$ for $i\in I$ (which is a set since $\s{A}$ is Grothendieck), we note that $\{\widetilde{N}_j\}_j$ is also stable under kernels, cokernels and extensions.  We claim that  $ \{\widetilde{N}_j\}_j$ generate  $\mathrm{coLMod}_{U}(\s{A})$. For that, let $M$ be a left $U$-comodule in $\s{A}$,  we need to show that there is some $j $  and a non-zero map $\widetilde{N}_j\to M$. We claim that the natural map of $U$-comodules $M\to GF(M)=U(M)$ is injective. Indeed, since $F$ is conservative it suffices to show that  $F(M)\to FGF(M)$ is injective, which follows from the fact that  the composition of this map with the conuint $FGF(M)\xrightarrow{\mu(F)} F(M)$ is the identity. Now, we can write $M=\varinjlim_k M_k$ with $M_k\in \{N_i\}_{i}$, thus, as $U$ commutes with colimits, we find that  $U(M)= \varinjlim_k U(M_k)$.  Since $M\to U(M)$ is  an injection of $U$-modules,  there is some $k$ such that $\mathrm{ker}(U(M_k)\oplus M\to U(M))\neq 0$ is non-zero, one can then find a sub $U$-comodule   of   $U(M_k)$  of the form $\widetilde{N}_j$ and an injection of $U$-comodules $\widetilde{N}_j\hookrightarrow M$, proving what we wanted.

We next prove the derived statements.  By construction of the comodule category,  the forgetful functor $F\colon \mathrm{coLMod}_U(\s{A})\to \s{A}$ is $t$-exact and comonadic, with $t$-exact right adjoint $G\colon \s{A}\to \mathrm{coLMod}_U(\s{A})$ (given by the cofree $U$-comodule),  such that $U=FG$ as comonads.  Now, consider the category $\mathrm{coLMod}_U(\s{D}^+)$ of $U$-comodules in $\s{D}^+$, it gives rise to an adjunction 
\[
 F \colon \mathrm{coLMod}_U(\s{D}^+) \rightleftharpoons \s{D}^+ \colon G
\] 
with $U=FG$.  We define  a $t$-structure on $\mathrm{coLMod}_U(\s{D}^+)$  making an object $M$ connective (resp. coconnective) if its image in $\s{D}^+$ is connective (resp. coconnective).   Equivalently, we have 
\[
\tau^{\geq 0 }\mathrm{coLMod}_U(\s{D}^+) = \mathrm{coLMod}_U(\s{D}^{\geq 0})
\]
and 
\[
\tau^{\leq 0} \mathrm{coLMod}_U(\s{D}^+) = \mathrm{coLMod}_U(\s{D}^{\leq 0}).
\]
Since $U$ is a $t$-exact comonad, given $M\in \mathrm{coLMod}_U(\s{D}^+)$ one sees that $\tau^{\geq 0} M$ and $\tau^{\leq 0} M$ have natural structures of  $U$-comodules on $\s{D}^+$. To show that the previous  decomposition defines a $t$-structure, we need to prove that for $M\in \mathrm{coLMod}_U(\s{D}^{\geq  0})$ and $N\in \mathrm{coLMod}_U(\s{D}^{\leq 0})$ the anima of maps 
\[
\mathrm{Map}_{\mathrm{coLMod}_U(\s{D}^+)}(N,M[-1])=0
\]
of $U$-comodules vanishes. 
For that, let $T=GF$ denote the monad of the adjunction, the cobar construction on the split cosimplicial diagram $(T^{\bullet+1}M)_{\Delta}$ yields a natural equivalence in $\mathrm{coLMod}_U(\s{D}^{+})$
\[
M= \varprojlim_{\Delta} T^{\bullet +1} M.
\]
 Since $T$ is $t$-exact, the terms in the totalization  are coconnective. Thus, it suffices to show that 
\[
\mathrm{Map}_{\mathrm{coLMod}_U(\s{D}^+)}(N,T(M)[-1])=0,
\]
but by adjunction we have that 
\[
\mathrm{Map}_{\mathrm{coLMod}_U(\s{D}^+)}(N,T(M)[-1])=\mathrm{Map}_{\mathrm{coLMod}_U(\s{D}^+)}(N,GF(M)[-1]) =\mathrm{Map}_{\s{D}^+}(FN,FM[-1])=0
\]
where in the last equivalence we used the $t$-structure on $\s{D}^+$.

It is clear from the construction that the heart of $\mathrm{coLMod}_U(\s{D}^+)$ is $\mathrm{coLMod}_U(\s{A})$, and that its left completion is precisely $\mathrm{coLMod}_U(\widehat{\s{D}})$. Hence, by \cite[Theorem 1.3.3.2]{HigherAlgebra} we have a natural morphism of left bounded derived categories 
\[
\s{D}^+(\mathrm{coLMod}_U(\s{A})) \to \mathrm{coLMod}_U(\s{D}^+)
\]
and by passing to left completions a $t$-exact map
\begin{equation}\label{MapLeftCompleteComoduleCategories}
\widehat{\s{D}}(\mathrm{coLMod}_U(\s{A})) \to \mathrm{coLMod}_U(\widehat{\s{D}}). 
\end{equation}

We want to show that \eqref{MapLeftCompleteComoduleCategories} is an equivalence of categories. By left completeness, it suffices to show that it is an equivalence when restricted to coconective objects
\[
\s{D}^{\geq 0}(\mathrm{coLMod}_U(\s{A})) \xrightarrow{\sim} \mathrm{coLMod}_U(\s{D}^{\geq 0}). 
\]
For this last claim, note that the adjunction between $F$ and $G$ with comonad $U=FG$ at abelian level produces an adjunction on the coconnective  derived category 
\[
 F \colon \mathrm{coLMod}_U(\s{D}^{\geq 0}) \rightleftharpoons \s{D}^{\geq 0} \colon G
\]
with comonad $U$. It suffices to show that the previous adjunction is comonadic, for what we will apply \cite[Theorem 4.7.3.5]{HigherAlgebra}. We need to show that $F$ is conservative and preserves $F$-split totalizations. We know that $F$ is conservative on abelian categories, by $t$-exactness we deduce that $F$ is conservative in the coconnective derived categories.  On the other hand, the compatibility of $F$ with $F$-split totalizations  follows from the exact same argument of  \cite[Proposition D.6.4.6]{LurieSpectralAlg}.

\end{enumerate}
\end{proof}

\begin{proof}[Proof of Theorem \ref{TheoremLocAnRepStack}]
   We first show that $\Mod_{\n{K}_{\sol}}^{qc}([*/G^{la}])$ is the left completion of the  derived category of its heart, we do this in two different steps. In order to shorten the proof of the theorem we will use basic facts of the theory of analytic stacks and six functor formalisms of Clausen and Scholze,  for example  as discussed in \cite{camargo2026notessolidgeometry}. Let $G_0\subset G$ be an open compact subgroup. Consider the maps of anlaytic stacks 
   \[
   \AnSpec \n{K}_{\sol}\xrightarrow{f}   [\AnSpec \n{K}_{\sol}/G^{la}_0]\xrightarrow{h}   [\AnSpec \n{K}_{\sol} /G^{la}]. 
   \]
   We first prove the claim for $G^{la}_0$. For that,  we note that $G^{la}_0=\AnSpec C^{la}(G_0,K)$ is the analytic spectrum of the ring $C^{la}(G_0,K)$ endowed with the induced analytic ring structure from $\n{K}_{\sol}$. This implies that the functor $f$ is prim and therefore that  $f^*$ admits a colimit preserving  $\n{K}_{\sol}$-linear right adjoint $f_*$ which is compatible with base change and satisfies the projection formula. Since in addition $f$ satisfies $*$-descent by design, the functor $f^*$ is comonadic  by (the opposite of) \cite[Theorem 4.7.5.2]{HigherAlgebra} applied to the \v{C}ech nerve of $f$, where the condition (*) follows from the compatibility of $f_*$ under base change. Since $f^*$ is $\Mod(\n{K}_{\sol})$-linear, proper base change along the pullback square 
   \[
   \begin{tikzcd}
   G_0^{la}\ar[r] \ar[d] & \AnSpec \n{K}_{\sol} \ar[d] \\ 
  \AnSpec \n{K}_{\sol} \ar[r]  & {[\AnSpec \n{K}_{\sol}/G^{la}]}
   \end{tikzcd}
   \] shows that the comonad $U=f^*f_*$ has underlying functor given by tensoring with $C^{la}(G_0,K)$. In particular $U$ is a colimit preserving  $t$-exact comonad and, by Lemma  \ref{LemmaDerivedCategoriesHeart} (3) and (4),  $\Mod^{qc}_{\n{K}_{\sol}}([*/G^{la}_0])$ is the left completion of the derived category of its heart.  
   
Now, we prove the claim for general $G$.   Consider  the map $h$. It is \'etale since the quotient $G^{la}/G_0^{la}=\bigsqcup_{g\in G/G_0} \AnSpec \n{K}_{\sol} \cdot g$ is a discrete set over  $\AnSpec \n{K}_{\sol} $. In particular, we have an equivalence $h^*=h^!$ and the functor $h^*$ admits a $\n{K}_{\sol}$-linear left adjoint $h_!$ compatible with base change and satisfying the projection formula. Since $h$  is surjective, it satisfies $*$-descent and by \cite[Theorem 4.7.5.2]{HigherAlgebra} it is monadic (where the condition (*) follows from the base change property of $h_!$). We show that the (colimit preserving) monad $T=h_!h^*$ is $t$-exact for the natural $t$-structure of $\Mod^{qc}_{\n{K}_{\sol}}([*/G_0^{la}])$, if this holds, then the category $\Mod^{qc}_{\n{K}_{\sol}}([*/G^{la}])$ will have a natural $t$-structure making the functors $h^*$ and $f^*$ $t$-exact and, by Lemma \ref{LemmaDerivedCategoriesHeart} (2) and (4), it is equal to the left completion of the  derived category of its heart.  It suffices to show that $f^*T$ is $t$-exact.   We have a pullback square of analytic stacks 
 \[
 \begin{tikzcd}
 G/G_0\times \AnSpec \n{K}_{\sol} \ar[r,"\widetilde{h}"] \ar[d,"\widetilde{p}"] & \AnSpec \n{K}_{\sol} \ar[d,"p"]  \\ 
 {[\AnSpec \n{K}_{\sol}/G_0^{la} ]} \ar[r,"h"] & {[\AnSpec \n{K}_{\sol}/G^{la}]}
 \end{tikzcd}.
 \]
Therefore, it suffices to prove that $\widetilde{p}^*T$ is $t$-exact. By proper  base change we have natural equivalences of functors
\[
\widetilde{p}^*T = \widetilde{p}^* h^*h_! = \widetilde{h}^*p^* h_! \cong    \widetilde{h}^* \widetilde{h}_! \widetilde{p}^*.
\]
We already know from the case of $G_0$ that  $\widetilde{p}^*$ is $t$-exact (being a disjoint union of copies of  the map $*\to [*/G^{la}_0]$), hence it suffices to see that both  $\widetilde{h}^*$ and $ \widetilde{h}_!$ are $t$-exact, but this is clear since $\widetilde{h}^*$ just sends $V\mapsto (V)_{g\in G/G_0}$ and $\widetilde{h}_!$ sends  $(W_g)_{g\in G/G_0}\mapsto \bigoplus_{g\in G/G_0} W_g$ which are clearly $t$-exact functors.

To conclude the theorem,  we must  identify the heart of $\Mod^{qc}_{\n{K}_{\sol}}([*/G^{la}])$ with the heart of  $\Rep^{la}_{\n{K}_{\sol}}(G)$.  The heart of  $\Mod^{qc}_{\n{K}_{\sol}}([*/G^{la}])$ is the abelian category of descent datum along the diagram 
\[
\begin{tikzcd}
G^{2,la} \ar[r, shift left= 2ex] \ar[r, shift right  =2ex] \ar[r] & G^{la} \ar[r,shift left =1 ex] \ar[r, shift right =1ex]  \ar[l, shift left = 1ex] \ar[l, shift right = 1ex]& * \ar[l]
\end{tikzcd}
\]
given by the simplicial diagram $(G^{n})_{[n]\in \Delta^{\mathrm{op}}_{\leq 2}}$ of Definition \ref{DefinitionCosimplicial} (1), with simplices of dimension $\leq 2$.  An explicit and standard  computation   shows that the abelian category of descent datum above is precisely the category of $C^{la}(G,-)$-comodules on solid $\n{K}_{\sol}$-vector spaces, see for example the proof of  \cite[Lemma 3.4.26]{MannSix}. The theorem follows  by applying Lemma \ref{LemmaComoduleAndRepLa}.
\end{proof}

\begin{corollary}
\label{CoroQuasiCohStackLocAn}
    Let $G$ be a compact $p$-adic Lie group over $L$, then we have natural equivalences of stable $\infty$-categories
    \[
    \Mod^{qc}_{\n{K}_{\sol}}(\n{D}^{la}(G,K)) = \Rep^{la}_{\n{K}_{\sol}}(G) = \Mod^{qc}_{\n{K}_{\sol}}([*/G^{la}]).
    \]
\end{corollary}

\subsection{Classifying stack of rank one $(\varphi, \Gamma)$-modules and locally analytic representations of $\GL_1$} \label{SectpadicLL}

In this section, we explore an interesting application of Theorem \ref{TheoremEquivalenceLocAn1} for the ($L$-analytic) $p$-adic Lie group $\n{O}_L^{\times}$ to the locally analytic categorical $p$-adic Langlands correspondence for $\GL_1$ as formulated in \cite{EGH}. 

We let $\n{X}_1$ be the classifying stack of rank $1$ $(\varphi, \Gamma)$-modules over the Robba ring on affinoid Tate algebras over $\n{K}=(K,K^+)$, c.f. \cite[\S 5]{EGH}. This stack is represented (c.f. \cite[\S 7.1]{EGH}) by the quotient 
\[ [(\widetilde{\n{W}} \times \bb{G}_m^{an}) / \bb{G}_m^{an}]\] with trivial action of $\bb{G}_m^{an}$, where $\widetilde{\n{W}}$ is the rigid analytic weight space of $\n{O}_L^\times$ whose points on an affinoid ring $A$ are given by continuous (eq. $\bb{Q}_p$-locally analytic) characters $\Hom(\n{O}_L^\times, A^{\times})$, and where $\bb{G}_m^{an}$ denotes the rigid analytic multiplicative group. Let $L^{\times}_{\bb{Q}_p}$ be the restriction of scalars of the $p$-adic Lie group $L^{\times}$ from $L$ to $\bb{Q}_p$.  In \cite{EGH}, the authors conjecture that the natural functor 
\begin{equation}
\label{EquationLL}
\f{LL}_p^{la} : \Rep_{\n{K}_\sol}^{la}(L^{\times}_{\bb{Q}_p}) \to \Mod_\sol^{qc}(\n{X}_1) 
\end{equation}
given by $\f{LL}_p^{la}(\pi) = \n{O}_{\n{X}_1} \otimes^L_{\n{D}^{la}(L^\times_{\bb{Q}_p}, K)} \pi$ (c.f. \cite[Equation (7.1.3)]{EGH}) is fully faithful when restricted to a suitable category of ``tempered'' (or finite slope) locally analytic representations. 

On the other hand, for the functor $\f{LL}_p^{la}$ to be fully faithful without restricting to a smaller subcategory of $\Rep^{la}_{\n{K}_\sol}(L^\times_{\bb{Q}_p})$, one can also modify the stack $\n{X}_1$, namely, we consider
\[ \n{X}_1^{mod}:=[\widetilde{\n{W}}\times \bb{G}_m^{alg}/\bb{G}_{m}^{alg}]\]   where $\bb{G}_m^{alg}$ is the analytic space attached to the ring $(K[T^{\pm 1}], K^+)_{\sol}=\n{K}_{\sol}\otimes_{\bb{Z}} \bb{Z}[T^{\pm 1}]$.
To lighten notation we will use the version of $\n{X}_1$ and $\n{X}_1^{mod}$ involving the space $\n{W}\subset \widetilde{\n{W}}$ of $L$-locally analytic characters, and the group $L^{\times}$ instead. The same arguments will hold for the spaces defined over $\bb{Q}_p$. We recall the reader the well known identification, due to Amice, between the algebra $\n{D}^{la}(\n{O}_L^\times, K)$ of locally analytic distributions on $\n{O}_L^\times$ and the the rigid analytic functions $\s{O}(\n{W})$ on the weight space $\n{W}$.

To describe the category of solid quasi-coherent sheaves of the original stack $\n{X}_1$ in terms of representation theory  we need to introduce a certain algebra of ``tempered sequences'' on $\bb{Z}$. 

\begin{definition} \leavevmode
\begin{enumerate}
   \item  We let $\ell^{temp}_{\bb{Z},K} \subset \prod_{\bb{Z}} K$ be the subalgebra with respect to the pointwise multiplication consisting of sequences $(a_n)_{n\in \bb{Z}}$ such that there exists $r>0$ such that $\sup_{n\in \bb{Z}}\{ | a_n| p^{-r|n|}\} <\infty$. Equivalently, let $\s{O}(\bb{G}_{m}^{an})= \varprojlim_{n\to \infty} K\langle  p^nT ,  \frac{p^{n}}{T}\rangle$, then $\ell^{temp}_{\bb{Z},K}= \s{O}(\bb{G}_m^{an})^{\vee}$. We let $\bb{Z}^{temp}$ denote the analytic space defined by the algebra $\ell^{temp}_{\bb{Z},K}$.

\item We let $L^{\times,temp}$ be the analytic space associated with the algebra $C^{temp}(L^{\times},K):= C^{la}(\n{O}_{L}^{\times},K)\otimes_{\n{K}_{\sol}}^L \ell^{temp}_{\bb{Z},K}$ of tempered locally analytic functions on $L^\times$.  Note that, by the compatibility of duality with tensor products of \cite[Theorem 3.40]{SolidLocAnRR},  we have 
\[
C^{temp}(L^{\times},K) = \n{D}^{la}(\n{O}_L^\times, K)^\vee \otimes_{\n{K}_\sol}^L \s{O}(\bb{G}_m^{an})^\vee =  \s{O}(\n{W}\times \bb{G}_{m}^{an})^{\vee}.
\]

\item We let $\Rep^{temp}_{\n{K}_\sol}(L^{\times}) := \Mod^{qc}_{\sol}([*/L^{\times,temp}])$ be the category of tempered (locally analytic) representations of $L^{\times}$.

\end{enumerate}
\end{definition}

\begin{remark}
    In \cite[Definition 13.5]{ClauseScholzeanalyticspaces} Clausen and Scholze  have introduced a notion of analytic space as certain sheaves in the category of analytic rings with respect to steady localizations. The analytic spaces $\bb{Z}^{temp}$ and $L^{\times, temp}$ can be considered in this category, or equivalently, as the presheaves on analytic rings corepresented by the corresponding algebra. 
\end{remark}

\begin{lemma}
    The spaces $\bb{Z}^{temp}$ and  $L^{\times,temp}$ have unique commutative group structures compatible with the natural maps $\bb{Z}\to \bb{Z}^{temp}$ and $L^{\times,la}\to L^{\times,temp}$.  
\end{lemma}

\begin{proof}
    A commutative group structure on $\bb{Z}^{temp}$ and $L^{\times,temp}$ is the same as a commutative Hopf algebra structure on their spaces of functions. Also note that these Hopf algebra structures will be determined by the Hopf algebra structures of $\s{O}(\bb{Z}) = \ell_{\Z, K} := \s{O}(\bb{G}_m^{alg})^\vee = \prod_\Z K$ and  $\s{O}(L^{\times, la}) = C^{la}(L^\times, K) = C^{la}(\n{O}_L^\times, K) \otimes^L_{\n{K}_\sol} \ell_{\Z, K}$, namely, this follows from the fact that the corresponding morphisms of functions are injective. But by definition $\ell_{\bb{Z}, K}^{temp}$ and $C^{temp}(L^{\times},K)$ are the duals of the global sections of $\bb{G}_m^{an}$ and $\n{W}\times \bb{G}_m^{an}$ which are themselves commutative groups, proving that $\ell_{\bb{Z}, K}^{temp}$ and $C^{temp}(L^{\times},K)$  have a natural structure of commutative Hopf algebras. 
\end{proof}

\begin{theorem} 
\label{TheoLLGL1Algebraic}
 There are natural equivalences of stable $\infty$-categories \footnote{In each classifying stack, the action of both groups is trivial. }
\begin{equation}
\label{eqEquivalenceLLALgebraic}
 \Mod^{qc}_{\n{K}_\sol}([\bb{Z}/ L^{\times,la}]) \xrightarrow{\sim} \Mod^{qc}_{\n{K}_\sol}(\n{X}_1^{mod}), \;\;\;  \Mod^{qc}_{\n{K}_\sol}([\bb{Z}^{temp}/ L^{\times, temp}])  \xrightarrow{\sim} \Mod^{qc}_{\n{K}_\sol}(\n{X}_1).
\end{equation}
Furthermore, the functor $\f{LL}_p^{la}$ defined in \eqref{EquationLL} induces equivalences
\begin{equation}
\label{eqCartierDualityLAlgebraic}
\Rep^{la}_{\n{K}_{\sol}}(L^{\times}) \xrightarrow{\sim} \Mod_{\n{K}_\sol}^{qc}(\n{W}\times \bb{G}_m^{alg}), \;\;\; \Rep^{temp}_{\n{K}_{\sol}}(L^{\times}) \xrightarrow{\sim} \Mod_{\n{K}_\sol}^{qc}(\n{W}\times \bb{G}_m^{an}).
\end{equation}

\end{theorem}

\begin{remark}
 The equivalences \eqref{eqEquivalenceLLALgebraic} and  \eqref{eqCartierDualityLAlgebraic}  of Theorem \ref{TheoLLGL1Algebraic}    should follow from a  Cartier duality theory for quasi-coherent sheaves in analytic spaces, this would imply that the natural symmetric monoidal structures are transformed in the convolution products via the Fourier-Mukai transform. In the cases of the theorem,  we will roughly prove that modules over the  Hopf algebras of the groups are equivalent to comodules of the dual Hopf algebras. 
\end{remark}

\begin{prop} \label{PropCartierDualityGm}
    Let $A$ be a flat solid $\n{K}_{\sol}$-algebra. Then there are natural equivalences 
    \[
    \Mod_{\n{K}_{\sol}}^{qc}([\SpecAn A / \bb{G}_m^{alg}]) = \Func(\bb{Z}, \Mod_{\n{K}_{\sol}}(A))
    \]
    and 
    \[
    \Mod_{\n{K}_{\sol}}^{qc}(\SpecAn A \times \bb{G}_m^{alg})= \Mod_{\n{K}_{\sol}}^{qc}([\SpecAn A / \bb{Z}])
    \]
 functorial with respect to base change $B\otimes_{A,\sol}^L-$.   In particular, the same statement holds for analytic spaces glued from flat $\n{K}$-algebras. 
\end{prop}
\begin{proof}
    By \cite[Proposition A.1.2]{MannSix}, the $\infty$-category $\Mod_{\n{K}_{\sol}}^{qc}([\SpecAn A /\bb{G}_m^{alg}])$ is the derived category of $A[T^{\pm 1}]$-comodules over $A$. The datum of a $A[T^{\pm 1}]$-comodule is the same as the datum of a $\bb{Z}$-graded $A$-algebra, namely, given $M$ an $A[T^{\pm 1}]$-comodule and $\s{O}: M \to M \otimes_{A} A[\pm 1]$ the comodule map, one has a graduation $M=\bigoplus_{i} M(i)$ by defining $M(i)=\s{O}^{-1}( M\otimes T^{-i})$. Conversely, if $M=\bigoplus_{i\in \bb{Z}}  M(i)$ one defines the comodule structure $\s{O}: M\to M\otimes A[T^{\pm 1}]$ by mapping $\s{O} : M(i) \xrightarrow{\sim} M(i)\otimes T^{-i}$. We have constructed a natural equivalence 
    \[
    \Mod^{\heartsuit}_{\n{K}_{\sol}}(\SpecAn A / \bb{G}_m) \xrightarrow{\sim} \Func(\bb{Z}, \Mod^{\heartsuit}_{\n{K}_{\sol}}(A)), 
    \]
    taking derived categories we get the first equivalence. 

    For the second one, the category $\Mod^{qc}_{\n{K}_{\sol}}(\SpecAn A\times \bb{G}_{m}^{alg})$ is by definition the derived category of $\n{K}_{\sol}$-solid $A[T^{\pm 1}] = A[\Z]$-modules, i.e. $\Z$-representations on solid $A$-modules. This gives a natural equivalence 
    \[
    \Mod^{\heartsuit}_{\n{K}_{\sol}}( A[T^{\pm 1}]) = \Mod^{qc, \heartsuit}_{\n{K}_{\sol}}([\SpecAn A / \bb{Z}]),
    \]
    taking derived categories one obtains the second equivalence of the lemma.  
\end{proof}

\begin{remark}
\label{RemarkSixFunctorNonsense}
    In the proof of the following lemma we are going to use some facts coming from a $6$-functor formalism for solid quasi-coherent sheaves  of  analytic stacks  over $\bb{Q}_{p,\sol}$ in the $\s{D}$-topology as in \cite[Definition 4.14]{SixFunctorsScholze}. This theory has been partially constructed in \cite{ClausenScholzeCondensed2019} and \cite{CondensesComplex} for schemes or complex analytic spaces, and the methods of \cite[Appendix A.5]{MannSix}, \cite[\S 5-9]{MannSix2} and \cite{SixFunctorsScholze} are enough to give proper foundations.  In particular, we assume that: 
    \begin{enumerate}

    \item The family $E$ of morphisms in the $6$-functor formalism (see \cite[Definition A.5.7]{MannSix}) contains all maps $f:X\to Y$ of rigid spaces. In particular, we have shriek functors $f_!$ and $f^!$ satisfying proper base change and projection formula, and compatible under compositions.

    \item Let $f:\n{A}\to \n{B}$ be  a map of analytic rings that defines a map of analytic spectra $f:\SpecAn \n{B} \to \SpecAn \n{A}$. If the pullback $f^*:\Mod_{\n{A}} \to \Mod_{\n{B}}$ is an open immersion in the sense of  \cite[Proposition 6.5]{CondensesComplex}, then $f\in E$ and $f_!$ is the left adjoint of $f^*$. Similarly, if $\n{B} =\n{B}_{\n{A}/}$ has the induced analytic ring structure, then $f\in E$ and $f_!=f_*$ is the right adjoint of $f^*$. 
    
    \item Smooth morphisms of rigid spaces are cohomologically smooth (c.f. \cite[Definition 5.1]{SixFunctorsScholze}). For partially proper smooth rigid spaces  over a point this follows from the proof of  \cite[Proposition 13.1]{CondensesComplex} for complex analytic spaces.  Moreover, given $f:X\to Y$  a smooth map of rigid spaces, we have that $f^!=f^! \s{O}_Y \otimes f^*$ and  we have a natural isomorphism $f^! \s{O}_{Y} = \Omega^{\dim X-\dim Y}_{X/Y}[\dim X-\dim Y]$, the last equality can be proven via the same argument as in \cite[Theorem 11.6]{ClausenScholzeCondensed2019}.

    \item Being cohomologically smooth is local in the target for the $\s{D}$-topology (see \cite[Definition 4.18 (2)]{SixFunctorsScholze}), this follows from arguments analogue to those of \cite[Lemma 8.7 (ii)]{MannSix2}. In particular,  if $\bb{G}$ is a   smooth rigid group over $\n{K}=(K,K^+)$, and $*= \SpecAn \n{K}_{\sol}$, then $f:*\to [*/\bb{G}]$ is cohomologically smooth. Indeed, by definition  $[*/\bb{G}]$ is the geometric realization of the \v{C}ech nerve $\{\bb{G}^{n}\}_{n\in \Delta^{op}}$, so that the map $*\to [*/\bb{G}]$ is a $\s{D}$-cover and $*\times_{[*/ \bb{G}]} * = \bb{G}$ which is cohomologically smooth over $*$ by (3).

    \item Being cohomologically proper is local in the  target for the $\s{D}$-topology, this follows from the same arguments of \cite[Lemma 9.8 (iii)]{MannSix2}. In particular, if $\bb{G}=\SpecAn A$ is the analytic  affinoid group  associated with a $\n{K}_{\sol}$-algebra with the induced analytic structure, then the map $*\to [*/\bb{G}]$ is cohomologically proper. 
    \end{enumerate}

In this section we do not pretend to give  proper foundations of the theory of analytic stacks or the $6$-functor formalism of solid quasi-coherent sheaves.  Instead, we only give an example of the power of these abstract tools, and their relation with our Theorem \ref{TheoremEquivalenceLocAn1} and categorical Langlands for $\GL_1$. This section is completely independent of the rest of the paper. 
\end{remark}

Before stating the next proposition, we explain how the formalism of categorified locales of \cite{CondensesComplex} allows us to treat $\bb{G}_m^{an}$ and $\bb{G}_m^{alg}$ on an equal footing. Let $\bb{P}^{1,an}_{K}$ be the projective space over $K$ with coordinates $[x,y]$ seen as a rigid space, let  $0=[0:1]$ and $\infty=[1:0]$ be marked points. Then $\bb{P}^{1,an}_{K}$ can be given the structure of categorified local as in \cite[Definition 11.14]{CondensesComplex}. We can identify $\bb{G}_m^{an}$ as the complement of $\{0,\infty\}$ in $\bb{P}^{1,an}_{K}$ as rigid analytic spaces. We can embed
\[ j:\bb{G}_{m}^{an}\subset \bb{G}_m^{alg}\]
as the open  subspace in the sense of categorified locales whose complement is the idempotent $K[T^{\pm 1}]$-algebra 
\[
C=K\{T\}[T^{-1}] \times K\{T^{-1}\}[T]
\]
where $K\{U\}=\varinjlim_{r\to \infty} K\langle \frac{U}{p^r }\rangle$ is the algebra of germs of functions of $\bb{A}^{1,an}_K$ at $0$.
Indeed, by \cite[Proposition 5.3 (4)]{CondensesComplex} the idempotent algebra defined by $\{0,\infty\}$ in $\bb{P}^{1,an}_K$ is equal to 
\[
D= K\{T\}\times K\{T^{-1}\}
\]
with $T=x/y$ (namely, we can write $\{0,\infty\}$ as the intersection over $r >0$ of $B(0, r) \cup B(\infty, r)$, where $B(0, r)$ and $B(\infty, r)$ denote, respectively, the discs around $0$ and $\infty$, with radius $r$). By \cite[Theorem 6.10]{CondensesComplex} we have a natural isomorphism of analytic spaces $\bb{P}^{1,an}_{K}=\bb{P}^{1,alg}_{K}$ between the rigid analytic and the schematic projective spaces (in the notation of \textit{loc. cit.} the rigid analytic and the schematic projective space correspond to $C(X,X)$ and $C(X)$ respectively). Taking pullbacks of $D$ through the map  $\bb{G}_{m}^{alg}\to \bb{P}^{1,alg}_{K}$  one obtains that 
$C= K\{T\}[T^{-1}] \times K\{T^{-1}\}[T]$ is the complement idempotent algebra of $\bb{G}_m^{an}$ in $\bb{G}_{m}^{alg}$ as claimed. 

\begin{prop}
\label{PropCartierDualityGan}
    Let $A$ be an animated solid $\n{K}_{\sol}$-algebra. There are natural equivalences 
     \[
    \Mod^{qc}_{\n{K}_{\sol}}([\SpecAn A / \bb{G}^{an}_m]) = \Mod^{qc}_{\n{K}_{\sol}}( \SpecAn(A\otimes_{\n{K}_{\sol}} \ell^{temp}_{\bb{Z},K} )) 
    \]
    and
    \[
     \Mod^{qc}_{\n{K}_\sol}(\SpecAn A \times \bb{G}_m^{an}) = \Mod^{qc}_{\n{K}_\sol}([\SpecAn (A) / \Z^{ temp}]) 
    \]
      natural with respect to base change $B\otimes^L_{A}-$.   In particular, the same statement holds for analytic spaces glued from   animated $\n{K}_{\sol}$-algebras. 
\end{prop}
\begin{proof}
To simplify notation we will assume that  $A=K$,  the same arguments hold for general $A$.  Let $*=\SpecAn \n{K}_{\sol}$. We start with the proof of the first equivalence. Consider the map $f:*\to [*/\bb{G}_m^{an}]$  of stacks obtained as the geometric realization of the morphism of simplicial analytic spaces 
\begin{equation}
\label{eqCechGm}
f_{\bullet}: (\bb{G}_{m}^{an,n+1})_{n\in \Delta^{op}}\to (\bb{G}_{m}^{an,n})_{n\in \Delta^{op}},
\end{equation}
where the map $f_n: \bb{G}_{m}^{an,n+1}\to \bb{G}_{m}^{an,n}$ is the projection towards the first $n$ components. In particular, as $\bb{G}_m^{an}$ is cohomologically smooth, the map $f$ is cohomologicaly smooth. This implies that $f^! \cong f^* \otimes f^! 1 $ and $f^{!}1$ invertible, which shows that $f^*$ has a left adjoint given by $f_{\natural}= f_{!}(-\otimes f^! 1)$ (the homology). Then, $f^*$ is a conservative functor that preserves limits and colimits and, by Barr-Beck-Lurie theorem \cite[Theorem 4.7.3.5]{HigherAlgebra}, we have a natural equivalence 
\[
\Mod^{qc}_{\n{K}_{\sol}}([*/\bb{G}_m^{an}]) = \Mod_{f^*f_{\natural}}(\Mod(\n{K}_{\sol})). 
\]
By the projection formula, $f^*f_{\natural}$ is a $\Mod(\n{K}_{\sol})$-linear functor, this shows that $\Mod_{f^*f_{\natural}}(\Mod(\n{K}_{\sol}))=\Mod_{\n{K}_{\sol}}(f^*f_{\natural} (K))$ by \cite[Theorem 4.8.4.1 (3)]{HigherAlgebra}. By Lemma \ref{LemmaMonad} below we have that the object $f^* f_{\natural}(K)$ is naturally isomorphic to $\ell_{\Z, K}^{temp}$ as Hopf algebras, and hence we obtain 
\[
\Mod_{\n{K}_{\sol}}^{qc}([*/\bb{G}_m^{an}])= \Mod_{\n{K}_{\sol}}(\ell_{\bb{Z},K}^{temp})
\]
which shows the first part of the lemma.

For the second part, we consider the projection map $q : \bb{G}_m^{an} \to *$ and let $g:\bb{G}_{m}^{alg}\to *$ so that $q = g \circ j$, we also write  $h:*\to [*/\bb{Z}^{temp}]$.  It suffices to prove that the adjunction  $q_{\natural}\rightleftharpoons  q^*$ is comonadic.  Indeed, assuming this, by Barr-Beck-Lurie one has 
\[
\Mod^{qc}_{\n{K}_{\sol}}(\bb{G}_m^{an}) = \mathrm{CoMod}_{q_{\natural}q^* }(\Mod(\n{K}_{\sol})). 
\]
The projection formula implies that the functor $q_{\natural} q^*$ is  $\Mod(\n{K}_{\sol})$-linear so that by Lemma \ref{LemmaMonad} we have
\[
\mathrm{CoMod}_{q_{\natural}q^* }(\Mod(\n{K}_{\sol}))= \mathrm{CoMod}_{\ell^{temp}_{\bb{Z},K}}(\Mod(\n{K}_{\sol})).
\]
Finally, by \cite[Theorem 4.7.5.2 (3)]{HigherAlgebra} we have a natural equivalence
\[ \Mod^{qc}_{\n{K}_{\sol}}([*/\bb{Z}^{temp}]) = \mathrm{CoMod}_{h^*h_*}(\Mod(\n{K}_{\sol})). \] Indeed, the left adjointable condition is a consequence of proper base change as $h$ is a proper map (c.f. Remark \ref{RemarkSixFunctorNonsense} (5)).  Moreover, by projection formula and proper base change $h^*h_*$ is $\Mod(\n{K}_{\sol})$-linear and one has  $\mathrm{CoMod}_{h^*h_*}(\Mod(\n{K}_{\sol}))= \mathrm{CoMod}_{h^*h_*(K)}(\Mod(\n{K}_{\sol}))$, but \cite[Theorem 4.7.5.2 (2)]{HigherAlgebra}  implies that $h^*h_*(K)=\ell^{temp}_{\bb{Z},K}$ as coalgebra. Putting all this together we deduce an equivalence 
\[ \Mod^{qc}_{\n{K}_{\sol}}(\bb{G}_m^{an}) = \Mod^{qc}_{\n{K}_{\sol}}([*/\bb{Z}^{temp}]). \]
This finishes the proof of the second assertion of the proposition under the assumption that the adjunction  $q_{\natural}\rightleftharpoons  q^*$ is comonadic.

We are left to prove comonadicity of the adjunction $q_\natural \rightleftharpoons q^*$:
\begin{itemize}
\item The functor $q_{\natural}$ is conservative: it is (modulo a twist) the composition of the forgetful functor $g_*: \Mod_{\n{K}_{\sol}}(K[T^{\pm 1}])\to \Mod(\n{K}_{\sol})$ and  the fully faithful inclusion $j_!: \Mod^{qc}_{\n{K}_{\sol}}(\bb{G}^{an}_{m})\to \Mod_{\n{K}_{\sol}}(K[T^{\pm 1}])$.

\item The functor $q_{\natural}$ preserves $q_{\natural}$-split totalizations. Since $q_{\natural} = q_!(-\otimes \Omega^1_{\bb{G}_m^{an}}[1])$, it suffices to see that $q_!$ preserves $q_!$-split totalizations. Let $M\in \Mod_{\n{K}_{\sol}}(K[T^{\pm 1}])$, we can write 
\[
q_!( j^*M)= g_*(j_! j^* M )=  [K[T^{\pm 1}]\to C]\otimes_{K[T^{\pm 1}],\sol}^L M,
\]
and since $K[T^{\pm 1}]$ is a Hopf algebra, by Proposition \ref{PropositionHopfAlgebras} (4), we have that 
\begin{equation} \label{eqCohomologyBB}
\begin{aligned}
[K[T^{\pm 1}]\to C]\otimes_{K[T^{\pm 1}],\sol}^L M  & =( [K[T^{\pm 1}]\to C] \otimes_{\n{K}_{\sol}}^L M) \otimes_{K[T^{\pm 1}]}^L K 
\end{aligned}
\end{equation}
where $K[T^{\pm 1}]$ acts antidiagonally in $[K[T^{\pm 1}]\to C] \otimes_{\n{K}_{\sol}}^L M$ and $K$ is the trivial representation of $\bb{G}_m^{alg}$. Observe that $K$ is a perfect $K[T^{\pm 1}]$-module by the exact sequence 
\[
0\to K[T^{\pm 1}]\xrightarrow{T-1} K[T^{\pm 1}] \to K \to 0
\]
and hence the functor $- \otimes_{K[T^{\pm1}]}^L K$ commutes with limits. Let $(M_{n})_{[n]\in \Delta}$ be a cosimplicial diagram in $\Mod_{\n{K}_{\sol}}(K[T^{\pm 1}])$ such that $(j^* M_{n})_{[n]\in \Delta}$ is $q_!$-split. Then we have 
\[
\begin{aligned}
 q_!(\varprojlim_{[n]\in \Delta} j^* M_n) & = q_!( \varprojlim_{[n]\in \Delta}  j^* j_!j^* M_n) \\ 
  &=  q_!(j^* \varprojlim_{[n]\in \Delta}  j_!j^*M_n) \\ 
    & = ([K[T^{\pm1}]\to C]\otimes^L_{\n{K}_{\sol}}\varprojlim_{[n]\in \Delta} q_!j^* M_n ) \otimes_{K[T^{\pm 1}]}^L K \\ 
    & =(\varprojlim_{[n]\in \Delta} ([K[T^{\pm 1}]\to C]\otimes^L_{\n{K}_{\sol}} q_! j^* M_{n} ) )\otimes_{K[T^{\pm 1 }]}^L K \\ 
    & = \varprojlim_{[n]\in \Delta}( ([K[T^{\pm 1}\to C ]\otimes^L_{\n{K}_{\sol}} q_!j^* M_n ) \otimes^L_{K[T^{\pm 1 }]} K) \\ 
    & = \varprojlim_{[n]\in \Delta} q_! j^* M_{n}.
\end{aligned}
\]
In the first equivalence we used that $j^*j_!$ is the identity. In the second equivalence we used that $j^*$ commute with limits being the pullback of an open immersion. In the third equality we use \eqref{eqCohomologyBB}.
The fourth equivalence follows by the  $q_!$-splitness of $(j^*M_n)_{[n]\in \Delta}$ and \cite[Example 3.11]{MathewDescent} which guarantees that the totalization tower is pro-constant, and \cite[Example 3.13]{MathewDescent} which allows to commute the limit with the tensor. The fifth follows since the functor $- \otimes_{K[T^{\pm1}]}^L K$ commutes with limits. The last equality follows from \eqref{eqCohomologyBB} again.

\end{itemize}

\end{proof}

\begin{lemma} \label{LemmaMonad}
Consider the cartesian square 
\begin{equation*}
\begin{tikzcd}
    \bb{G}^{an}_m \ar[r, "q"]  \ar[d, "q"] & * \ar[d, "f"] \\
    * \ar[r, "f"] & {[}*/ \bb{G}_{m}^{an}{].}
\end{tikzcd}
\end{equation*} Then $ f^*f_{\natural}K= q_{\natural}q^* K$ is canonically isomorphic  to $\ell_{\bb{Z}, K}^{temp}$ as Hopf algebras.   
\end{lemma}
\begin{proof}

Let $j: \bb{G}_{m}^{an}\subset \bb{G}_{m}^{alg}$ and $g:\bb{G}_{m}^{alg}\to *$. We have that  
\[
\begin{aligned}
f^*f_{\natural}(K)& =q_{\natural} q^*(K)\\ 
& = q_{!}( \Omega_{\bb{G}_m^{an}}^{1}[1]) \\ 
& \cong q_{!} (\s{O}_{\bb{G}_m^{an}} [1]) \\
& = g_*( j_!  (\s{O}_{\bb{G}_{m}^{an}}  )) [1]\\ 
& = [K[T^{\pm 1}]\to C] [1] \\ 
& \cong  \ell^{temp}_{\bb{Z},K}.
\end{aligned}
\]
The first equality follows from proper base change.  The second one follows from the identity $q_{ \natural} = q_{!}(- \otimes (q! K))$ and Remark \ref{RemarkSixFunctorNonsense} (3). The third one follows since $\Omega^1_{\bb{G}_m^{an}} \cong \s{O}_{\bb{G}_m^{an}}$ by taking the   differential $dT/T$ as a basis. The fourth one follows since $q=g\circ j$, $j$ is an open immersion, and $\bb{G}_m^{alg}= \SpecAn \n{K}[T^{\pm 1}]$ has the induced analytic structure from $\n{K}_{\sol}$, see Remark \ref{RemarkSixFunctorNonsense} (2). The fifth one follows from the formula for $j_!$ for an open immersion given in \cite[Lecture V]{CondensesComplex}. In the last isomorphism   we write $[K[T^{\pm 1}]\to C][1] = TK\{T\}\bigoplus K \bigoplus T^{-1}K\{T^{-1}\}$ to identify it with $\ell^{temp}_{\bb{Z},K}$. This shows that $f^*f_{\natural}(K)$ is a solid $\n{K}_{\sol}$-vector space that is abstractly isomorphic to $\ell^{temp}_{\bb{Z},K}$, which is an $LB$ space of compact type.

We claim moreover that they are actually naturally isomorphic. For this,  by the duality of $LB$ and Fr\'echet spaces of compact type (see \cite[Theorem 3.40]{SolidLocAnRR}) it suffices to see that their duals are naturally isomorphic. Indeed 
\[
R\iHom_{K}(f^*f_{\natural} K,  K ) = R\iHom_{K}(K, f^*f_* K)
\]
and $f^*f_* K$ is naturally isomorphic to $q_*q^* K=  \s{O}(\bb{G}_{m}^{an})$ by smooth base change \cite[Proposition 8.5 (ii.b)]{MannSix2}. This shows that $f^*f_{\natural} K$ is naturally isomorphic to the (abelian) dual of $\s{O}(\bb{G}_{m}^{an})$ which by definition is $\ell_{\bb{Z},K}^{temp}$.

It is left to see that the Hopf algebra structure of $f^*f_{\natural} K= q_{\natural} q^*K$ is identified with the Hopf algebra structure of $\ell^{temp}_{\bb{Z},K}$. The proof of this fact is probably standard but we include it for completeness.  Similarly as before, it suffices to identify the Hopf algebra structure of its dual.  Consider the simplicial diagram 
\[
(\bb{G}_m^{an,n})_{[n]\in \Delta^{\mathrm{op}}}
\]
and let $q_{n}\colon \bb{G}_m^{an,n}\to *$. The algebra structure on $f^*f_{\natural} K\cong q_{\natural}q^* K$ arises from the corresponding simplicial diagram of homology $(q_{n,\natural} \s{O}_{\bb{G}_m^{an,n}})_{[n]\in \Delta^{\mathrm{op}}} = ( (q_{\natural}q^* K)^{\otimes n})_{[n]\in \Delta^{\mathrm{op}}}$ where we identify  $q_{n,\natural} \s{O}_{\bb{G}_m^{an,n}}$ with the $n$-th fold tensor product of $q_{\natural}q^* K$ thanks to the projection formula. In particular, the multiplication of $f^*f_{\natural} K\cong q_{\natural}q^* K$ arises from the homology of the multiplication map 
\[
m\colon \bb{G}_m^{an}\times \bb{G}_m^{an}\to \bb{G}_m^{an},
\]
and this map is precisely the dual to the comultiplication map $\s{O}(\bb{G}_m^{an})\to \s{O}(\bb{G}_m^{an})\otimes_{\n{K}_{\sol}} \s{O}(\bb{G}_m^{an}) = \s{O}(\bb{G}_m^{an,2})$, proving that  $f^*f_{\natural} K \cong \ell^{temp}_{\bb{Z},K}$ as algebras as desired.

Next, we show that the previous identification also holds as coalgebras, for that it suffices to observe that the  coalgebra structure of $q_{\sharp} q^* K $   is dual to the natural algebra structure of the object $q_*q^*K= \s{O}(\bb{G}_m^{an})$, namely, this holds from the fact that  the comonad $q_{\sharp} q^*$ is the dual of the monad  $q_*q^*$. This proves that the natural isomorphism   $q_{\sharp} q^* K \cong \ell^{temp}_{\bb{Z},K}$ is also as coalgebras. 
\end{proof}

We now show the analogue of Proposition \ref{PropCartierDualityGan} for the weight space $\n{W}$.

\begin{prop}
\label{PropositionDualityOL}
    Let $A$ be an animated  $\n{K}_{\sol}$-algebra. Then there are natural equivalences 
    \[
    \Mod_{\n{K}_{\sol}}^{qc}(\SpecAn A \times \n{W}) = \Mod_{\n{K}_{\sol}}( [\SpecAn A/ \n{O}_L^{\times,la}])
    \]
    and 
    \[
    \Mod_{\n{K}_{\sol}} (\SpecAn A \times \n{O}_{L}^{\times,la}) = \Mod_{\n{K}_{\sol}}([\SpecAn A/ \n{W}])
    \]
    natural with respect to base change $B\otimes^L_{A}-$.  In particular, the same statement holds for analytic spaces glued from animated  $\n{K}_{\sol}$-algebras. 
\end{prop}
\begin{proof}
We just mention how to modify the main points of the proof of  Proposition \ref{PropCartierDualityGan}. Since $\n{W}$ is a smooth group over $K$, the only difference with the case of $\bb{G}_m^{an}$ is to find a replacement for $\bb{G}_{m}^{an}\subset \bb{G}_m^{alg}$. We first claim that the pullback map $j^*: \Mod_{\n{K}_{\sol}}(\n{D}^{la}(\n{O}_L^{\times},K))\to \Mod_{\n{K}_{\sol}}^{qc}(\n{W})$ is an open localization as in \cite[Proposition 6.5]{CondensesComplex}. Indeed, by Theorem \ref{TheoremEquivalenceLocAn1} $j^*$ has a fully faithful left adjoint 
\[
j_!:\Mod_{\n{K}_{\sol}}^{qc}(\n{W})\to \Mod_{\n{K}_{\sol}}(\n{D}^{la}(\n{O}_L^{\times},K))
\]
such that 
\begin{equation}
\label{eqLowerShierkOL}
j_!j^*M = M^{Rla} = \iota(C^{la}(\n{O}_{L}^{\times},K)) \otimes \omega^{-1} \otimes_{\n{D}^{la}(\n{O}_L^{\times},K)}^L M, 
\end{equation}
where $\omega$ is a suitable dualizing sheaf and where the last equivalence follows using Corollary \ref{CorollaryFiniteCohoDim} (1). This implies that $j_!$ satisfies the projection formula and that $j^*$ is indeed an open localization.  Then, replacing $K[T^{\pm 1}]$ with $\n{D}^{la}(\n{O}_L^{\times},K)$ and \eqref{eqCohomologyBB} with \eqref{eqLowerShierkOL} the same proof of Proposition \ref{PropCartierDualityGan} holds in this situation.
\end{proof}

\begin{remark}
   The first equivalence of Proposition \ref{PropositionDualityOL} for $A=K$ is Theorem \ref{TheoremEquivalenceLocAn1} together with Theorem \ref{TheoremLocAnRepStack} for $G= \n{O}_L^{\times}$. It should be possible to give a proof of Theorem \ref{TheoremEquivalenceLocAn1} using the more categorical approach of the Proposition \ref{PropositionDualityOL}; we shall do this in a future work. 
\end{remark}

\begin{remark}
The flatness assumption in Proposition \ref{PropCartierDualityGm} can be dropped by using the same arguments as in Proposition \ref{PropCartierDualityGan}. Indeed, one considers the Cartesian square
\[
\begin{tikzcd}
 \bb{Z} \ar[r,"1"] \ar[d,"q"] & * \ar[d, "f"] \\ 
 *\ar[r,"f"] & {[}*/\bb{Z}{]},
\end{tikzcd}
\]
then one computes that the Hopf algebra $f^*f_{\natural}K$ is naturally isomorphic to $K[T^{\pm 1} ]$, and that the conditions of the (co)monadicity theorem are satisfied. 
\end{remark}

We can finally move to the proof of the main result of this section.

\begin{proof}[Proof of Theorem \ref{TheoLLGL1Algebraic}]
We start with the proof of the first equivalence. By Propositions \ref{PropCartierDualityGm} and  \ref{PropositionDualityOL} we have natural equivalences 
\begin{eqnarray*}
\Mod^{qc}_{\n{K}_\sol}([(\n{W} \times \bb{G}_m^{alg}) / \bb{G}_m^{alg}]) &=& \Mod^{qc}_{\n{K}_\sol}(\n{W} \times \bb{G}_m^{alg} \times \Z) \\
&=& \Mod^{qc}_{\n{K}_\sol}([(\n{W} \times \Z) / \Z]) \\
&=& \Mod^{qc}_{\n{K}_\sol}([\Z / (\n{O}_L^{\times, la} \times \Z)]) \\
&=& \Mod^{qc}_{\n{K}_\sol}([\Z / L^{\times, la}]).
\end{eqnarray*}
Observe that, in the third equivalence, we used that
\begin{eqnarray*}
\Mod_{\n{K}_\sol}([(\n{W} \times \Z) / \Z]) &=& \varprojlim_{[n] \in \Delta} \Mod_{\n{K}_\sol}(\n{W} \times \Z \times \Z^n) \\
&=& \varprojlim_{[n] \in \Delta} \Mod_{\n{K}_\sol}([(\Z \times \Z^n) / \n{O}^{\times, la}_L]) \\
&=& \varprojlim_{[n] \in \Delta} \varprojlim_{[m] \in \Delta} \Mod_{\n{K}_\sol}(\Z \times \Z^n \times (\n{O}^{\times, la}_L)^m) \\
&=& \varprojlim_{[n] \in \Delta} \Mod_{\n{K}_\sol}(\Z \times \Z^n \times (\n{O}^{\times, la}_L)^n) \\
&=& \Mod_{\n{K}_\sol}([\Z / (\n{O}^{\times, la}_L \times \Z)]). \\
\end{eqnarray*}

Analogously, Propositions  \ref{PropCartierDualityGan} and   \ref{PropositionDualityOL} show that
\begin{eqnarray*}
\Mod^{qc}_{\n{K}_\sol}([(\n{W} \times \bb{G}_m^{an}) / \bb{G}_m^{an}]) &=& \Mod^{qc}_{\n{K}_\sol}(\n{W} \times \bb{G}_m^{an} \times \Z^{temp}) \\
&=& \Mod^{qc}_{\n{K}_\sol}([(\n{W} \times \Z^{temp}) / \Z^{temp}]) \\
&=& \Mod^{qc}_{\n{K}_\sol}([\Z^{temp} / (\n{O}_L^{\times, la} \times \Z^{temp})]) \\
&=& \Mod^{qc}_{\n{K}_\sol}([\Z^{temp} / L^{\times, temp}]).
\end{eqnarray*}
This finishes the proof of the first equivalences. The second equivalences follow from the exact same arguments and Theorem \ref{TheoremLocAnRepStack}, the fact that the functor defining the equivalence if given by $\f{LL}_p^{la}$ follows from construction and the adjunction of  $j_!$ and $j^*$ in Theorem \ref{TheoremEquivalenceLocAn1}.
\end{proof}

\section{Solid smooth representations}

Let $G$ be a $p$-adic Lie group over a finite extension $L$ of $\bb{Q}_p$ and let $\n{K}=(K, K^+)$ be a complete non-archimedean field extension of $L$.  In this section we construct the $\infty$-category of smooth representations of $G$ on $\n{K}_{\sol}$-vector spaces and study its main properties.

\subsection{Solid smooth representations} \label{SectionSmoothRepresentations}

Let $G$ be a locally profinite group, and let $\Mod_{\n{K}_\sol}(\n{D}^{sm}(G,K))$ be the derived ($\infty$-)category of $\n{D}^{sm}(G,K)$-modules on $\n{K}_{\sol}$-vector spaces. In this paragraph we will define the category of smooth representations of $G$ on  $\n{K}_{\sol}$-vector spaces as a suitable full subcategory of $\Mod_{\n{K}_{\sol}}(\n{D}^{sm}(G,K))$.

\subsubsection{Smooth functions valued in solid vector spaces}

\begin{definition}[{\cite[Definition 3.4.7]{MannSix}}]
Let $S$ be a profinite set and $V\in \Mod(L_{\sol})$, the space of smooth functions from $S$ to $V$ is the solid $L_{\sol}$-vector space given by 
 \[
 C^{sm}(S,V)= \underline{\Cont}(S,\bb{Z})\otimes_{\bb{Z}_\sol} V. 
 \]
 In particular, since $\bb{Z}$ is discrete, we have  $C^{sm}(S,V)= \varinjlim_{i} \underline{\Cont}(S_i, V)$ where $S=\varprojlim_{i} S_i$ is written as a limit of finite subsets. 
\end{definition}

\begin{lemma} [{\cite[Lemma 3.4.8]{MannSix}}] Let $S$ be a profinite set and $V \in \Mod^{\heartsuit}(L_{\sol})$. The following hold
\begin{enumerate}
\item The values of $C^{sm}(S,V)$ at a profinite set $T$ are given by 
 \[
C^{sm}(S,V)(T) := \Cont(S,V(T)), 
 \]
 where $V(T)$ is discrete.

 \item The natural map $C^{sm}(S,V) \to \underline{\Cont}(S,V)$ is injective.

 \end{enumerate}
\end{lemma}

\begin{definition}[{\cite[Definition 3.4.9]{MannSix}}]
    Let $G$ be a locally profinite group and $V\in \Mod(L_{\sol})$ a solid $L$-vector space. We define the space of smooth functions of $G$ with values in $V$ to be the solid $L$-vector space with values at a profinite $T$ given by 
    \[
    C^{sm}(G,V)(T) = \Cont(G, V(T)).
    \]
    Equivalently, if $H\subset G$ is an open compact subgroup, we have that 
    \[
    C^{sm}(G,V) = \prod_{g\in G/H} C^{sm}(gH, V). 
    \]

\end{definition}

\subsubsection{Smooth vectors}

\begin{definition} \leavevmode
\begin{enumerate}
\item Let $V\in \Mod^{\heartsuit}_{\n{K}_{\sol}}(\n{D}^{sm}(G,K))$, the smooth vectors of $V$ are defined by 
\[
V^{sm} = \varinjlim_{H\subset G} V^{H} := \varinjlim_{H\subset G} \iHom_{\n{D}^{sm}(G,K)}(\n{K}_{\sol}[G/H],V)
\]
 where $H$ runs over all the open compact subgroups of $G$. We say that $V$ is a smooth representation of $G$ if the natural map $V^{sm}\to V$ is an isomorphism.

\item We let $(-)^{Rsm}: \Mod_{\n{K}_{\sol}}(\n{D}^{sm}(G,K)) \to \Mod_{\n{K}_{\sol}}(\n{D}^{sm}(G,K))$ be the functor of derived smooth vectors 
\[
V^{Rsm}= \varinjlim_{H\subset G} V^{RH} := \varinjlim_{H\subset G} R\iHom_{\n{D}^{sm}(G,K)}(\n{K}_{\sol}[G/H], V), 
\]
where $H$ runs over all the open compact subgroups of $G$. We say that an object in $\Mod_{\n{K}_{\sol}}(\n{D}^{sm}(G,K))$ is smooth if the natural arrow $V^{Rsm}\to V$ is an equivalence. We let $\Rep^{sm}_{\n{K}_\sol}(G) \subset \Mod_{\n{K}_{\sol}}(\n{D}^{sm}(G,K))$ be the full subcategory consisting of smooth objects.

\end{enumerate}

\end{definition}

\begin{remark}
In (1) of the previous definition we defined smooth vectors for a module over the smooth distribution algebra. One can of course give a similar definition for a solid $G$ representation, namely, if $V\in \Mod^{\heartsuit}(\n{K}_{\sol}[G])$ one defines 
\[
V^{sm}= \varinjlim_{H\subset G} V^{H} := \varinjlim_{H\subset G} \iHom_{\n{K}_{\sol}[G]}(\n{K}_{\sol}[G/H], V).
\]
If $V$ is in addition a $\n{D}^{sm}(G,K)$-module, then both definitions coincide by base change and the fact that $\n{D}^{sm}(G, K) \otimes_{\n{K}_\sol[G]} \n{K}_\sol[G / H] = \n{K}_\sol[G/H]$.  However, at derived level it turns out that the smooth distribution algebra is better suited to define derived smooth representations, e.g., the derived smooth representations will embed fully faithfully into $\Mod_{\n{K}_\sol}(\n{D}^{sm}(G, K))$, but not into $\Mod(\n{K}_\sol[G])$, see \S \ref{s:AdjunctionsCoho} for a more concrete explanation of this fact. 
\end{remark}

We start by proving some basic facts on smooth representations.

\begin{lemma} \label{LemmaSmoothVectorsIndependent}
Let $V \in \Mod_{\n{K}_\sol}(\n{D}^{sm}(G, K))$. Then, for any open subgroup $G' \subseteq G$,  the natural map $V^{Rsm} \to V$ induces an equivalence of $\n{D}^{sm}(G',K)$-modules 
\[ V^{Rsm}|_{G'} \xrightarrow{\sim} (V|_{G'})^{Rsm}. \]
Moreover, we have $(V^{Rsm})^{Rsm}=V^{Rsm}$.   In particular, the derived category $\Rep^{sm}_{ \n{K}_{\sol}}(G) \subset \Mod_{\n{K}_{\sol}}(\n{D}^{sm}(G,K))$ is stable under all colimits. 
\end{lemma}

\begin{proof}
For $V \in \Mod(\n{D}^{sm}(G, K))$ and any open compact $H \subseteq G'$, since $\n{K}_\sol[G / H] =\n{D}^{sm}(G,K)\otimes^L_{\n{D}^{sm}(G',K)} \n{K}_{\sol}[G'/H]$, by base change we have an equivalence
\[ V^{Rsm} = \varinjlim_{H \subset G} R\iHom_{\n{D}^{sm}(G, K)} (\n{K}_{\sol}[G/H], V) = \varinjlim_{H \subset G'} R\iHom_{\n{D}^{sm}(G', K)}(\n{K}_{\sol}[G'/H], V ) = (V|_{G'})^{Rsm} \]
of $\n{D}^{sm}(G', K)$-modules. This proves the first claim.

For the second one, we need to show that the natural map $(V^{Rsm})^{Rsm} \to V^{Rsm}$ of $\n{D}^{sm}(G, K)$-modules is an equivalence. We can assume by the first assertion that $G$ is compact. Then we have
\begin{align*}
(V^{Rsm})^{Rsm} &= \varinjlim_{H \subset G} R \iHom_{\n{D}^{sm}(G)}(\n{K}_\sol[G / H], \varinjlim_{H' \subset G} R \iHom_{\n{D}^{sm}(G, K)}(\n{K}_\sol[G / H'], V)) \\
&= \varinjlim_{H \subset G} \varinjlim_{H' \subset G} R \iHom_{\n{D}^{sm}(G)}(\n{K}_\sol[G / H], R \iHom_{\n{D}^{sm}(G, K)}(\n{K}_\sol[G / H'], V)) \\
&= \varinjlim_{H \subset G} R \iHom_{\n{D}^{sm}(G)}(\n{K}_\sol[G / H] \otimes^L_{\n{D}^{sm}(G, K)}\n{K}_\sol[G / H], V) \\
&= \varinjlim_{H \subset G} R \iHom_{\n{D}^{sm}(G)}(\n{K}_\sol[G / H], V) \\
&= V^{Rsm},
\end{align*}
where the first and last equalities follow from definition, the second one from the fact that $\n{K}_\sol[G / H]$ is a compact $\n{D}^{sm}(G, K)$-module, the third one follows from adjunction, and the fourth one follows since $\n{K}_\sol[G / H]$ is idempotent over $\n{D}^{sm}(G, K)$ (c.f. Lemma  \ref{CoroSmoothAlgebraMatrix} (2)). 

Finally, for the last statement, let $\{V_{i}\}_{i\in I}$ be a colimit diagram of smooth representations, to check that $\varinjlim_{i} V_i$ is smooth we can restrict to $G$ compact, in this case we have that 
\begin{align*}
(\varinjlim_{i} V_i)^{Rsm} & = \varinjlim_{H\subset G} R\iHom_{\n{D}^{sm}(G,K)}(\n{K}_{\sol}[G/H], \varinjlim_{i}  V_i) \\ 
& =  \varinjlim_{i} \varinjlim_{H\subset G} R\iHom_{\n{D}^{sm}(G,K)}(\n{K}_{\sol}[G/H], V_i) \\ 
& = \varinjlim_{i} V_i^{Rsm}= \varinjlim_{i} V_i,
\end{align*}
where in the second equality we used again the compactness of the $\n{D}^{sm}(G_0, K)$-module $\n{K}_\sol[G_0 / H]$
\end{proof}

The following two lemmas describe the smooth vectors in a similar way as we have previously defined continuous and locally analytic vectors (c.f. \cite{SolidLocAnRR}).

\begin{lemma} \label{LemmaExactSmooth}
The functor $V\mapsto C^{sm}(G,V)$ of smooth functions induces a $t$-exact functor of derived categories 
\[
C^{sm}(G,-): \Mod(\n{K}_{\sol}[G]) \to \Mod(\n{K}_{\sol}[G^3])
\]
and 
\[
C^{sm}(G,-): \Mod_{\n{K}_{\sol}}(\n{D}^{sm}(G,K)) \to \Mod_{\n{K}_{\sol}}(\n{D}^{sm}(G^{3},K))
\]
where $(g_1,g_2,g_3)$ acts on a function $f:G\to V$ by $((g_1,g_2,g_3)\cdot f)(h)= g_3f(g_1^{-1}hg_2)$.
\end{lemma}
\begin{proof}
    Let $V\in \Mod^{\heartsuit}(\n{K}_{\sol})$. If $G_0$ is compact we have that $C^{sm}(G_0,V)= \varinjlim_{H\subset G_0} \Hom_{K}(\n{K}_{\sol}[G_0/H],V)$. One deduces that the functor $V\mapsto C^{sm}(G_0,V)$ is exact and that it is a $\n{D}^{sm}(G_0,K)$-module for the left and right regular actions. This implies the lemma for $G=G_0$ compact.  For general $G$ and $V\in \Mod^{\heartsuit}(\n{K}_{\sol})$, by definition we have that $C^{sm}(G, V)= \prod_{g\in G/G_0} C^{sm}(gG_0, V)=\iHom_{\n{D}^{sm}(G_0,K)}(\n{D}^{sm}(G,K), C^{sm}(G_0,V))$ for both the left or right regular action of $G_0$ on $C^{sm}(G_0,V)$. Therefore the functor $V\mapsto C^{sm}(G,V)$ is exact and the left and right regular actions of $\n{K}_\sol[G]$  factor through the  left and right regular actions of $\n{D}^{sm}(G,K)$, proving the lemma. 
\end{proof}

\begin{lemma}
\label{LemmaSmoothVectorsViaDist}
Let $V\in \Mod_{\n{K}_{\sol}}(\n{D}^{sm}(G,K))$. Then, for any open subgroup $G' \subseteq G$ we have the following equivalence of $\n{D}^{sm}(G', K)$-modules:
\[ V^{Rsm}|_{G'} = R\iHom_{\n{D}^{sm}(G', K)}(K, C^{sm}(G', V)_{\star 1,3}). \]  
\end{lemma}

\begin{proof}
We start by proving the result for a compact subgroup. Let $G_0 \subset G$ be a compact open and let $V \in \Mod_{\n{K}_{\sol}}(\n{D}^{sm}(G_0, K))$. We recall that we have  
\begin{equation} \label{Eqtmp1}
C^{sm}(G_0,V) = \varinjlim_{H \subset G_0} C(G_0 / H, V) = \varinjlim_{H \subset G_0} R\iHom_{K}(\n{K}_{\sol}[G_0/H], V)
\end{equation}
where $H$ runs over all the normal  open compact subgroups. Note that the $\star_{1,3}$-action on the LHS translates to the contragradient action of the RHS, i.e. the diagonal action of $\n{D}^{sm}(G, K)$ acting on $V$ and on $\n{K}_\sol[G_0 / H]$ via composition of the map induced by the antipode and left multiplication (heuristically we have  $g\cdot f (x)= gf(g^{-1}x)$ for $f\in R\iHom_{\n{K}}(\n{K}_{\sol}[G_0/H], V)$ and $x\in \n{K}_{\sol}[G_0/H]$).  Taking $G_0$-invariants in Equation \eqref{Eqtmp1} (c.f. Proposition \ref{PropositionHopfAlgebras} (4)) and since $K$ is a direct summand of  $\n{D}^{sm}(G_0, K)$, we obtain
\begin{align*}
R\iHom_{\n{D}^{sm}(G_0,K)}(K, C^{sm}(G_0, V)_{\star 1,3}) &= R\iHom_{\n{D}^{sm}(G_0, K)}(K, \varinjlim_{H \subset G_0} R\iHom_{K	}(\n{K}_{\sol}[G_0/H], V)) \\
&= \varinjlim_{H \subset G_0} R\iHom_{\n{D}^{sm}(G_0, K)}(K, R\iHom_{K}(\n{K}_{\sol}[G_0/H], V)) \\
&= \varinjlim_{H \subset G_0} R\iHom_{\n{D}^{sm}(G_0, K)}(\n{K}_{\sol}[G_0/H], V) \\
&= V^{Rsm},
\end{align*}
where the first and fourth equivalences hold by definition, the second one uses the fact that $K$ is a compact $\n{D}^{sm}(G, K)$-module, and the third equivalence follows from Prposition \ref{PropositionHopfAlgebras} (4).

We now treat the general case. By Lemma \ref{LemmaSmoothVectorsIndependent} we can assume $G' = G$. First observe that for $V\in \Mod^{\heartsuit}_{\n{K}_{\sol}}(\n{D}^{sm}(G,K))$ we have a natural isomorphism
\begin{align*} \iHom_{\n{D}^{sm}(G_0, K)}(\n{D}^{sm}(G,K), C^{sm}(G_0,V)_{\star_{1,3}}) &= \iHom_{\n{D}^{sm}(G_0, K)}( \bigoplus_{g \in G_0 \backslash G} \n{D}^{sm}(G_0,K) \cdot g, C^{sm}(G_0,V)) \\
&= \prod_{g \in G_0 \backslash G} \iHom_{\n{D}^{sm}(G_0, K)}( \n{D}^{sm}(G_0,K) \cdot g, C^{sm}(G_0,V)) \\
&= \prod_{g \in G_0 \backslash G} C^{sm}(G_0 g, V) = C^{sm}(G, V)_{\star_{1,3}}
\end{align*}
where the $G$-action on the first term is induced by the right action on $\n{D}^{sm}(G,K)$ \footnote{Heuristically, the inverse $C^{sm}(G,V)_{\star_{1,3}} \to \iHom_{\n{D}^{sm}(G_0, K)}(\n{D}^{sm}(G,K), C^{sm}(G_0,V)_{\star_{1,3}})$ is given by sending a smooth function $f:G \to V$ to the map $\tilde{f}: G \to C^{sm}(G_0,V)$ given by  $\tilde{f}(g)=(g\star_{1,3} f)|_{G_0}$.}. We deduce, using Lemma \ref{LemmaExactSmooth}, a natural equivalence 
\[
C^{sm}(G,V)_{\star_{1,3}}\xrightarrow{\sim} R\iHom_{\n{D}^{sm}(G_0,K)}(\n{D}^{sm}(G,K),C^{sm}(G_0,V)_{\star_{1,3}}) 
\]
for all $V\in \Mod_{\n{K}_{\sol}}(\n{D}^{sm}(G,K))$. Hence, we get 
\begin{align*}
    R\iHom_{\n{D}^{sm}(G, K)}(K, C^{sm}(G, V)_{\star_{1,3}}) &= R\iHom_{\n{D}^{sm}(G, K)}(K, R\iHom_{\n{D}^{sm}(G_0, K)}(\n{D}^{sm}(G,K), C^{sm}(G_0,V)_{\star_{1,3}})) \\ 
&= R\iHom_{\n{D}^{sm}(G_0, K)}( K, C^{sm}(G_0,V)_{\star_{1,3}}),
\end{align*}   
hence the result reduces to the compact case.
\end{proof}

\begin{lemma} \label{LemmaSmExact}
    Let $V\in \Mod_{\n{K}_{\sol}}(\n{D}^{sm}(G,K))$, then $H^i(V)^{sm}= H^i(V^{Rsm})$ for all $i\in \bb{Z}$, i.e., taking smooth vectors is exact in the abelian category of solid $\n{D}^{sm}(G,K)$-modules. 
\end{lemma}
\begin{proof}
  Taking smooth vectors is independent of the open subgroup of $G$ by Lemma \ref{LemmaSmoothVectorsIndependent}, so we can assume that $G$ is compact. In this case we can write $V^{Rsm}= \varinjlim_{H \subset G} R\iHom_{\n{D}^{sm}(G, K)}(\n{K}_{\sol}[G/H], V)$ where $H$ runs over all the normal open compact subgroups of $G$. But $\n{K}_{\sol}[G/H]$ is a projective $\n{D}^{sm}(G,K)$-algebra (being a direct summand of it, see the proof of Lemma \ref{CoroSmoothAlgebraMatrix}), and hence $\iHom_{\n{D}^{sm}(G, K)}(\n{K}_{\sol}[G/H], -)$ is exact. The lemma follows since taking filtered colimits is exact. 
\end{proof}

\begin{proposition} \label{PropositionRepsmoothGrothendieck}
An object $V\in \Mod_{\n{K}_{\sol}}(\n{D}^{sm}(G,K))$ is smooth if and only if $H^i(V)$ is smooth for all $i\in \bb{Z}$. Therefore, the natural $t$-structure of $\Mod_{\n{K}_{\sol}}(\n{D}^{sm}(G,K))$ induces a $t$-structure on $\Rep^{sm}_{\n{K}_\sol}(G)$. Moreover, $\Rep^{sm,\heartsuit}_{\n{K}_{\sol}}(G)$ is a Grothendieck abelian category and $\Rep^{sm}_{\n{K}_{\sol}}(G)$ is the derived category of its heart. 
\end{proposition}
\begin{proof}
    An object $V\in \Mod_{\n{K}_{\sol}}(\n{D}^{sm}(G,K))$ is smooth if and only if the natural map $V^{Rsm}\to V$ is an equivalence if and only if $H^i(V)^{sm}=H^i(V^{Rsm})=H^i(V)$ for all $i\in \bb{Z}$. The fact that the category $\Rep_{\n{K}_\sol}^{sm,\heartsuit}(G)$ is an abelian Grothendieck category is clear, c.f.  \cite[Lemma 3.4.10]{MannSix}. Note that a system of generators of the category is given by the objects $\n{K}_\sol[G/H]\otimes_{\n{K}_{\sol}} \n{K}_{\sol}[S]$ where $H$ runs over the open compact subgroups of $G$ and $S$ over the ($\kappa$-small) profinite sets.  Let $\s{C}$ be the derived category of $Rep^{sm, \heartsuit}_{\n{K}_{\sol}}(G)$. By \cite[Proposition 1.3.3.7]{HigherAlgebra} we have a natural morphism $\s{C} \to  \Rep^{sm}_{\n{K}_{\sol}}(G)$. To prove that this is an equivalence it suffices to show that for $V,W\in \Mod^{\heartsuit}_{\n{K}_{\sol}}(\n{D}^{sm}(G,K))$ smooth representations we have that 
    \[
    R\Hom_{\n{C}}(V,W)= R\Hom_{\n{D}^{sm}(G,K)}(V,W). 
    \]
This essentially follows from the fact that taking smooth vectors is exact. Indeed, let $I^{\bullet}$ be an injective resolution of $W$ as $\n{D}^{sm}(G,K)$-modules, then  $I^{\bullet, Rsm}= I^{\bullet,sm}$ (by Lemma \ref{LemmaSmExact})  is an injective resolution of $W$ in $\s{C}^{\heartsuit}=\Rep^{sm,\heartsuit}_{\n{K}_{\sol}}(G)$. We have that 
    \begin{align*}
        R\Hom_{\n{D}^{sm}(G,K)}(V,W) & = \Hom_{\n{D}^{sm}(G,K)}(V, I^{\bullet}) \\ 
        & = \Hom_{\n{D}^{sm}(G,K)}(V, I^{\bullet,sm}) \\ 
        & = \Hom_{\n{C}^{\heartsuit}}(V, I^{\bullet,sm}) \\ 
        & = R\Hom_{\n{C}}(V, W), 
    \end{align*}

finishing the proof of the result.    
\end{proof}

\begin{prop} \label{PropAdjunctionSmooth}
    The inclusion $\Rep^{sm}_{\n{K}_{\sol}}(G)\to \Mod_{\n{K}_{\sol}}(\n{D}^{sm}(G,K))$ has a right adjoint given by the smooth vectors functor $V\mapsto V^{Rsm}$. 
\end{prop}

\begin{proof}
Let $G$ be a locally profinite group and let $V\in \Rep^{sm}_{\n{K}_{\sol}}(G)$ and $W\in \Mod_{\n{K}_{\sol}}(\n{D}^{sm}(G,K))$. It suffices to show the adjunction at the level of abelian categories (c.f. \cite[\href{https://stacks.math.columbia.edu/tag/0FNC}{Tag 0FNC}]{stacks-project}), so we can assume both $V$ and $W$ to be in degree $0$. Moreover, since by Proposition \ref{PropositionRepsmoothGrothendieck} the abelian category of smooth representations is generated by $\n{K}_{\sol}[G/H]\otimes_{\n{K}_{\sol}} \n{K}_{\sol}[S]$ for $H\subset G$ open  compact and $S$ profinite, we can assume $V=\n{K}_{\sol}[G/H]\otimes_{\n{K}_{\sol}} \n{K}_{\sol}[S]$. Moreover, since we are computing the internal $\iHom$ we can even assume that $V= \n{K}_{\sol}[G/H]$. But then we have that $R \iHom_{\n{D}^{sm}(G,K)}(\n{K}_{\sol}[G/H],W)= W^{RH}=W^{H}$ are the $H$-invariant vectors which coincide with the $H$-invariant vectors of $W^{sm}$. Indeed, we have, for any open compact $H' \subset H$,
\begin{align*}
R\iHom_{\n{D}^{sm}(G,K)}(\n{K}_{\sol}[G/H],W) &= R\iHom_{\n{D}^{sm}(G,K)}(\n{K}_{\sol}[G/H'] \otimes_{\n{D}^{sm}(G, K)}^L \n{K}_{\sol}[G/H],W) \\
&= R\iHom_{\n{D}^{sm}(G,K)}(\n{K}_{\sol}[G/H],W^{RH'}),
\end{align*}
where in the first line we used adjunction, and the second line  the fact that  $ \n{K}_{\sol}[G/H] = \n{K}_{\sol}[G/H'] \otimes_{\n{D}^{sm}(G, K)}^L \n{K}_{\sol}[G/H]$ 
by the idempotency of  Lemma \ref{CoroSmoothAlgebraMatrix} (2). Taking colimits over $H'$ (and using that $\n{K}_\sol[G/H]$ is projective over $\n{D}^{sm}(G, K)$) shows that
\[ R\iHom_{\n{D}^{sm}(G,K)}(\n{K}_{\sol}[G/H],W) = R\iHom_{\n{D}^{sm}(G,K)}(\n{K}_{\sol}[G/H],W^{sm}), \]
proving what we wanted. 
\end{proof}

\subsection{Smooth representations as quasi-coherent $\n{D}^{sm}(G, K)$-modules} \label{SmoothQuasiCoherent}

In this section we will give a first geometric  description of the category of solid smooth representations of a profinite group $G$, analogous to those appearing in Theorem \ref{TheoremEquivalenceLocAn1} for solid locally analytic representations.

\begin{definition}
    Let $G$ be a profinite group, we define the category of solid quasi-coherent modules over $\n{D}^{sm}(G,K)$ as \[ \Mod^{qc}_{\n{K}_{\sol}}(\n{D}^{sm}(G,K))= \varprojlim_{H\subset G} \Mod_{\n{K}_{\sol}}(\n{D}^{sm}(G/H,K)), \] where $H$ runs over all the normal open subgroups and the transition maps are base changes. We let $j^*: \Mod_{\n{K}_{\sol}}(\n{D}^{sm}(G,K))\to \Mod^{qc}_{\n{K}_{\sol}}(\n{D}^{sm}(G,K))$ be the pullback functor  $j^*W= ( \n{D}^{sm}(G/H)\otimes_{\n{D}^{sm}(G)}^L W)_{H}$.
\end{definition}

\begin{proposition}
\label{PropositionQuasiCoherentSmooth}

Let $G$ be a profinite group. The pullback functor
\[ j^*: \Mod_{\n{K}_\sol}(\n{D}^{sm}(G,K))\to \Mod^{qc}_{\n{K}_{\sol}}(\n{D}^{sm}(G,K)) \] has a right adjoint $j_*( V_{H})_{H}= R \varprojlim_{H} V_{H}$ and a left adjoint $j_! V_{H} = (j_*V)^{Rsm}$. Furthermore, $j^*j_*V=j^*j_!V=V$ for $V\in \Mod^{qc}_{\n{K}_{\sol}}(\n{D}^{sm}(G,K))$ and $j_!j^*W= W^{Rsm}$ for $W\in \Mod_{\n{K}_{\sol}}(\n{D}^{sm}(G,K))$. The functor is a fully faithful embedding with essential image $\Rep^{sm}_{\n{K}_{\sol}}(G)$. 
\end{proposition}
\begin{proof}
 Let $V=(V_H)_{H}\in \Mod^{qc}_{\n{K}_{\sol}}(\n{D}^{sm}(G,K))$ and $W\in \Mod_{\n{K}_{\sol}}(\n{D}^{sm}(G,K))$. One has

\begin{align*}
    R\Hom_{\n{D}^{sm}(G,K)}(W,j_* V) & = R\varprojlim_{H} R\Hom_{\n{D}^{sm}(G,K)}(W,V_H) \\ 
    & = R\varprojlim_H R\Hom_{\n{D}^{sm}(G,K)}( K[G/H] \otimes_{\n{D}^{sm}(G,K)}^L W,V_H)  \\
    &= R\Hom_{\Mod^{qc}_{\n{K}_{\sol}}(\n{D}^{sm}(G,K))}(j^*W,V),
    \end{align*}
   where $H$ runs over open compact subgroups of $G$, the first three Hom spaces are as solid $\n{D}^{sm}(G,K)$-modules,  and  the last Hom is in $\Mod^{qc}_{\n{K}_{\sol}}(\n{D}^{sm}(G,K))$. The previous shows that  $j_*$ is the right adjoint of $j^*$. The other statements of the proposition follow easily by unraveling the definitions, $\otimes$-$\iHom$ adjunction and using the fact that $K[G/H]$ is a direct summand of $\n{D}^{sm}(G,K)$, so in particular compact and dualizable. 
\end{proof}

\subsection{Smooth dualizing functors}

Let $G$ be a compact $p$-adic Lie group over $L$.  The following result answers a question raised by Schneider and Teitelbaum in \cite[p. 26]{SchTeitDuality} on the extension of the contragradient functor for smooth representations to the category of locally analytic representations.

\begin{remark} We first observe that, if $V \in \Mod_{\n{K}_\sol}(\n{D}^{sm}(G, K))$ is a solid smooth representation, then is is also locally analytic. Indeed, for any open compact $H \subset G$, we have
\begin{align*}
(V^{Rsm})^{Rla} &= \varinjlim_{H} (V^{RH})^{Rla} \\
&= \varinjlim_{h, H} R\iHom_{\n{D}^{la}(G,K)}(\n{D}^h(G,K), R \iHom_{\n{D}^{sm}(G, K)}(\n{K}_\sol[G/H], V)) \\
&= \varinjlim_{h,H} R\iHom_{\n{D}^{sm}(G, K)}(\n{D}^h(G,K) \otimes_{\n{D}^{la}(G,K)} \n{K}_\sol[G/H], V) \\
&= \varinjlim_{H} R\iHom_{\n{D}^{sm}(G, K)}(\n{K}_\sol[G/H], V) \\
&= V^{Rsm}
\end{align*}
where the first equivalence follows since taking locally analytic vectors commutes with colimits, the second is just by definition, the third one follows by adjunction, and the fourth one follows from the idempotency of $\n{D}^h(G,K)$ over $\n{D}^{la}(G,K)$ and the fact that $\n{K}_{\sol}[G/K]$ is a $\n{D}^h(G,K)$-module for $h\gg 0$. This provides a natural forgetful functor from the category of solid smooth representations to solid locally analytic representations.
\end{remark}

\begin{proposition} \label{PropContragradientSmoothLocAn}
Let $V \in \Rep^{sm}_{\n{K}_\sol}(G)$. Then 
\[ (V^\vee)^{Rsm} = (V^\vee)^{Rla}. \]
In other words, there is a commutative diagram
\[
\begin{tikzcd}
    \Rep^{sm}_{\n{K}_\sol}(G) \ar[r] \ar[d, "((-)^\vee)^{Rsm}"] & \Rep^{la}_{\n{K}_\sol}(G) \ar[d, "((-)^\vee)^{Rla}"] \\
    \Rep^{sm}_{\n{K}_\sol}(G) \ar[r] & \Rep^{la}_{\n{K}_\sol}(G),
\end{tikzcd}
\]
where the horizontal functors are the natural forgetful functors.
\end{proposition}

\begin{proof}
This is a consequence of Remark  \ref{CoroDualityLocAn} and the analogous calculation for smooth representations, which follow from \cite[Corollary 3.7]{SchTeitDuality}. Indeed these statements assert that both functors are given by the same duality functor in the category $\Mod_{\n{K}_\sol}(\n{D}^{la}(G, K))$. But we give a direct proof. We can and do assume that $G$ is compact, or even a uniform pro-$p$-group. We have
\begin{align*}
(V^\vee)^{Rla} &= \varinjlim_h R \iHom_{\n{D}^{la}(G, K)}(\n{D}^{h}(G, K), R \iHom_{K}(V, K)) \\
&= \varinjlim_h R \iHom_K(\n{D}^{h}(G, K) \otimes_{\n{D}^{la}(G, K)} \n{D}^{sm}(G, K) \otimes_{\n{D}^{sm}(G, K)} V, K) \\
&= \varinjlim_h R \iHom_K(K[G / \mathring{\bb{G}}_{h^+}(L)] \otimes_{\n{D}^{sm}(G, K)} V, K) \\
&= \varinjlim_h R \iHom_{\n{D}^{sm}(G, K)}(K[G / \mathring{\bb{G}}_{h^+}(L)], V^\vee) \\
&= (V^\vee)^{Rsm},
\end{align*}
where the first, second and fourth equalities follow from definition and adjunction, and the third one follows from the equality $\n{D}^{h}(G, K) \otimes_{\n{D}^{la}(G, K)} \n{D}^{sm}(G, K) = K[G / \mathring{\bb{G}}_{h^+}(L)]$ of Lemma \ref{LemmaTensorProducthansmoothDistributions} (we refer to \S \ref{ss:FunctionsSolidV} for the notations). The fifth one follows since the groups $\bb{G}^{(h^+)}(L)$ form a cofinal system of open neighbourhoods of the identity in $G$.
\end{proof}

\subsection{Smooth representations as comodules over $C^{sm}(G, K)$}

We now explain the analogue equivalence of Theorem \ref{TheoremLocAnRepStack} for smooth representations. 

\begin{definition}
    Let $G$ be a locally profinite group and $G_0\subset G$ an open compact subgroup. We let 
    \[
    \Mod^{qc}_{\n{K}_{\sol}}(G^{sm})= \prod_{g\in G/G_0} \Mod_{\n{K}_{\sol}} (C^{sm}(gG_0,K)).
    \]
    We define the quasi-coherent modules of $[*/G^{sm}]$ to be 
    \[
    \Mod^{qc}_{\n{K}_{\sol}}([*/G^{sm}])= R\varprojlim_{[n]\in \Delta} \Mod^{qc}_{\n{K}_{\sol}}(G^{n,sm}).
    \]
\end{definition}

\begin{proposition}
\label{PropositionStackSmooth}
    There is a natural equivalence of symmetric monoidal stable  $\infty$-categories 
    \[
    \Rep^{sm}_{\n{K}_{\sol}}(G)= \Mod^{qc}_{\n{K}_{\sol}}([*/G^{sm}]).
    \]
    In particular, if $G$ is compact, we have natural equivalences of stable $\infty$-categories
    \[
    \Mod^{qc}_{\n{K}_{\sol}}(\n{D}^{sm}(G,K))= \Rep^{sm}_{\n{K}_{\sol}}(G)= \Mod^{qc}_{\n{K}_{\sol}}([*/G^{sm}]).
    \]
\end{proposition}
\begin{proof}
    This follows by the same proof of Theorem \ref{TheoremLocAnRepStack}, the only thing to verify is that the abelian category of smooth representations is naturally equivalent to the abelian category of comodules $V\to C^{sm}(G,V)$, which is obvious. 
\end{proof}

\subsection{Locally algebraic representations of reductive groups}

In this last subsection we introduce a category of solid locally algebraic representations for the $L$-points of a reductive group $\mathbf{G}/L$. Let $C^{\alg}(\bbf{G},K)$ be the ring of algebraic functions of $\bbf{G}$, i.e., the global sections of the affine group scheme $\bbf{G}_K$.  For $G_0\subset \bbf{G}(L)$ a compact open subgroup we define the space of locally algebraic functions of $G_0$ (relative to $\bbf{G}$) to be 
\[
C^{\lalg}(G_0,K)= C^{sm}(G_0,K)\otimes_{K} C^{\alg}(\bbf{G},K). 
\]
We let $\n{D}^{\lalg}(G_0,K)= \iHom_{K}(C^{\lalg}(G_0,K),K)$ be the locally algebraic distribution algebra of $G_0$ and for any $G_0 \subset G \subset \bbf{G}(L)$ an open subgroup we denote
\[ \n{D}^{\lalg}(G,K):= \n{K}_{\sol}[G]\otimes_{\n{K}_{\sol}[G_0]} \n{D}^{\lalg}(G_0,K) \] the locally algebraic distribution algebra of $G$. 

\begin{definition}
    Let $V\in \Mod_{\n{K}_{\sol}}(\n{D}^{\lalg}(G,K))$.
    \begin{enumerate}
        \item We let $C^{\lalg}(G,V):= \prod_{g\in G/G_0} (C^{\lalg}(gG_0,K)\otimes_{\n{K}_{\sol}}^L V)$ be the space of locally algebraic functions of $G$ with values in $V$.  The space $C^{\lalg}(G,V)$ has three commuting actions of $\n{D}^{\lalg}(G,K)$ given by the left $\star_1$ and right $\star_2$ regular actions, and the action $\star_3$ on $V$.

        \item Define the functor of locally algebraic vectors $(-)^{R\lalg}: \Mod_{\n{K}_{\sol}}(\n{D}^{\lalg}(G,K)) \to \Mod_{\n{K}_{\sol}}(\n{D}^{\lalg}(G,K)) $ to be 
        \[
        V^{R\lalg}:= R\iHom_{\n{D}^{\lalg}(G,K)}(K, C^{\lalg}(G, V)_{\star_{1,3}})
        \]
        endowed with the $\star_{2}$-action of $\n{D}^{\lalg}(G,K)$. 

        \item We say that an object $V\in \Mod_{\n{K}_{\sol}}(\n{D}^{\lalg}(G,K))$ is locally algebraic if the natural map $V^{R\lalg}\to V$ is an equivalence. We let $\Rep^{\lalg}_{\n{K}_{\sol}}(G)\subset \Mod_{\n{K}_{\sol}}(\n{D}^{\lalg}(G,K))$ be the full subcategory of locally algebraic representations. 
    \end{enumerate}
\end{definition}

\begin{lemma}
\label{LemmaLocAlgDescription}
    Let $G$ be a compact open subgroup of $\bbf{G}(L)$. We have natural isomorphisms of $\n{D}^{\lalg}(G^2,K)$-modules (for the actions $\star_1$ and $\star_2$)
    \[
C^{\lalg}(G,K)= \bigoplus_{\pi,\lambda} (\pi \otimes V^{\lambda} )_{\star_1} \otimes (\pi \otimes V^{\lambda})^{\vee}_{\star_2}
    \]
    and
    \[
\n{D}^{\lalg}(G,K):= \prod_{\pi, \lambda} (\pi \otimes V^{\lambda})^{\vee}_{\star_1} \otimes  (\pi\otimes V^{\lambda})_{\star_2} 
    \]
   where $\pi$ runs over all the smooth irreducible representations of $G$, and $V^{\lambda}$ over all the irreducible algebraic representations of $\bbf{G}$. 
\end{lemma}
\begin{proof}
    This follows from Lemma \ref{CoroSmoothAlgebraMatrix} and \cite[Theorem 4.2.7]{GoodmanWallach} describing the algebra of functions of $\bbf{G}$ in terms of irreducible algebraic representations. That is, taking duals on Lemma \ref{CoroSmoothAlgebraMatrix}  implies that the smooth regular representation $C^{sm}(G,K)$, seen as by $\n{D}^{sm}(G,K)$-module,  is the direct sum of $\pi\otimes \pi^{\vee}$ where $\pi$ runs over all irreducible finite dimensional smooth representations of $G$. Similarly, \cite[Theorem 4.2.7]{GoodmanWallach} states that the regular algebraic representation $C^{\mathrm{alg}}(\mathbf{G}, K)$ is isomorphic, as $\mathbf{G}\times \mathbf{G}$-representation, to the direct sum of $V^{\lambda}\otimes  V^{\lambda,\vee}$ where $V^{\lambda}$ runs over all the irreducible algebraic representations of $\mathbf{G}$. 
\end{proof}

\begin{proposition}
    The following assertions hold.
    \begin{enumerate}
        \item Let $V\in \Mod_{\n{K}_{\sol}}(\n{D}^{\lalg}(G,K))$, the natural map $(V^{R\lalg})^{R\lalg}\to V^{R\lalg}$ is an equivalence. 

        \item The functor $(-)^{R\lalg}$ commute with colimits.

        \item  Let $V,W\in \Mod_{\n{K}_{\sol}}(\n{D}^{\lalg}(G,K))$, then $(V^{R\lalg}\otimes_{\n{K}_{\sol}}^L W)^{R\lalg} = V^{R\lalg}\otimes_{\n{K}_{\sol}}^L W^{R\lalg}$. In particular, $\Rep^{\lalg}_{\n{K}_{\sol}}(G)$ has a natural symmetric monoidal structure. 

        \item The functor $(-)^{R\lalg}$ is the right adjoint of the inclusion $\Rep^{\lalg}_{\n{K}_{\sol}}(G) \subset \Mod_{\n{K}_{\sol}}(\n{D}^{\lalg}(G,K))$. 

         \item The functor $(-)^{\lalg} := (-)^{R^0\lalg}$ is exact in the abelian category $\Mod^{\heartsuit}_{\n{K}_{\sol}}(\n{D}^{\lalg}(G,K))$. In particular, $\Rep^{\lalg}_{\n{K}_{\sol}}(G)$ has a natural $t$-structure. 

         \item The $\infty$-category $\Rep^{\lalg}_{\n{K}_{\sol}}(G,K)$ is the derived category of its heart. 
    \end{enumerate}
\end{proposition}

\begin{proof}
This follows the same arguments of Propositions \ref{PropositionStableColimitsAndTensor}, \ref{PropPropertiesCategoryLocAn} and \ref{PropReplaGrothendieck} in the locally analytic case, or the Propositions \ref{PropositionRepsmoothGrothendieck} and \ref{PropAdjunctionSmooth} in the smooth case.   We give a sketch for completeness.  Let $G_0\subset G$ be an open compact subgroup, by adjunction  we have that 
\[
W^{R\lalg} = R\iHom_{\n{D}^{\lalg}(G,K)}(K, C^{\lalg}(G,W)) = R\iHom_{\n{D}^{\lalg}(G_0,K)}(K, C^{\lalg}(G_0,W)),
\]
then for (1)-(3) and (5) we can assume that $G$ is compact. By Lemma \ref{LemmaLocAlgDescription} any finite dimensional representation of $G$ is a direct summand of $\n{D}^{\lalg}(G, K)$, in particular they are projective.  This implies that $(-)^{R\lalg}$ is an exact functor in the abelian category and that it commutes with colimits. Moreover, we have that 
\[
\begin{aligned}
W^{R\lalg} & = R\iHom_{\n{D}^{\lalg}(G,K)}(K, (C^{\lalg}(G,K)\otimes_{\n{K}_{\sol}}^L W)_{\star_{1,3}}) \\ 
& = \varinjlim_{\pi,\lambda} R\iHom_{\n{D}^{\lalg}(G,K)}(K, (\pi\otimes V^{\lambda})\otimes  ((\pi\otimes V^{\lambda})^{\vee} \otimes W )_{\star_{1,3}} ) \\ 
& = \varinjlim_{\pi,\lambda} R\iHom_{\n{D}^{\lalg}(G,K)}((\pi \otimes V^{\lambda})^{\vee}, W) \otimes (\pi\otimes V^{\lambda})^{\vee}. 
\end{aligned}
\]
Then, to prove that the functor $(-)^{\lalg}$ is idempotent, by Lemma \ref{LemmaLocAlgDescription}, it suffices to prove it for the representations of the form $W=(\pi\otimes V^{\lambda})^{\vee}$. But the $\n{D}^{\lalg}(G,K)$-modules$\pi\otimes V^{\lambda}$ are irreducible and projective (being a direct summand of $\n{D}^{\lalg}(G, K)$). Using the above formula, we then deduce that, for $W = (\pi \otimes V^\lambda)^\vee$,
\[ W^{R \lalg} =  \varinjlim_{\pi',\lambda'} \iHom_{\n{D}^{\lalg}(G,K)}((\pi' \otimes V^{\lambda'})^{\vee}, W) \otimes (\pi' \otimes V^{\lambda'})^{\vee} = W, \]
proving that $(-)^{\lalg}$ is idempotent.
So far we have proven parts (1), (2) and (5). For part (3) we can assume that $W=C^{\lalg}(G,K)$ in which case we can untwist the diagonal action of $C^{\lalg}(G,K)\otimes C^{\lalg}(G,K) \otimes V$ to a representation where $\n{D}^{\lalg}(G,K)$ acts trivially on the first factor. Taking invariants by $\n{D}^{\lalg}(G,K)$ one gets that 
\[
(C^{\lalg}(G,K)\otimes W)^{R\lalg}= C^{\lalg}(G,K)\otimes W^{R\lalg}.
\]
Parts (4) and (6) follow the same lines of their analogues for smooth representations, see Propositions \ref{PropositionRepsmoothGrothendieck} and \ref{PropAdjunctionSmooth}.    
\end{proof}

\section{Adjunctions and cohomology}
\label{s:AdjunctionsCoho}

In this final section, we show how the cohomology comparison theorems of \cite[\S 5.2]{SolidLocAnRR} are explained in terms of adjunctions. 

\subsection{Geometric solid representations}

Following the interpretation of the categories of locally analytic and smooth representations as quasi-coherent sheaves of ``classifying stacks of $G^{la}$ and $G^{sm}$'', one can introduce a different category of ``continuous geometric'' representations where now $G$ is the analytic space defined by the algebra of its continuous functions. In this section, we will fix a $p$-adic Lie group $G$ over $\Q_p$.

\begin{definition}
Let $G_0\subset G$ be an open compact subgroup.
    \begin{enumerate}
        \item Let $V\in \Mod(\n{K}_{\sol})$, we define the space of ``geometric continuous'' functions of $G$ with values in  $V$ to be 
        \[
        C^{geo}(G,V):= \prod_{g\in G/G_0}(C(gG_0,K)\otimes_{\n{K}_{\sol}}^L V).
        \]
    
        \item We define the category of quasi-coherent sheaves of  the  underlying locally profinite group $G^{cont}$ to be  $\Mod^{qc}_{\n{K}_\sol}(G^{cont})=\prod_{g\in G/G_0} \Mod_{\n{K}_{\sol}}(C(gG_0,K))$. 
        \item We define the category of ``continuous geometric'' representations of $G$ to be the limit 
        \[
        \Rep^{geo}_{\n{K}_{\sol}}(G):= \Mod^{qc}_{\n{K}_\sol}([*/G^{cont}]):= \varprojlim_{[n] \in \Delta} \Mod^{qc}_{\n{K}_\sol}(G^{cont, n}).
        \]
    \end{enumerate}
\end{definition}

\begin{lemma}
\label{LemmaInjectionGEoInCont}
    Let $V\in \Mod^{\heartsuit}(\n{K}_{\sol})$ and $S$ a profinite set. Then the natural map $C(S, K)\otimes_{\n{K}_{\sol}} V \to C(S,V)$ is an injection.
\end{lemma}

\begin{proof}
    It is enough to take $K=\bb{Q}_p$. Since any solid  $\bb{Q}_p$-vector space is a colimit of quotients of compact projective $\bb{Q}_{p, \sol}$-vector spaces, we can assume that $V$ fits in a short exact sequence $0 \to  \bb{Q}_{p,\sol}[S] \to \bb{Q}_{p,\sol}[S'] \to V \to 0$. Taking lattices $0\to \bb{Z}_{p}[S] \to \bb{Z}_{p}[S'] \to Q \to 0$ (after rescaling if necessary), it suffices to show that the map 
    \[
    C(S,\bb{Z}_p)\otimes_{\bb{Z}_p} Q \to C(S,Q)
    \]
    is injective.  But both objects are $p$-adically complete, so it suffices to show that their reduction modulo $p^n$ are injective, i.e. that we have monomorphisms 
    \[
    C^{sm}(S, Q/p^n) \to C(S,Q/p^n).
    \]
    This is Lemma 3.4.8 (iii) of \cite{MannSix}.
\end{proof}

\begin{lemma}
\label{LemmaComodulesGeoRep}
    Let $\n{A}$ be the  category of comodules $V\to C^{geo}(G,V)$ with $V\in \Mod^{\heartsuit}(\n{K}_{\sol})$. Then $\n{A}$ is a Grothendieck abelian full  subcategory  of $\Mod^{\heartsuit}_{\n{K}_{\sol}}(\n{K}_{\sol}[G])$ with left completed derived $\infty$-category naturally equivalent to   $\Rep_{\n{K}_{\sol}}^{geo}(G)$.  
\end{lemma}
\begin{proof}
    The fact that $\n{A}$ is an abelian category follows from the fact that $V\mapsto C^{geo}(V,K)$ is an exact functor. We have a natural functor $\n{A} \to \Mod^{\heartsuit}_{\n{K}_{\sol}}(\n{K}_{\sol}[G])$ sending the comodule $V$ to the representation defined by the orbit map $V\to C^{geo}(G,V) \to C(G,V)$. It is clear that for $V,W\in \n{A}$  one has $\Hom_{\n{A}}(V,W)\subset \Hom_{\n{K}_{\sol}[G]}(V,W)$. Conversely, let $f:V\to W$ be a morphism of $\n{K}_{\sol}[G]$-modules. We have a  diagram  whose ambient square is commutative
    \[
    \begin{tikzcd} 
    V \ar[d] \ar[r,"f"]& W \ar[d] \\
         C^{geo}(G,V) \ar[d] \ar[r] & C^{geo}(G,W) \ar[d] \\ 
        C(G,V) \ar[r] & C(G, W)
    \end{tikzcd}
    \]
    and such that the lower vertical arrows are injective by Lemma \ref{LemmaInjectionGEoInCont}, then the  upper square must be commutative proving that $\Hom_{\n{A}}(V,W) = \Hom_{\n{K}_{\sol}[G]}(V,W)$.

To finish the proof of the lemma, we  need to show the following facts:

\begin{enumerate}

\item The category $\Rep^{geo}_{\n{K}_{\sol}}(G)$ has a natural $t$-structure making the pullback along $*\to */G^{cont}$ a $t$-exact functor. 

\item  The heart $\Rep^{geo,\heartsuit}_{\n{K}_{\sol}}(G)$ is a Grothendieck abelian category. 

\item $\Rep^{geo}_{\n{K}_{\sol}}(G)$ is the left completion of the derived category its heart. 

\item There is a natural equivalence of symmetric monoidal categories 
\[
\s{A}= \Rep^{geo,\heartsuit}_{\n{K}_{\sol}}(G). 
\]

\end{enumerate}
    
These facts are proven in the exact same way as in Theorem \ref{TheoremLocAnRepStack}, we give a summary of the main points of the proof for completeness. Let $G_0\subset G$ be a compact open subgroups and consider the maps of analytic stacks 
\[
\AnSpec \n{K}_{\sol}\xrightarrow{f} [\AnSpec \n{K}_{\sol}/ G^{cont}_0]\xrightarrow{h} [\AnSpec \n{K}_{\sol}/ G^{cont}]
\]     
Since  $G^{cont}_0=\AnSpec C(G_0,K)$ is affinoid with induced analytic structure from $\n{K}_{\sol}$, the map $f$ is prim. This implies that $f^*$ admits a colimit preserving right adjoint $f_*$ which satisfies the projection formula and is compatible with base change, as $f$ satisfies $*$-descent one deduces from \cite[Theorem 4.7.5.2]{HigherAlgebra} that $f^*$ is comonadic. By looking at the pullback square 
\[
\begin{tikzcd}
G_0^{cont}  \ar[r] \ar[d] & \AnSpec \n{K}_{\sol}  \ar[d] \\ 
\AnSpec \n{K}_{\sol} \ar[r] &  {[\AnSpec \n{K}_{\sol}/ G_0^{cont}]}
\end{tikzcd}
\]
one deduces that the underlying functor of the comonad $f^*f_*$ is just tensoring with $C(G_0,K)$, and in particular $t$-exact.   Then, Lemma \ref{LemmaDerivedCategoriesHeart} (1) and (3) shows that $\Rep^{geo}_{\n{K}_{\sol}}(G_0)$ has a natural $t$-structure making $f^*$ a $t$-exact map, that its heart is a Grothendieck abelian category, and that it is the left completion of the derived category of its heart, namely, points (1)-(3) above.  

One then extends the properties (1)-(3) from $G_0$ to $G$ using the map  $h$ instead: this map is \'etale since the quotient $G^{cont}/G^{cont}_0=\bigsqcup_{g\in G/G_0} \AnSpec \n{K}_{\sol}$ is a finite disjoint union of points over $\AnSpec \n{K}_{\sol}$. This guarantees that $h^*$ admits a left adjoint $h_{\sharp}$ satisfying projection formula and compatible with base change, since $h$ satisfies $*$-descent one deduce by \cite[Theorem 4.7.5.2]{HigherAlgebra} that $h^*$ is monadic. One shows that the monad $h^*h_{\sharp}$ is $t$-exact and colimit preserving, which by   Lemma \ref{LemmaDerivedCategoriesHeart} (1) and (2) implies that $ \Rep^{geo}_{\n{K}_{\sol}}(G)$ satisfies the points (1)-(3) above.  

Finally, to conclude the lemma, it suffices to identify the heart  $\Rep^{geo}_{\n{K}_{\sol}}(G)$ with the category $\s{A}$ of comodules for $C^{geo}(G,-)$. This is proven exactly as in Theorem \ref{TheoremLocAnRepStack} by explicitly describing the abelian descent datum along the diagram
\[
\begin{tikzcd}
G^{2,cont} \ar[r, shift left= 2ex] \ar[r, shift right  =2ex] \ar[r] & G^{cont} \ar[r,shift left =1 ex] \ar[r, shift right =1ex]  \ar[l, shift left = 1ex] \ar[l, shift right = 1ex]& * \ar[l].
\end{tikzcd}
\]
\end{proof}

We have a natural morphism of coalgebras $C^{la}(G,K)\to C(G,K)$ which heuristically should induce a group homomorphism $G^{cont} \to G^{la}$ and as consequence a morphism of their classifying stacks $f:[*/G^{cont}] \to [*/ G^{la}]$. We can define a pullback functor $f^*: \Mod^{qc}_{\n{K}_\sol}([*/G^{la}]) \to \Mod^{qc}_{\n{K}_\sol}([*/G^{cont}])$ which corresponds to a forgetful functor $F: \Rep^{la}_{\n{K}_{\sol}}(G) \to \Rep^{geo}_{\n{K}_{\sol}}(G)$ sending the comodule $V\to C^{la}(G,V)$ to the comodule $V\to C^{la}(G,V)\to C^{geo}(G,V)$. The functor $f^*$ preserves colimits, so it admits a right adjoint that we can call the pushforward $f_*: \Mod^{qc}_{\n{K}_\sol}([*/G^{cont}])  \to \Mod^{qc}_{\n{K}_\sol}([*/G^{la}])$. At the level of representations we can think of $f_*$ as a locally analytic vectors functor $(-)^{Rla}: \Rep^{geo}_{\n{K}_{\sol}}(G) \to \Rep^{la}_{\n{K}_{\sol}}(G)$. 

\begin{definition}
    We define the ``continuous geometric'' cohomology $R\Gamma^{geo}(G,-):  \Rep^{geo}_{\n{K}_{\sol}}(G) \to \Mod(\n{K}_{\sol})$ to be the right adjoint of the trivial representation functor $\Mod(\n{K}_{\sol}) \to \Rep^{geo}_{\n{K}_{\sol}}(G)$. 
\end{definition}

We have the following proposition. 

\begin{prop}
    The forgetful functor $f^*: \Mod^{qc}_{\n{K}_\sol}([*/G^{la}]) \to \Mod^{qc}_{\n{K}_\sol}([*/G^{cont}])$ is fully faithful. The right adjoint of $f^*$ on a geometric representation $V$ can be computed as 
    \[
    f_*V =  R\Gamma^{geo}(G, C^{la}(G,V)_{\star_{1,3}}),
    \]
as $\n{D}^{la}(G,K)$-modules, via the inclusion  $\Mod^{qc}_{\n{K}_\sol}([*/G^{la}]) = \Rep^{la}_{\n{K}_{\sol}}(G) \subset \Mod_{\n{K}_\sol}(\n{D}^{la}(G,K)) $,    where $\n{D}^{la}(G,K)$ acts on the right term via the right regular action. 
\end{prop}
\begin{proof}

    By Lemma  \ref{LemmaComoduleAndRepLa} the category $\Rep^{la}_{\n{K}_{\sol}}(G)$ is the derived category of comodules of the functor $C^{la}(G,-)$. Similarly, by Lemma \ref{LemmaComodulesGeoRep} the category $\Rep^{geo}_{\n{K}_{\sol}}(G)$ is the derived category of the abelian category of comodules of $C^{geo}(G,-)$. Moreover, we have fully faithful inclusions of abelian categories  $\Rep^{la,\heartsuit}_{\n{K}_{\sol}}(G) \subset \Rep^{geo,\heartsuit}_{\n{K}_{\sol}}(G) \subset \Mod^{\heartsuit}(\n{K}_{\sol}[G])$. This implies that the right adjoint of the first inclusion is given by the locally anlaytic vectors functor that can be computed as $C^{la}(G,V)^{G_{\star_{1,3}}}$. Taking right derived functors we see that $f_*V= R\Gamma^{geo}(G, C^{la}(G,V))$ for any $V\in \Rep^{geo}_{\n{K}_{\sol}}[G]$, where $\n{D}^{la}(G,K)$ acts on the right term via the $\star_2$-action.

    It is left to show that the unit map $1\to f_*f^*$ is an equivalence.  Let $G_0\subset G$ be a compact open subgroup, we have a commutative diagram of morphisms of stacks 
    \[
    \begin{tikzcd}
  {[}*/G^{cont}_0{]}  \ar[r,"\tilde{f}"] \ar[d, "\tilde{g}"]  &  {[}*/G^{la}_0{]}  \ar[d,"g"]\\ 
 {[}*/G^{cont}{]}   \ar[r,"f"]   & {[}*/G^{la}{]}.
    \end{tikzcd}
    \]
    The pullback functors correspond to forgetful functors, and the vertical pushforward functions are given by inductions. Indeed, we can check this at the level of abelian categories where the right adjoint of a forgetful functor is clearly an induction. As a  consequence one deduces that 
    \[
    R\Gamma^{geo}(G, C^{la}(G,V)) = R\Gamma^{geo}(G_0, C^{la}(G_0,V)).
    \]
    Thus, we can assume without loss of generality that $G$ is compact. In this case $C^{la}(G,V)=C^{la}(G,K)\otimes^{L}_{\n{K}} V$ and $C^{geo}(G,V)=C(G,K)\otimes^{L}_{\n{K}} V$.

 Notice that for $V\in \Rep^{la}_{\n{K}_{\sol}}(G)$ we have a natural equivalence of representations $C^{la}(G,V)_{\star_{1,3}} \xrightarrow{\sim }   C^{la}(G,V)_{\star_{2}}$. Indeed, by taking derived functors it suffices to produce such equivalence for $V$ in the heart, in that case the equivalence arises from the orbit map  (which is $\n{D}^{la}(G,K)$-equivariant)
 \[
O_V\colon  V\to  C^{la}(G,V)_{\star_2}
 \]
which induces a $C^{la}(G,K)_{\star_2}$-linear isomorphism  
\[
C^{la}(G,V)_{\star_{1,2}}=C^{la}(G,K)\otimes_{\n{K}_{\sol}} V \xrightarrow{\sim} C^{la}(G,V)_{\star_2}
\]
which is $\n{D}^{la}(G,K)$-equivariant for the $\star_{1,2}$-action on the left, and the $\star_2$-action on the right. Informally, this isomorphism sends a locally analytic function $f\colon G\to V$ to the function $\widetilde{f}\colon G\to V$ with $\widetilde{f}(g) =g \cdot f(g)$. Applying the antiinvolution on $G$,  one obtains a $\n{D}^{la}(G,K)$-equivariant  isomorphism  $C^{la}(G,V)_{\star_{1,3}}\cong C^{la}(G,V)_{\star_{1,2}}$ which yields the equivalence $C^{la}(G,V)_{\star_{1,3}} \xrightarrow{\sim }   C^{la}(G,V)_{\star_{2}}$ as desired.

  Thus, it suffices to show that for a trivial representation $V$ one has $R\Gamma^{geo}(G, C^{la}(G,V)_{\star_{2}})=V$, or equivalently by applying the antiinvolution of $G$, that $R\Gamma^{geo}(G,C^{la}(G,V)_{\star_{1}})=V$. Writing $V$ as limit of canonical and stupid truncations we can assume that $V$ is a solid $\n{K}_{\sol}$-vector space in degree $0$. But by Proposition \ref{PropComparisonCochainsAndCoho} down below one can compute this geometric cohomology using  geometric cochains, i.e. $R\Gamma^{geo}(G, C^{la}(G,V)_{\star_1})$ is represented by the bar complex of geometric cochains 
 \[
 [C^{la}(G,V)\to C^{geo}(G,  C^{la}(G,V)) \to C^{geo}(G^2, C^{la}(G,V)) \to \cdots],
 \]
which is the same as the tensor product of the bar complex
\[
[C^{la}(G,K) \to C(G,K)\otimes^L_{\n{K}} C^{la}(G,K) \to \cdots ]\otimes_{\n{K}_{\sol}}^L V.
\]
But $C^{la}(G,K)$ is a nuclear $\n{K}_{\sol}$-vector space, so that the geometric bar complex  of $C^{la}(G,K)$ (obtained by cochains of the form $C(G^{\bullet},K)\otimes_{\n{K}_{\sol}}C^{la}(G,K)$) is equal to  the solid bar complex (obtained by cochains of the form $C(G^{\bullet},C^{la}(G,K))=\iHom_K(\n{K}_{\sol}[G^{\bullet}],C^{la}(G,K))$) which computes $R\Hom_{\n{K}_{\sol}[G]}(K, C^{la}(G,K))=R\Hom_{\n{D}^{la}(G,K)}(K, C^{la}(G,K))=K$, where we used Corollary \ref{coroIdemGQp} for the first equivalence. This finishes the proof. 
\end{proof}

\begin{remark}
Using the  six-functor formalism for analytic stacks the previous proof simplifies considerably. Let $f: [*/G^{geo}] \to [*/G^{la}]$ be the natural map of stacks, it suffices to prove that the natural map $\mathrm{id}\to f_*f^*$ is an equivalence.  The map $f$ is prim  as the fibers are isomorphic to $[G_0^{la}/G_0^{geo}]$ for $G_0$ any compact open subgroup. Thus, the pullback square 
\[
\begin{tikzcd}
{[G^{la}/G^{cont}]} \ar[r] \ar[d] & {[*/G^{cont}]} \ar[d] \\ 
{[*/G^{cont}]}\ar[r] & {[*/ G^{la}]}
\end{tikzcd}
\]
and proper base change for $f_*$ implies that it is computed as $R\Gamma^{geo}(G, C^{la}(G,V)_{\star_{1,3}})$. Moreover,  by the projection formula, in order to prove    fully faithfulness of $f^*$ it suffices to show that the natural map    $1_{[*/G^{la}]}\to f_* 1_{[*/G^{geo}]}$ is an equivalence, this follows from the explicit computation using the bar complexes and the previous explicit description of $f_*$. 
\end{remark}

\subsection{Adjunctions}
\label{ss:AdjunctionsRepresentations}

Let $G$ be as always a $p$-adic Lie group defined over a finite extension $L$ of $\qp$,  $ \n{K}_{\sol}=(K, K^+)$ a complete non-achimedean extension of $L$ and $\n{K}_\sol = (K, K^+)_\sol$. To avoid any confusion, when talking about locally analytic representations, in this section we will note $G_L = G$ to stress that we see the group $G$ defined over $L$ and we denote by $G_\qp$ the $p$-adic Lie group $G$ viewed over $\qp$. For continuous and smooth representations this disctinction is unnecessary since their definition is independent of the $L$-analytic structure, and we will simply use the notation $G$. We have the following diagram of categories.

\begin{equation} \label{EquationDiagramRepresentationCategories}
    \Rep^{sm}_{\n{K}_\sol}(G) \xrightarrow{F_1} \Rep^{la}_{\n{K}_\sol}(G_L) \xrightarrow{F_2} \Rep^{la}_{\n{K}_\sol}(G_\qp) \xrightarrow{F_3}  \Rep_{\n{K}_\sol}(G),
\end{equation}
where $\Rep_{\n{K}_\sol}(G) = \Mod_{\n{K}_\sol}(\n{K}_\sol[G])$ denotes the category of solid representations of $G$, and where the natural functors $F_i$ are just the forgetful functors. Since all these functors commute with colimits, they all have right adjoints and the purpose of this section is to calculate each of them.

\begin{proposition} \label{PropositionAdjunctions} \leavevmode
\begin{enumerate}
    \item The right adjoint of $F_1$ is given by Lie algebra cohomology $R \Gamma(\f{g}_L, -) := R \iHom_{U(\f{g}_L)}(K, -)$.
    \item The right adjoint of $F_2$ is given by $R \Gamma(\f{k}, -) := R \iHom_{U(\f{k})}(K, -)$, where $\f{k} = \ker(\f{g}_\qp \otimes_\qp L \to \f{g}_L)$.
    \item The right adjoint of $F_3$ is given by the functor of locally analytic vectors $(-)^{Rla}$.
\end{enumerate}
\end{proposition}

\begin{proof}
Let $V \in \Rep^{la}(G_L)$ and $W \in \Rep^{la}(G_\qp)$. Then
\begin{align*}
R \iHom_{\n{D}^{la}(G_\qp, K)}(V, W) &= R \iHom_{\n{D}^{la}(G_\qp, K)}( \n{D}^{la}(G_L, K)\otimes_{\n{D}^{la}(G_L, K)}^L  V, W) \\
&= R \iHom_{\n{D}^{la}(G_L, K)}(V, R \iHom_{\n{D}^{la}(G_\qp, K)}( \n{D}^{la}(G_L, K), W) \\
&= R \iHom_{\n{D}^{la}(G_L, K)}(V, R \iHom_{\n{D}^{la}(G_\qp, K)}( \n{D}^{la}(G_{\Q_p}, K) \otimes_{\n{D}^{la}(\f{k}, K)} K, W)) \\ 
&= R \iHom_{\n{D}^{la}(G_L, K)}(V, R \iHom_{\n{D}^{la}(\f{k}, K)}( K, W)) \\ 
&= R \iHom_{\n{D}^{la}(G_L, K)}(V, R \iHom_{U(\f{k})}(K, W)),
\end{align*}
where the first two and the fourth equalities are trivial, the third one follows Lemma \ref{LemmmaBaseChangeDistAlgebras} (1), and the last one by idempotency of $\n{D}^{la}(\f{k}, K)$ over $U(\f{k})$, which follows from Lemma \ref{PropKoszulResolutionDist}. This proves $(2)$.

Recall from Lemma \ref{PropPropertiesSmoothDistributionAlgebrab} that $\n{D}^{sm}(G_L, K) = K \otimes^L_{\n{D}^{la}(\f{g}_L,K)} \n{D}^{la}(G_L, K)$. Then, using the exact same argument as in the proof of $(2)$, we have, for $V \in \Rep^{sm}_{\n{K}_\sol}(G)$ and $W \in \Rep^{la}_{\n{K}_\sol}(G_L)$,
\[ R \iHom_{\n{D}^{la}(G_L, K)}(V, W) = R \iHom_{\n{D}^{sm}(G_L, K)}(V, R \iHom_{\n{D}^{la}(\f{g}_L, K)}(K, W)), \]
proving $(1)$.

By Corollary \ref{CoroRightAdjointLocAnGroup}, the right adjoint to the fully faithful inclusion $\Rep^{la}_{\n{K}_\sol}(G_{\bb{Q}_p}) \to \Mod_{\n{K}_\sol}(\n{D}^{la}(G_{\bb{Q}_p},L))$ is given by the functor $(-)^{Rla}$. Since the (fully faithful) inclusion $\Mod_{\n{K}_\sol}(\n{D}^{la}(G_{\bb{Q}_p}, K)) \to \Mod_{\n{K}_\sol}(\n{K}_\sol[G])$ has a right adjoint given by $R \iHom_{\n{K}_\sol[G]}(\n{D}^{la}(G_{\bb{Q}_p}, K), -)$, the third assertion follows formally: 
\begin{align*} 
\big( R \iHom_{\n{K}_\sol[G]}(\n{D}^{la}(G_{\bb{Q}_p}, K), W) \big)^{Rla}  & =\varinjlim_h R\iHom_{\n{D}^{la}(G_{\bb{Q}_p}, K)}(\n{D}^{h}(G,K),  \iHom_{\n{K}_\sol[G]}(\n{D}^{la}(G_{\bb{Q}_p}, K),W))   \\ 
&= \varinjlim_h R\iHom_{\n{K}_{\sol}[G]}(\n{D}^{h}(G,K),W)  \\
&=W^{Rla}
\end{align*}
where the first and third equivalences follow from the definition,  and the second equivalence form an adjunction. 
\end{proof}

\begin{remark}
Consider the following sequence  of adjunctions 
\[
\begin{tikzcd}
    \Rep^{sm}_{\n{K}_\sol}(G) \ar[r, "F_1", shift left= 1ex,rightharpoonup] &  \Rep^{la}_{\n{K}_\sol}(G_L) \ar[r, "F_2",shift left= 1ex,rightharpoonup] \ar[l, "R\Gamma(\f{g}_L{,}-) ",  shift left= 1ex,rightharpoonup] &  \Rep^{la}_{\n{K}_\sol}(G_\qp ) \ar[r, "F_3",shift left= 1ex,rightharpoonup] \ar[l, "R\Gamma(\f{k}{,}-) ",  shift left= 1ex,rightharpoonup]  &   \Rep_{\n{K}_\sol}(G) \ar[l, "(-)^{Rla} ",  shift left= 1ex,rightharpoonup].
\end{tikzcd}
\]
One can define functors of smooth or locally analytic vectors from different categories of representations as right adjoint of forgetful functors. For example, let $F$ be the composite forgetful functor $\Rep^{sm}_{\n{K}_{\sol}}(G) \to \Rep_{\n{K}_{\sol}}(G)$, then its right adjoint can be computed as the composite of the right adjoints of the forgetful functors 
\[
\Rep^{sm}_{\n{K}_{\sol}}(G) \to \Mod_{\n{K}_{\sol}}(\n{D}^{sm}(G,K)) \to \Mod_{\n{K}_{\sol}}(\n{K}_{\sol}[G]) =\Rep_{\n{K}_{\sol}}(G).
\]
This can be computed by applying simple adjunctions as follows: if $V \in \Rep^{sm}_{\n{K}_\sol}(G)$ and $W \in  \Rep_{\n{K}_\sol}(G)$, then 
    \begin{align*}
        R \iHom_{\n{K}_\sol[G]}(V, W) &= R \iHom_{\n{D}^{sm}(G,K)}( V, R \iHom_{\n{K}_\sol[G]}(\n{D}^{sm}(G, K), W) ) \\
        &= R \iHom_{\n{D}^{sm}(G,K)}( V, \big(R \iHom_{\n{K}_\sol[G]}(\n{D}^{sm}(G, K), W) \big)^{Rsm} ) 
    \end{align*}
    where the first equality follows using the fact that $V$ is a $\n{D}^{sm}(G, K)$-module and adjunction, the second one by Proposition \ref{PropAdjunctionSmooth}. Thus, the right adjoint of $F$ is 
    \begin{align*}
    W^{Rsm} & :=\big( R\iHom_{\n{K}_\sol[G]}(\n{D}^{sm}(G, K), W) \big)^{Rsm} \\ 
     & =  \varinjlim_{H\subset G} R\iHom_{\n{D}^{sm}(G,K)}(\n{K}_{\sol}[G/H], R\iHom_{\n{K}_\sol[G]}(\n{D}^{sm}(G, K), W)) \\
    & = \varinjlim_{H\subset G} R\iHom_{\n{K}_{\sol}[G]}(\n{K}_{\sol}[G/H], W).
    \end{align*}
\end{remark}

\subsection{Cohomology and comparison theorems}

We now introduce all the cohomology theories we are interested in, namely, Lie algebra, smooth, locally $L$ and $\qp$-analytic, and solid group cohomologies. We will first define them and show that these definitions recover the usual ones at abelian level. Finally, we will show how they compare to each other by some formal adjunctions. 

There is a natural map from the category $\Mod(\n{K}_\sol)$ to each of the categories appearing in \eqref{EquationDiagramRepresentationCategories} (and to $\Mod_{\n{K}_\sol}(U(\f{g}))$) given by trivial representations.

\begin{definition}
We define
\begin{itemize}
\item Solid group cohomology $R \Gamma(G, -) : \Rep_{\n{K}_\sol}(G) \to \Mod(\n{K}_\sol)$,
\item ($\qp$-)Locally analytic group cohomology $R \Gamma^{la}(G_\qp, -) : \Rep^{la}_{\n{K}_\sol}(G_\qp) \to \Mod(\n{K}_\sol)$,
\item ($L$-)Locally analytic group cohomology $R \Gamma^{la}(G_L, -) : \Rep^{la}_{\n{K}_\sol}(G_L) \to \Mod(\n{K}_\sol)$,
\item Smooth group cohomology $R \Gamma^{sm}(G, -) : \Rep^{sm}_{\n{K}_\sol}(G) \to \Mod(\n{K}_\sol)$
\item Lie algebra cohomology $R \Gamma(\f{g}, -) : \Mod_{\n{K}_\sol}(U(\f{g})) \to \Mod(\n{K}_\sol)$,
\end{itemize}
as the right adjoint functor of the trivial representation functor from $\Mod(\n{K}_\sol)$ to the corresponding category.
\end{definition}

\begin{remark}
As the categories $\Rep^{sm}_{\n{K}_\sol}(G)$, $\Rep^{la}_{\n{K}_\sol}(G_L)$ and $\Rep^{la}_{\n{K}_\sol}(G_\qp)$ embed fully faithfully in the categories $\Mod_{\n{K}_\sol}(\n{D}^{sm}(G, K))$, $\Mod_{\n{K}_\sol}(\n{D}^{la}(G_L, K))$  and $\Mod_{\n{K}_\sol}(\n{D}^{la}(G_\qp, K))$ respectively, we also have that \begin{gather*} R \Gamma^{la}(G_\qp, V) = R \iHom_{\n{D}^{la}(G_\qp, K)}(K, V),  \\  R \Gamma^{la}(G_L, V) = R \iHom_{\n{D}^{la}(G_L, K)}(K, V), \\ R \Gamma^{sm}(G, V) = R \iHom_{\n{D}^{sm}(G, K)}(K, V). \end{gather*} Moreover, since the categories $\Rep^{sm}_{\n{K}_{\sol}}(G)$ and $\Rep^{la}_{\n{K}_{\sol}}(G_L)$ are the derived categories of their heart, the smooth  and locally analytic cohomology functors can be computed as the right derived functors of the $G$-invariants of  their respective representation categories. 
\end{remark}

By \cite[Corollary 3.4.17]{MannSix}, smooth cohomology can be computed using smooth cochains. We prove the same for geometric, solid and locally analytic representations. 

\begin{prop}
\label{PropComparisonCochainsAndCoho}
Let $\Rep^{?}_{\n{K}_{\sol}}(G)$ denote the category of smooth, $L$-locally analytic, geometric or solid representations of $G$, and let $R\Gamma^{?}(G,-)$ denote their corresponding cohomology functor.   Let $V\in \Rep^{?, \heartsuit}_{\n{K}_\sol}(G)$ be a representation in degree $0$  and let $[C^{?}(G^{\bullet},V), d^n]$ be the bar complex in $\Mod(\n{K}_{\sol})$  with $n$-th term $C^{?}(G^{n}, V)$ and  $n$-th boundary map 
\begin{align*}
d^n(f)(g_1,\ldots, g_{n+1}) &= g_1f(g_2,\ldots, g_{n+1}) \\
&+ \sum_{i=1}^{n}(-1)^{i} f(g_1,\ldots, g_{i-1}, g_ig_{i+1},  \ldots, g_{n+1}) \\
&+ (-1)^{n+1} f(g_1,\ldots, g_n).
\end{align*}
Then there is a natural equivalence 
\[
R\Gamma^{?}(G,V) = [C^{?}(G^{\bullet},V), d^
{\bullet}]. 
\]
\end{prop}
\begin{proof}
We follow the same proof of \cite[Lemma 3.4.15]{MannSix}. Let $?-\mathrm{Ind}: \Mod(\n{K}_{\sol})\to \Rep^{?}_{\n{K}_\sol}(G)$ be the right adjoint of the forgetful functor,  and let $r$ denote the composition of the forgetful functor of $\Rep^{?}_{\n{K}_{\sol}}(G)$ with $?-\mathrm{\Ind}$, for $n\geq 0$ we let $r^n(-)$ denote the application of $n$-times $r$. By adjunction, we have natural transformations $r^{n}(-)\to  r^{n+1}(-)$ for all $n\geq 0$. For $M\in \Rep^{?,\heartsuit}_{\n{K}_{\sol}}(G)$, we claim that the complex 
\begin{equation}
\label{eqComplexResolutionLocAn}
0\to M \to r(M)\to r^2(M)\to \cdots
\end{equation}
is exact and that $r^{n}(M)= C^{?}(G^n,M)$. First, we claim  that for any $W\in \Mod(\n{K}_{\sol})$ one has $?-\mathrm{Ind}(W)=C^{?}(G,W)$.  Note that since the forgetful functor is $t$-exact, its right adjoint is left exact and equal to  the right derived functor of its restriction to abelian level. Thus,   it suffices to take $W\in \Mod^{\heartsuit}(\n{K}_{\sol})$, in which case we need to compute the right adjoint of the forgetful functor of abelian categories $\Rep^{?,\heartsuit}_{\n{K}_{\sol}}(G)\to \Mod^{\heartsuit}(\n{K}_{\sol})$.  For ? solid one has $\Rep^{?}_{\n{K}_{\sol}}(G)= \Mod^{\heartsuit}(\n{K}_{\sol}[G])$ and the induction is just $C(G,V)$. For $?$ being  smooth, locally analytic or geometric, the category $\Rep^{?,\heartsuit}_{\n{K}_{\sol}}(G)$ is the category of comodules of the exact functor $C^{?}(G,-)$, and one easily checks that the right adjoint of the forgetful functor is simply $V\mapsto C^{?}(G,V)$ proving the claim. 

Now, unraveling the definitions, one has that the sequence \eqref{eqComplexResolutionLocAn} is given by the usual bar complex of the respective representation category, which is an exact complex as they are constructed functorially from the  augmented cosimplicial object $(G^{n+1})_{n\in \Delta^{op}} \xrightarrow{\epsilon} *$. To conclude the proof we need to show that $R\Gamma^{?}(G, ?-\mathrm{Ind}(W))=W$ for any $W\in \Mod(\n{K}_{\sol})$, but the functor $R\Gamma^{?}(G,?-\mathrm{Ind}(-))$ is the right adjoint of the composite of the trivial representation and the forgetful functor which is the identity on $\Mod(\n{K}_{\sol})$, so it is equivalent to the identity. This finishes the proof. 
\end{proof}

All our comparison results are subsumed in the following statement, which generalizes in particular our main results \cite[Theorem 5.3 and Theorem 5.5]{SolidLocAnRR} from the case of a compact $p$-adic Lie group defined over $\qp$ to that of a (non-necessarily compact) $p$-adic Lie group defined over a finite extension of $\qp$ .

\begin{theorem} \label{TheoremComparisonCohomologies}
We have the following commutative diagram:
\begin{small}
\[
\begin{tikzcd}
& \Rep^{la}_{\n{K}_\sol}(G_\qp) \ar[rr, bend left = 10, "R \Gamma(\f{k}{,} -)"] \arrow[dddr, "R \Gamma^{la}(G_\qp {,} -)" near start] & & \Rep^{la}_{\n{K}_\sol}(G_L) \ar[dr, bend left = 10, "R \Gamma(\f{g}{,} -)"] \ar[dddl, "R \Gamma^{la}(G_L{,}-)"' near start] & \\
\Rep_{\n{K}_\sol}(G) \ar[ur, bend left = 10, "(-)^{Rla}"] \arrow[ddrr, "R \Gamma(G{,} -)"'] & & & & \Rep^{sm}_{\n{K}_\sol}(G) \ar[ddll, "R \Gamma^{sm}(G{,}-)"] \\
& & & & \\
& & \Mod(\n{K}_\sol) & & 
\end{tikzcd}
\]
\end{small}
Moreover, since the embedding $\Rep^{la}_{\n{K}_\sol}(G_\qp)$ in $\Rep_{\n{K}_\sol}(G)$ is fully faithful then, for any $V \in \Rep_{\n{K}_\sol}(G)$, we have $R \Gamma(G, V) = R \Gamma(G_\qp, V^{Rla})$. In particular, we have that
\[ R \Gamma(G, V) = R \Gamma(G_\qp, V^{Rla}) = R \Gamma^{la}(G_L, V^{R L-la}) = R \Gamma^{sm}(G, R \Gamma(\f{g}, V^{R L-la})), \]
where $V^{R L-la} := R \Gamma(\f{k}, V^{Rla})$ are the $L$-locally analytic vectors of a solid representation of $G$ (see Remark \ref{RemarkLocAnvectorsSolidRep}).
\end{theorem}

\begin{proof}
It follows by the adjunctions of Proposition \ref{PropositionAdjunctions}.
\end{proof}

\subsection{Homology and duality}

We conclude with some applications to duality between cohomology and homology. We let $G$ be a  compact $p$-adic Lie group over $L$.  The following result is the infinitesimal analogue of \cite[Theorem 5.19]{SolidLocAnRR}.

\begin{proposition} \label{PropDualityInfinitesimal}
Let $V\in \Mod_{\n{K}_\sol}(U(\f{g}))$. Then we have
\[ R \Gamma(\f{g}, V) = K(\chi)[-d] \otimes^L_{U(\f{g})} V. \]
In particular, if $V \in \Rep_{\n{K}_\sol}(G)$, then
\[ R \Gamma(G, V) = R \Gamma^{sm}(G, K(\chi)[-d] \otimes^L_{U(\f{g})} V^{R L-la})). \]
\end{proposition}

\begin{proof}
This follows exactly the same argument as in \cite[Theorem 5.19]{SolidLocAnRR} replacing the Lazard-Serre resolution by the Chevalley-Eilenberg resolution to calculate cohomology. The last assertion follows from the first one and Theorem \ref{TheoremComparisonCohomologies}.
\end{proof}

\appendix

\section{Auxiliary results}

In this appendix, we collect some abstract results on idempotent algebras and Hopf algebras that are used all along the text.

\subsection{Idempotent algebras}

We will let $\n{K}=(K,K^+)$ denote a complete non-archimedean extension of $\bb{Q}_p$, and let $\n{K}_{\sol}=(K,K^+)_{\sol}$ be the analytic ring attached to the Huber pair as in \cite[\S 3.3]{Andreychev}. Given an algebra $D$ in $\Mod(\n{K}_{\sol})$, we  endow $D$ with the induced analytic ring structure from $\n{K}_{\sol}$, and let $-\otimes_{D}^L-$ (or sometimes $-\otimes_{D, \sol}^L-$) denote the relative tensor product of $D$-modules in $\n{K}_{\sol}$-vector spaces.

Idempotent maps of associative algebras are defined as follows:

\begin{lemma}\label{lemmaIdempotentAlgebras}
Let  $\n{C}$ be a presentable symmetric monoidal stable $\infty$-category, let $f:A\to B$ be a morphism of associative algebras in $\n{C}$. The following conditions are equivalent:

\begin{enumerate}
    \item  The multiplication map $B\otimes_{A}B\xrightarrow{m} B$ of $B$-bimodules is an isomorphism.

    \item  The map $B\xrightarrow{\mathrm{id} \otimes 1} B\otimes_{A} B$ of left $B$-modules is an isomorphism. 

    \item The map $B\xrightarrow{1\otimes \mathrm{id}} B\otimes_A B$ of right $B$-modules is an isomorphism.

    \item The forgetful map $\mathrm{LMod}_{B}(\n{C}) \to \mathrm{LMod}_{A}(\n{C})$ of left modules if fully faithful. 

       \item The forgetful map $\mathrm{RMod}_{B}(\n{C}) \to \mathrm{RMod}_{A}(\n{C})$ of right modules if fully faithful. 
    
\end{enumerate}

If any of the previous conditions holds we say that the morphism of algebra $A\to B$ is idempotent. 

\end{lemma}
\begin{proof}
Let $A\to B$ be a map of associative algebra in $\n{C}$, the morphism $B\xrightarrow{\mathrm{id} \otimes 1} B\otimes_A B$ (resp, $B\xrightarrow{1\otimes \mathrm{id}} B\otimes_A B$) is a section of the multiplication map $B\otimes_A B \xrightarrow{m} B$, thus conditions (1), (2) and (3) are equivalent. 

We now prove that (1) is equivalent to (4); one  proves that (1) is equivalent to (5) in the same way. Suppose that $m: B\otimes_A B \to B$ is an equivalence. The forgetful map 
\[
G: \mathrm{LMod}_{B}(\n{C})\to \mathrm{LMod}_{A}(\n{C})
\]
has a left adjoint given by the base change functor $F:=B\otimes_{A}-$. We want to show that $G$ is fully faithful, for proving this it suffices to show that the counit 
\[
FG\to \mathrm{id}
\]
is an equivalence. More precisely, we need to show that for all $M\in \mathrm{LMod}_{B}(\n{C})$ the map 
\begin{equation}\label{eqorijqjwee}
B\otimes_{A} M\to M
\end{equation}
is an equivalence. But since $\n{C}$ is presentable, the tensor product is associative and the map \eqref{eqorijqjwee} is equivalent to
\[
B\otimes_A M = B\otimes_{A} (B\otimes_B M) = (B\otimes_A B)\otimes_B M \xrightarrow{m\otimes 1} B\otimes_B M=M.
\]
But by (1) the map $m$ is an equivalence proving what we wanted. 

Conversely, suppose that (4) holds, then the counit $FG\to \mathrm{id}$ is an equivalence, applying this to $B$ we get that $B\otimes_A B\xrightarrow{m} B$ is an equivalence which is precisely (1). 
\end{proof}

\begin{lemma}\label{lemmaIdempotentBaseChange}
    Let  $\n{C}$ be a presentable symmetric monoidal stable $\infty$-category, let $f:A\to B$ be a morphism of associative algebras in $\n{C}$. Consider a commutative square of associative algebras in $\n{C}$
    \[
    \begin{tikzcd}
        A \ar[r] \ar[d] & B \ar[d] \\
        C \ar[r] & D.
    \end{tikzcd}
    \]
    Suppose that $A\to B$ is idempotent and that the map $B\otimes_A C \to D$ is an isomorphism. Then $C\to D$ is idempotent. 
\end{lemma}
\begin{proof}
  The diagram of algebras gives rise to a commutative square of left modules
  \[
  \begin{tikzcd}
\mathrm{LMod}_{A}(\n{C}) \ar[r, "f"] \ar[d,"h"] & \mathrm{LMod}_B(\n{C}) \ar[d,"h'"] \\ 
\mathrm{LMod}_C(\n{C}) \ar[r,"f'"] & \mathrm{LMod}_{D}(\n{C}).
  \end{tikzcd}
  \]
Taking right adjoint of the vertical maps $h^{R}$ and $h'^{R}$ respectively, we get a natural transformation of functors $\mathrm{LMod}_{C}(\n{C}) \to \mathrm{LMod}_{B}(\n{C})$
  \[
  \alpha: fh^R\to  h'^Rf'
  \]
The map $\alpha$ is an isomorphism since we have an equivalence of $(B,C)$-bimodules 
\[
 B\otimes_A C \xrightarrow{\sim} D.
\]
Applying $\alpha$ to $D$ one sees that the map 
\[
  B \otimes_{A}D \xrightarrow{\sim} D\otimes_C D
\]
is an isomorphism of $(B,D)$-bimodules. Since $D$ is a $B$-module then the map $D\xrightarrow{1\otimes  \mathrm{id}} B\otimes_A D$ is an isomorphism, then property (3) of Lemma \ref{lemmaIdempotentAlgebras} holds and $D$ is an idempotent algebra over $C$. 
\end{proof}

\subsection{Hopf algebras}

Finally, we address the following proposition that will be used in different parts of the paper.

\begin{prop}
\label{PropositionHopfAlgebras}
Let $\n{R}$ be a static commutative analytic ring such that $-\otimes^L_{\n{R}}-$ is the left derived functor of $-\otimes_{\n{R}}-$.  Let $\n{A}$ be a static  $\n{R}$-Hopf algebra over $\n{R}$ with the induced analytic structure. Suppose that $\n{A}$ is cocommutative and that its antipode is an anti-involution, i.e. $s^2=\id$. Suppose that the self tensor products of analytic rings $\n{A}^{\otimes_{\n{R}} n}$ are static for all $n\in \bb{N}$. Then the following assertions hold: 

\begin{enumerate}
    \item The tensor product $-\otimes_{\n{R}}^L-$ defines a  symmetric monoidal structure on $\LMod_{\n{R}}(\n{A})$ obtained by restriction of scalars along the comultiplication $\Delta: \n{A} \to \n{A}\otimes_\n{R} \n{A}$.

    \item  ($\otimes$-$R \iHom$ adjunction) The derived internal $\iHom$ over $\n{R}$ induces a natural functor 
    \[
    R\iHom_{\n{R}}(-,-)_{\star_{1,3}}: \LMod_{\n{R}}(\n{A}) \times \LMod_{\n{R}}(\n{A}) \to \LMod_{\n{R}}(\n{A})
    \]
  given by precomposing the natural $\n{A}^{op}\otimes_{\n{R}} \n{A}$-module structure with the map $\n{A} \xrightarrow{\Delta} \n{A} \otimes_{\n{R}} \n{A} \xrightarrow{s\otimes 1} \n{A}^{op}\otimes_{\n{R}} \n{A}$,  where $s:\n{A} \xrightarrow{\sim } \n{A}^{op}$ is the antipode.  Furthermore, $R\iHom_{\n{R}}(-,-)_{\star_{1,3}}$ is a right adjoint of the internal tensor product $-\otimes_{\n{R}}^L-$. 

    \item (Twisting/untwisting) There are natural equivalences
    \[
    \begin{gathered}
     \Psi: \n{A} \otimes_{\n{R}}^L - \xrightarrow{\sim} \n{A}\otimes^L_{\n{R}} (-)_0   \\  \Phi: R\iHom_{\n{R}}(\n{A}, -)_{\star_{1,3}} \xrightarrow{\sim } R\iHom_{\n{R}}(\n{A}, -)_{\star_1} := R\iHom_{\n{R}}(\n{A}, (-)_{0}),
     \end{gathered}
    \]
    of endofunctors of $\LMod_{\n{R}}(\n{A})$, where $(-)_0$ is the trivial $\n{A}$-module structure obtained by restricting scalars along the composition $\n{A}\xrightarrow{\nu} \n{R} \xrightarrow{\mu} \n{A}$.

    \item Let $\iota: \LMod_{\n{R}}(\n{A}) \xrightarrow{\sim} \RMod_{\n{R}}(\n{A})$ be the precomposition with the antipode of $\n{A}$. We have  natural equivalences of functors 
    \begin{gather*}
    \iota(N)\otimes^L_{\n{A}} M = \n{R} \otimes^L_{\n{A}}(N\otimes^L_{\n{R}} M)  \\ 
    R\iHom_{\n{A}}(N,M) = R\iHom_{\n{A}}(\n{R}, R\iHom_{\n{R}}(N,M)_{\star_{1,3}})
    \end{gather*}
    for any $N,M\in \LMod_{\n{R}}(\n{A})$, where $\n{R}$ is equipped with an $\n{A}$-module structure through the counit.

    \item  Let $\n{B}$ be a static $\n{R}$-Hopf algebra satisfying the same hypothesis as $\n{A}$ and let $\n{A}\to \n{B}$ be a morphism of $\n{R}$-Hopf algebras. Then $\n{B}$ is an idempotent $\n{A}$-algebra if and only if $\n{B}\otimes_{\n{A}}^L \n{R} = \n{R}$. 
    
\end{enumerate}

\end{prop}

\begin{proof}
\begin{enumerate}
 \item First, let $\n{C}= \Mod_{\n{R}}$  be the symmetric monoidal $\infty$-category of $\n{R}$-modules, and let $\n{C}^{op}$ be its opposite category. Then, $\n{A}$ defines a commutative  Hopf algebra in the symmetric monoidal category $\n{C}^{op}$. Therefore,  the category $\mathrm{CoMod}_{\n{A}}(\n{C}^{op})= \varprojlim_{[n]\in \Delta} \n{A}^{\otimes_{\n{R}} n}\mbox{-}\Mod(\n{C}^{op})$ of (left) comodules over $\n{A}$ in $\n{C}^{op}$ is symmetric monoidal, with symmetric monoidal structure given by $-\otimes^L_{\n{R}}-$ on underlying objects. Part (1) follows since $\LMod_{\n{R}}(\n{A})= (\mathrm{CoMod}_{\n{A}}(\n{C}^{op}))^{op}$, and since the opposite of a symmetric monoidal category is symmetric monoidal.

 \item Given $N,M\in \LMod_{\n{R}}(\n{A})$,  we see $R\iHom_{\n{R}}(M,N)$ as an $\n{A}$-module via the forgetful functor through the algebra homomorphism   $A\xrightarrow{\Delta} A\otimes_{\n{R}} A \xrightarrow{s\otimes 1} A^{op}\otimes_{\n{R}} A$. To prove the $\otimes$-$R\iHom$ adjunction, since both functors arise as derived functors of suitable abelian categories with enough projectives and injectives (after fixing the cardinal $\kappa$), it suffices to  know the non-derived $\otimes$-$\iHom$ adjunction  of the underlying abelian categories, which is  \cite[Example 1.2.2 (3)]{TannakaDuality}.

 \item Let $\n{C}=\Mod_{\n{R}}$. We have an equivalence of symmetric monoidal categories $\Mod_{\n{R}}(\n{A})= \mathrm{CoMod}_{\n{A}}(\n{C}^{op})^{op}$. Let $f^*: \mathrm{CoMod}_{\n{A}}(\n{C}^{op}) \to   \n{C}^{op}$ be the forgetful functor taking the underlying object in $\n{C}^{op}$, and let $f_*$ be its right adjoint. In the opposite category $f^{*, op}$ is the forgetful from $\n{A}$ to $\n{R}$-modules, and $f_{*}^{op}= \n{A} \otimes_{\n{R}}^L-$. The functor $f^*$ is symmetric monoidal, we then have a natural transformation 
 \[
  f_* \n{R} \otimes M \to f_*f^*M
 \]
for $M\in \mathrm{CoMod}_{\n{A}}(\n{C}^{op})$. In the opposite category this translates to a natural transformation 
\[
\n{A} \otimes_{\n{R}} M_0 \to \n{A} \otimes^L_{\n{R}} M. 
\]
We claim that it is an isomorphism. By writing $M$ as filtered  colimits of projective generators, and since $\n{A}[S]= \n{A}\otimes_{\n{R}} \n{R}[S]$, one is reduced to the case when $M=\n{A}$. Following the construction, the map of $\n{A}$-modules  $\n{A}\otimes_{\n{R}} \n{A}_0 \to  \n{A}\otimes_{\n{R}} \n{A}$ is adjoint to the map
\[
\n{A}_0\xrightarrow{1\otimes \mathrm{id}} \n{A} \otimes_{\n{R}} \n{A}. 
\]
An inverse of this map can be given explicitly by the composite 
\[
\n{A} \otimes_{\n{R}} \n{A} \xrightarrow{\Delta\otimes \mathrm{id}} \n{A} \otimes_{\n{R}} \n{A}  \otimes_{\n{R}} \n{A}  \xrightarrow{\mathrm{\id}\otimes s \otimes \mathrm{\id}} \n{A}  \otimes_{\n{R}} \n{A} \otimes_{\n{R}} \n{A} \xrightarrow{\id \otimes m}  \n{A} \otimes_{\n{R}} \n{A}_0, 
\]
where $m: \n{A} \otimes_{\n{R}} \n{A} \to \n{A}$ is the multiplication map. Finally, the untwisting map $\Phi$ for the internal Hom follows from adjunction and the untwisting map $\Psi$. 

\item The natural transformation for the tensor product is a consequence of the following   natural equivalences for $N,M,Y\in \Mod_{\n{R}}(\n{A})$.
\[
\begin{aligned}
    R\iHom_{\n{R}}(\n{R} \otimes^L_{\n{A}}(N\otimes_{\n{R}} M), Y) & = R\iHom_{\n{A}}(N\otimes_{\n{R}} M, Y) \\ 
    & = R\iHom_{\n{A}}(M, R\iHom_{\n{R}}(N, Y)_{\star_{1,3}}) \\ 
    & = R\iHom_{\n{A}}(M,  R\iHom_{\n{R}} (\iota(N), Y)) \\  
    & = R\iHom_{\n{R}}( \iota(N) \otimes_{\n{A}}^L M , Y).
\end{aligned}
\]
The natural equivalence for the internal Hom's follows by the adjunction of point (2).

\item Suppose that $\n{B}$ is an idempotent $\n{A}$-algebra. Then we have that 
\[
\n{B} \otimes_{\n{A}}^L \n{R} = \n{B} \otimes_{\n{A}}^L (\n{B} \otimes_{\n{B}}^L \n{R}) =  (\n{B} \otimes^L_{\n{A}} \n{B}) \otimes_{\n{B}}^L \n{R} = \n{B} \otimes^L_{\n{B}} \n{R} =\n{R}.
\]
Conversely, suppose that $\n{B} \otimes_{\n{A}}^L \n{R} = \n{R}$, then by (the version for right modules of) part (4)  we have 
\[
\begin{aligned}
\n{B} \otimes_{\n{A}}^L \n{B}  & = (\n{B} \otimes_{\n{R}} \iota(\n{B})) \otimes_{\n{A}}^L \n{R} \\ 
& = (\n{B}_0 \otimes_{\n{R}}  \iota(\n{B} )) \otimes_{\n{A}}^L \n{R} \\ 
& = \n{B}_0 \otimes_{\n{R}} (\n{B} \otimes_{\n{A}}^L \n{R}) \\
& = \n{B},
\end{aligned}
\]
in the third equality we used the antipode $s: \n{B}^{op} \xrightarrow{\sim} \n{B}$ to identify the right  and left actions of $\n{A}$ on $\n{B}$. An explicit diagram chase shows that the resulting map $\n{B}\otimes_{\n{A}}^L \n{B} \to \n{B}$ is the multiplication map, proving that $\n{B}$ is an idempotent $\n{A}$-algebra.
\end{enumerate}
\end{proof}

\bibliographystyle{alpha}
\bibliography{biblioSolAnII}

\end{document}